%% file: mainfile.tex
\documentstyle[12pt,diagram]{article}
\input mssymb
\newcommand{\ga}{\alpha}
\newcommand{\gb}{\beta}
\renewcommand{\gg}{\gamma}
\newcommand{\gd}{\delta}
\renewcommand{\ge}{\varepsilon}
\newcommand{\gz}{\zeta}
\newcommand{\gth}{\theta}
\newcommand{\gk}{\kappa}
\newcommand{\gl}{\lambda}
\newcommand{\gm}{\mu}
\newcommand{\gn}{\nu}

\newcommand{\gp}{\pi}
\newcommand{\gr}{\rho}
\newcommand{\gs}{\sigma}
\newcommand{\gt}{\tau}
\newcommand{\gf}{\varphi}
\newcommand{\go}{\omega}
\newcommand{\gG}{\Gamma}

\newcommand{\gP}{\Pi}
\newcommand{\gF}{{\Phi}}
\newcommand{\gS}{\Sigma}
\newcommand{\seq}[1]{{\langle #1 \rangle}}

\newcommand{\implies}{\Longrightarrow}

\newtheorem{definition}{Definition}[section]
\newtheorem{theorem}{Theorem}[section]
\newenvironment{proof}{\noindent{\bf Proof:}}
{\nopagebreak\mbox{}\newline\makebox[\textwidth]{\hfill
$\blacklozenge$}\par\bigskip}
\newtheorem{lemma}{Lemma}[section]
\newtheorem{corollary}{Corollary}[section]
\newtheorem{claim}{Claim}[section]

\catcode`\@=11
\def\@begintheorem#1#2{\rm \trivlist \item[\hskip \labelsep{\bf #1\ #2:}]}
\def\@opargbegintheorem#1#2#3{\rm \trivlist
      \item[\hskip \labelsep{\bf #1\ #2\ (#3):}]}
\catcode`\@=12

\newcommand{\ult}{{\hbox{\rm ult}}}
\newcommand{\re}{{\restriction}}
\newcommand{\br}{{\hbox{\bf r}}}
\newcommand{\brr}{{\hbox{${\br}^-$}}}
\newcommand{\chunk}[1]{{[\![#1]\!]}}
\newcommand{\hull}{{\hbox{\rm Hull}}}
\newcommand{\On}{\hbox{\rm Ord}}
\newcommand{\rt}{\hbox{\rm root}}

\begin{document}
\title{How to win some simple iteration games}
\author{Alessandro Andretta\thanks{Partially supported by
the Italian CNR and by 40\% Funds of MURST}\\
Universit\'a di Camerino, Italy
\and John Steel\thanks{Partially supported by NSF Grant
DMS92-06946}\\
University of California at Los Angeles}
\date{September 1993}
\maketitle

\centerline{\bf Abstract} 

\noindent
We introduce two new iteration games: the game $\cal G$,
which is a strengthening of the weak iteration game,
and the game ${\cal G}^+$, which is somewhat stronger than
$\cal G$ but weaker than the full iteration 
game of length $\go_1$.
For a countable $M$ elementarily embeddable in some
$V_{\eta}$, we can show that $II$ wins ${\cal G}(M,\go_1)$
and that $I$ does not win the ${\cal G}^+(M)$.

\section{Introduction}
\input introduction.tex

\section{Iteration games}
\input games

\section{Preliminaries}
\input preliminaries

\section{Strong past one Woodin cardinal}
\input strongpast

\section{How not to lose}
\input hownotlose

\section{The tree $\cal U$}
\input bigtree.tex

\section{The enlargement}
\input enlargement

\end{document}

%% file: introduction.tex
%

Iterability results, that is theorems ensuring the 
existence of wellfounded branches in iteration trees, 
are the main technical tool used in proving
the comparison theorem for inner models for
large cardinals.
The main iterability result of \cite{Iterationtrees}, 
Theorem 4.3, shows that any countable iteration tree 
$\cal T$ on a countable $M\preceq V_{\ga}$, has a 
maximal wellfounded branch, and this is enough 
to prove a comparison theorem for the canonical 
inner model for one Woodin cardinal.
In fact in \cite{Manywoodins} this result is used
to prove a comparison theorem for countable 
{\it tame\/} premice $J^{\vec{\cal E}}_{\ga}$, i.e. 
structures in the sense of \cite{FSIT} 
satisfying \lq\lq$\gd$ is not Woodin" for 
every $(\gk,\gl)$-extender $E$ on the 
coherent sequence $\vec{\cal E}$, with $\gk<\gd<\gl$. 
Tame premice can have many Woodins, but cannot 
satisfy the sentence: {\em there is $\gk$
which is $\gd+1$-strong and $\gd$ is Woodin.\/}
On the other hand the absence of more powerful 
iterability results has been the main obstacle 
towards extending the existing theory to core 
models with larger cardinals.
The Cofinal Branch Hypothesis (CBH) 
(for the definition of this or other notions 
see \cite{Iterationtrees} or \S2 below) 
is the single most important open question 
in this area, and a proof of it (if true) 
would almost certainly yield a comparison 
lemma for mice with, say, superstrong cardinals.
Barring CBH, the next best thing we could hope 
to prove is the Strategic Branch Hypothesis (SBH), 
which is a weakening of CBH. 
As the name suggests SBH asserts that player $II$
has a winning strategy in the {\em full iteration game\/} 
on $V$ of length $\nu$, ${\cal I\cal G}(V,\nu)$, 
for any $\nu$. In this game
the two players cooperatively build 
in $\gn$ rounds an iteration tree on $V$,
with $II$ on the move at limit rounds choosing 
a cofinal wellfounded branch.
Just as with CBH, SBH is pretty much open.
Theorem 4.3 of \cite{Iterationtrees} implies that 
player $II$ has a winning strategy in the
{\em weak iteration game\/} on countable 
$M\preceq V_{\ga}$, ${\cal W\cal G}(M)$.
This is a weaker game than ${\cal I\cal G}(M)$,
in the sense that if $II$ wins ${\cal I\cal G}(M)$
then he wins ${\cal W\cal G}(M)$.
On the other hand the weak iteration game seems
of little or no use in proving a comparison 
theorem for non-tame premice. 

In this paper we prove two new iterability results
which yield a comparison lemma for non tame mice.
The extent to which our results civilize these \lq\lq wild"
mice is not clear, but it should fall somewhere between
the hypotheses: {\em a strong cardinal below
a Woodin\/} and {\em a Woodin limit of Woodins\/}.

Our first result, proved in \S 4, says that player $II$ wins 
a certain game which we call
${\cal G}(M,\go_1+1)$, when $M\preceq V_{\ga}$ is 
countable. The game $\cal G$ is stronger than 
$\cal W\cal G$, but much weaker than $\cal I\cal G$.
This is just about the best
we are able to show in the line of proving directly
that player $II$ has a winning strategy for 
games approximating $\cal I\cal G$.

The second result, which takes up the rest of the paper
\S\S 5, 6 and 7, deals with an iteration game ${\cal G}^+$
which is a much closer approximation to $\cal I\cal G$. 
It is played like $\cal I\cal G$ except for the fact that 
$I$ has to play distinct integers on the side. 
The game is over once $I$ runs out of integers,
provided none of the players has lost by that time.
We are able to prove that $I$ does not have 
a winning strategy in ${\cal G}^+(M)$ for 
countable $M\preceq V_{\ga}$.
(So perhaps this paper could have been more aptly entitled:
How not to lose a short iteration game.)
By results of Steel and Woodin, ${\cal G}^+(M)$ is determined, 
modulo supercompact cardinals, hence $II$ wins the game.

We think that both proofs present interesting new features. 
In a way these are more important than the statements of
the theorems themselves.
Both results seem likely to admit further generalizations,
although at this time we do not know how to do it.
One drawback to our present approach is the use 
of $2^{\aleph_0}$-closed extenders in the proofs.
In fact in the proof of Theorem \ref{little}
we must also assume that the iteration trees 
are non overlapping. 
This is not too great a restriction 
if the goal is to construct an inner model
$L[\vec{\cal E}]$ with many Woodins assuming 
the existence of such cardinals in the universe,
as the extenders witnessing Woodiness in $V$
can be taken to be as closed as we want. 
Of course this would be a problem were we not to assume
the existence of large cardinals in $V$ in building
$L[\vec{\cal E}]$, as done in core model theory.

This paper is fairly self-contained, but the 
reader is assumed to be acquainted with iteration 
trees and extenders. Sections \S\S 1, 3 and parts of \S 5
of \cite{Iterationtrees} would do.
No knowledge of fine structure or inner model theory is required.

\bigskip \noindent
{\bf Aknowledgements.}
These results were obtained while the first 
author was visiting the UCLA Mathematics Department 
during the academic year 91-92 and the summer of 93. 
He would like to thank the Mathematics Department of 
UCLA for its hospitality. Both authors would like to 
thank  Tony Martin and Philip Welch for making helpful 
suggestions.

%% file: games.tex
%

The most general iteration game 
is the {\it full iteration game} and was defined 
in \S 5 of \cite{Iterationtrees}.
In the full iteration game of length $\nu$ on a premouse $M$, 
${\cal I \cal G}(M,\nu)$, players $I$ and $II$ cooperatively
construct a plus-2 normal iteration tree $({\cal T},M)$: $I$ plays 
at successor rounds, while $II$ plays at limit rounds. 
At round $\ga+1<\nu$, $I$ plays an extender 
$E_{\ga}\in M_{\ga}$ and an ordinal 
$\gr_{\ga}$ such that $M_{\ga}\models$\lq\lq$E_{\ga}$ 
is $\gr_{\ga}+2$ strong." 
Let $P=\ult(M_{\gb},E_{\ga})$, where $\gb$ is least 
such that crit$(E_{\ga})\leq\gr_{\gb}$.
If $P$ is illfounded then $I$ wins, otherwise let
$M_{\ga+1}=P$ and we move to the next round $\ga+1$.
At limit rounds $\gl<\nu$, $II$ plays a cofinal 
wellfounded branch $b$ of the iteration tree
built insofar, and set $M_{\gl}=M_b$.
(At round 0, neither player does anything.) The first player 
that cannot make a legal move loses. 
If neither player has lost by round $\nu$, then $II$ wins.
[The reader should keep in mind that, as we are dealing
with normal iteration trees, the game described above is 
slightly more restrictive than the game described in
\cite{Iterationtrees}.] 

The Strategic Branch Hypothesis (SBH) asserts that 
$V$ is strategically iterable, i.e. player
$II$ has a winning strategy for ${\cal I\cal G}(V,\gn)$,
for all $\nu$. 
It is a weaker form of the Cofinal Branch Hypothesis
(CBH), asserting that: if $\cal T$ is an iteration tree  
on $V$ then if $\cal T$ is of limit length 
it has a cofinal wellfounded branch, and
if $\cal T$ is of successor length, we do not run
into problems by taking an ultrapower and extending
the tree one more step. 
Note that SBH is preserved by going to elementary substructures: 
if $M$ is a countable premouse elementarily embeddable 
in some $V_{\ga}$ via $\pi:M\to V_{\ga}$ 
and $\gS$ is a strategy for $II$
in ${\cal I\cal G}(V,\gn)$, then a strategy for
$II$ in ${\cal I\cal G}(M,\nu)$ is obtained
by copying via $\pi$ and following $\gS$.

The argument above does not apply, though, to CBH.
Theorem 4.3 of \cite{Iterationtrees} shows 
that every countable tree on a countable $M\preceq V_{\ga}$ 
has a maximal wellfounded branch. 
On the other hand the analogous statement
on $V$ is open even for trees of height $\go$.

\medskip\noindent
{\bf Open problem 1:}
Does every countable iteration tree $\cal T$ on $V$ 
have a maximal wellfounded branch? In particular: 
does every iteration tree of height $\go$ have a (necessarily
cofinal) wellfounded branch?

\medskip\noindent
That the answer is affirmative 
for trees $\cal T$ where all extenders are 
$2^{\aleph_0}$-closed in the model they appear, 
is the content of 

\medskip
\noindent {\bf Theorem 5.6 of \cite{Iterationtrees}:} 
{\em Suppose $\cal T$ is a countable iteration 
tree on a premouse $N$, ${}^{2^{\aleph_0}}N\subseteq N$ 
and for all $\ga+1<lh({\cal T})$, 
$M_{\ga}^{\cal T}\models\lq\lq\mbox{\em ult}(V,E_{\ga})$ is
$2^{\aleph_0}$-closed," then there is a maximal 
wellfounded branch $b$ of $\cal T$.} 

\medskip

This result will be used in \S 4 of this paper.
A related conjecture is the Unique Branch Hypothesis
(UBH) asserting that every iteration tree on $V$ 
of limit length has at most one cofinal 
wellfounded branch.
Woodin, in unpublished work, has shown the
consistency of $\neg(\mbox{UBH}+\mbox{CBH})$ 
assuming the existence of a non-trivial
$j:L(V_{\gl+1})\to L(V_{\gl+1})$, 
where $\gl=\sup_n j^n(\gk)$ and $\gk=\mbox{crit}(j)$.
Thus it is quite possible that CBH is consistently false,
although at this point we have no reason to believe either way.

The {\it weak iteration game\/} of length $\nu$ on 
a premouse $M$, ${\cal W\cal G}(M,\nu)$, is a weakening
of the full iteration game, with player
$I$ playing only at successor rounds
and player $II$ playing at every round. 
At round $\ga<\nu$, $\seq{({\cal T}_{\gb},P_{\gb})\mid\gb<\ga}$
and $\seq{j_{\gb,\gg}\mid \gb<\ga}$ are given such that

\begin{enumerate}

\item
$P_0=M$, each $({\cal T}_{\gb},P_{\gb})$ 
is an iteration tree of successor length $\gth_{\gb}+1<\go_1$,
$j_{\gb,\gb+1}=i_{0,\gth_{\gb}}^{{\cal T}_{\gb}}$ and
$j_{\gg,\gb+1}=j_{\gb,\gb+1}\circ j_{\gg,\gb}$;

\item
$({\cal T}_{\gb+1},P_{\gb+1})\|({\cal T}_{\gb},P_{\gb})$;
that means $P_{\gb+1}=M_{\gth_{\gb}}^{{\cal T}_{\gb}}$ 
is the last model of ${\cal T}_{\gb}$,
$\gr_0^{{\cal T}_{\gb+1}}\geq\sup
\{\gr_{\gg}^{{\cal T}_{\gb}}\mid\gg+1\leq\gth_{\gb}\}$
and the first model $E_0^{{\cal T}_{\gb+1}}$
can be applied to is $P_{\gb+1}=M_0^{{\cal T}_{\gb+1}}$;

\item
if $\gg<\ga$ is limit, then $P_{\gg}$ is the 
direct limit of the $P_{\gb}$'s and $j_{\gb,\gg}$ are the
limit maps, for $\gb<\gg$.

\end{enumerate}
So $\seq{({\cal T}_{\gb},P_{\gb})\mid\gb<\ga}$
forms an iteration tree $({\cal T},M)$, with
$({\cal T}_{\gg}, P_{\gg})$ stacked on top of
$({\cal T}_{\gb}, P_{\gb})$ for $\gb<\gg<\ga$.
If $\ga$ is limit, then $II$ is to move and
must play a cofinal wellfounded branch of the
tree constructed so far. There is not much 
choice in this case as there is only one cofinal
branch of $\cal T$: if the direct limit 
of the $P_{\gb}$'s is illfounded
then $II$ loses, otherwise that will be $P_{\ga}$.
If $\ga$ is successor, $\ga=\gb+1$, then $I$ 
plays a putative iteration tree $({\cal S}_{\ga}, P_{\ga})$,
with $P_{\ga}=M_{\gth_{\gb}}^{{\cal T}_{\gb}}$ such that, 
extending ${\cal T}$ via ${\cal S}_{\ga}$,
we still have a putative iteration tree on $M$.
[A putative iteration tree is
an object obeying all the usual rules for
ordinary iteration trees except for the fact 
that the last model can be illfounded.] 
$II$ responds by playing either:

\begin{enumerate}

\item
({\bf accept}), if ${\cal S}_{\ga}$ is of successor 
length and its last model is wellfounded, that is:
${\cal S}_{\ga}$ really is an iteration tree on $P_{\ga}$; 
then set ${\cal T}_{\ga}={\cal S}_{\ga}$. Or

\item
$(\mbox{\bf accept},b)$, if ${\cal S}_{\ga}$ 
is of limit length and $b$ is a cofinal 
wellfounded branch; then let ${\cal T}_{\ga}$ 
be ${\cal S}_{\ga}$ extended via $b$ and
$\gth_{\ga}=lh({\cal S}_{\ga})$. Or

\item $(\mbox{\bf reject},b)$, where $b$ is a maximal
wellfounded branch of ${\cal S}_{\ga}$; let 
${\gth}_{\ga}=\sup(b)$, and ${\cal T}_{\ga}$ be 
${\cal S}_{\ga}\re\gth_{\ga}$ extended via $b$.

\end{enumerate}
It is easy to see that if $II$ wins 
${\cal I \cal G}(M,\gn)$ for any countable $\nu$, then $II$
wins ${\cal W\cal G}(M,\go_1+1)$.
Let us recall the main result of \S4 of \cite{Iterationtrees}.

\medskip

\noindent {\bf Theorem 4.3 of \cite{Iterationtrees}:}
{\it If $N$ is a countable premouse,
$\gf:N\to V_{\eta}$ is elementary and $\cal T$ is 
a countable iteration tree on $N$, then there exists
a maximal branch $b$ and an elementary embedding 
$\gt$ such that $\gt\circ i_{0,b}^{\cal T}=\gf$.\/}

\medskip

This easily implies that $II$ 
wins ${\cal W\cal G}(M,\go_1+1)$, for
$M$  countable and embeddable in some $V_{\eta}$
via $\pi:M\to V_{\eta}$:
it is enough to maintain inductively that 
at round $\ga<\nu$ we have elementary embeddings
$\gs_{\gg}:P_{\gg}\to V_{\eta}$, for $\gg\leq\ga$,
$\gs_0=\pi$, so that, for $\gb<\gg\leq\ga$, the diagram
$$
\begin{diagram}
\node{P_{\gg}}
\arrow[2]{se,t}{\gs_{\gg}}\\
\node{P_{\gb}}
\arrow{n,l}{j_{\gb,\gg}}\arrow{ese,t}{\gs_{\gb}}\\
\node{M=P_0}\arrow{n,l}{j_{0,\gb}}\arrow[2]{e,b}{\pi}
\node[2]{V_{\eta}}
\end{diagram}
$$
commutes. If $\ga$ is limit then $\gs_{\ga}$ is the limit map,
and if $\ga=\gb+1$ then $\gs_{\ga}:P_{\ga}=
M_{\gth_{\gb}}^{{\cal T}_{\gb}}\to V_{\eta}$ 
is the map $\gt$ given by Theorem 4.3 when
$N= P_{\gb}$, $\gf=\gs_{\gb}$ and 
${\cal T}={\cal T}_{\gb}$.

The weak iteration game described above will not suffice
to ensure that the comparison process for non-tame mice
terminates. The reason is that at some round 
$\ga+1<\go_1$ we might be forced by the 
comparison process to apply an extender
$E\in M_{\gth_{\ga}}^{{\cal T}_{\ga}}$
to an earlier model $M_{\gg}^{{\cal T}_{\gb}}$,
$\gg\leq\gth_{\gb}$ and $\gb<\ga$. Suppose
$P_{\ga+1}=\ult(M_{\gg}^{{\cal T}_{\gb}}, E)$
and ${\cal S}_{\ga+1}$  is an iteration tree
on $P_{\ga+1}$ (rather than letting $P_{\ga+1}=
M_{\gth_{\ga}}^{{\cal T}_{\ga}}$, as in $\cal W\cal G$). 
When this happens we write
$({\cal T}_{\ga+1}, P_{\ga+1})\bot({\cal T}_{\ga},P_{\ga})$.
The technique used before, i.e. embedding the $P$'s back 
to $V_{\eta}$ does not apply here, because 
the embeddings $\gs$'s do not agree enough with 
one another to ensure that $P_{\ga+1}$ embeds back in $V_{\eta}$.

In \S 4 we introduce a new game ${\cal G}(M,\nu)$ 
in which $I$ is allowed to go back and construct 
$({\cal T}_{\ga+1},P_{\ga+1})\bot({\cal T}_{\ga},P_{\ga})$
infinitely often, and we will show that $II$ wins 
${\cal G}(M,\go_1+1)$, for countable $M$ embeddable in some 
$V_{\eta}$ (see Theorem \ref{little}).
In order to highlight the ideas in that proof, we
briefly describe the techniques needed to prove 
a simpler result.

Assume, as usual, that $M$ is countable and 
embeddable in $V_{\eta}$ via $\pi$.
Suppose that the game considered is just
like ${\cal W\cal G}$ except that player $I$ at any stage
$\ga+1$ may play $({\cal S}_{\ga+1}, P_{\ga+1},
E,\gb,\gg)$, where $E$ is an extender in 
$M_{\gth_{\ga}}^{{\cal T}_{\ga}}$, $\gb\leq\ga$
and $\gg<\min (\gth_{\gb}+1,\gth_{\ga})$,
$P_{\ga+1}=\ult(M_{\gg}^{{\cal T}_{\gb}},E)$
and ${\cal S}_{\ga+1}$ is an iteration tree on $P_{\ga+1}$.
But from this point on the game proceeds as in the weak game. 
In other words: we can go back, if we want, but only once.
The trick is to introduce an intermediate model 
$N$ between $M$ and $V_{\eta}$, so that 
$V_{\eta}$ is the background universe of $N$, 
and $N$ is the background universe of $M$.
As long as we play the weak iteration game 
we just copy the trees on $N$ and then choose 
the branches  by playing the weak game on $N$.
If we do go back at some stage $\ga+1$ and
take $P_{\ga+1}=\ult(M_{\gg}^{{\cal T}_{\gb}}, E)$,
we use the copy construction between $M$ and $N$ to 
embed $P_{\ga+1}$ back into $N$, 
and hence into $V_{\eta}$. From this point on we 
simply play the weak game on $M$.

Formally, let $\gk>\eta$ 
and let $\gs_0:N\to V_{\gk}$, where $N$ is
of size $2^{\aleph_0}$ and contains all reals.
Suppose also $\pi_0:M\to N\cap V_{\bar{\eta}}$,
some $\bar{\eta}\in N$, is such that
$\pi=\gs_0\circ\pi_0$. 

Let's make it as a rule that the extender played
are $2^{\aleph_0}$-closed.
We now start playing the game.
Suppose that until round $\ga+1$ the weak 
iteration game was played, so that 
$(\forall\gb+1\leq\ga){\cal T}_{\gb+1}\|{\cal T}_{\gb}$,
i.e. ${\cal T}_{\gb+1}$ is an iteration tree on $P_{\gb+1}$
which is the last model $M_{\gth_{\gb}}$
of the tree ${{\cal T}_{\gb}}$.
Suppose also that the concatenation of the ${\cal T}_{\gb}$'s 
can be copied via $\pi_0$ on $N$ and let
$$
\pi_{\gb}:P_{\gb}\to Q_{\gb}=^{def} M_0^{\pi_{\gb}{\cal T}_{\gb}}
$$
be the copy map. Suppose also that we are given embeddings
$\gs_{\gb}:Q_{\gb}\to V_{\gk}$ such that 
$\gs_{\gg}=\gs_{\gb}\circ j_{\gb,\gg}^N$, 
where $j_{\gb,\gg}^N:Q_{\gb}\to Q_{\gg}$ 
are the embeddings given by the copied trees.

If $I$ plays $({\cal S}_{\ga}, P_{\ga})$, i.e. if he keeps 
on playing the weak game, then we choose $\gth$ largest such that 
${\cal T}_{\ga+1}={\cal S}_{\ga+1}\re\gth$ 
can be copied on $V_{\gk}$ via $\gs_{\ga+1}\circ\pi_{\ga+1}$
and has no non-cofinal wellfounded branches.
Theorem 5.6 of \cite{Iterationtrees} guarantees 
the existence of a cofinal wellfounded branch.
Let $P_{\ga+2}$ and $Q_{\ga+2}$ be the models 
$M_b^{{\cal T}_{\ga+1}}$ and $M_b^{\pi_{\ga+1}{\cal T}_{\ga+1}}$ 
respectively. 
A tree argument enables us to replace the
copy map from $Q_{\ga+2}$ to (a rank of) 
$M_b^{\gs_{\ga+1}\pi_{\ga+1}{\cal T}_{\ga+1}}$ 
with a similar embedding 
belonging to the latter model.
So by elementarity we get 
$\gs_{\ga+2}:Q_{\ga+2}\to V_{\gk}$.

If, otherwise, $I$ decides to go back and
take $\ult(M_{\gg}^{{\cal T}_{\gb}},E)$
then a tree argument is used to replace the 
copy map from $P_{\ga+1}$ to $Q_{\ga+1}=\ult
(M_{\gg}^{\pi_{\gb}{\cal T}_{\gb}},\pi_{\ga}(E))$
with a similar map that belongs to $Q_{\ga+1}$
and then we pull it back to $N$.
In order to do this we need to know that 
$Q_{\ga+1}$ is wellfounded and that 
$P_{\ga+1}$ belongs to it.
If we assume, as we do, that the iteration 
trees are non-overlapping, then $Q_{\ga+1}$ is
wellfounded by Theorem 1.2 of \cite{Wfdd}
or Lemma \ref{successor}.
For $P_{\ga+1}\in Q_{\ga+1}$ we seem to need 
that $N$ contains $HC$. This on the other hands
forces $N$ to be uncountable and so the usual
tree argument will not apply to it.
To overcome this difficulty, we have to resort 
to the concept of {\em support} (cf. Definition \ref{definitionsupport})
and $2^{\aleph_0}$-closed
extenders.

To summarize: the copy maps $\pi_{\ga}$ are needed 
in order to be able to embed the ultrapower
$\ult(M_{\gg}^{{\cal T}_{\gb}},E)$ back into the
$V$-like model $N$, while the maps $\gs$'s are
needed to ensure that the direct limit models
of the tree copied on $N$ are wellfounded.
The $\gs$'s would be superfluous, were we able to
prove the following instance of (CBH).

\medskip\noindent
{\bf Open problem 2:}
Consider the following game. (For notational simplicity we
state the length $\go$ case only.)
$I$ plays iteration trees ${\cal T}_n$
and $II$ plays cofinal wellfounded branches
$b_n$ such that:
\begin{enumerate}

\item
${\cal T}_0$ is on $V$ and 
${\cal T}_{n+1}$ is on $M_{b_n}^{{\cal T}_n}$ and

\item
each ${\cal T}_n$ has no non-cofinal wellfounded branch
and all extenders used in ${\cal T}_n$ are 
$2^{\aleph_0}$-closed in the model they appear.
\end{enumerate}
The first player to violate the rules loses. If 
neither player has lost by the end of the game, then
$II$ wins iff the (only) cofinal branch of the 
resulting tree is wellfounded.

\noindent
Does $II$ win this game?

\medskip
By Theorem 5.6 of \cite{Iterationtrees} stated
above, $II$ does not 
lose at any finite position of the game, and
with some extra work it can be shown that $I$ 
does not have a winning strategy.

Although the game ${\cal G}(M,\go_1+1)$
of \S 4 ensures enough iterability to
prove a comparison theorem for inner models 
with a cardinal strong past one Woodin cardinal
(and slightly beyond), there seems to be 
genuine difficulties in generalizing ${\cal G}$ 
to handle stronger hypotheses.
As we always deal with countable iteration trees
on countable premice (hence objects that can be
coded as reals),
the various iteration games can be studied from the
point of view of descriptive set theory.
In particular, rather than trying to  prove outright that
$II$ has a winning strategy in a given iteration game, 
one can try to show that $I$ does not have a winning 
strategy and then appeal to determinacy.
Although real games (i.e. games in which the 
players play elements of ${}^{\go}\go$) of length
$\go_1$ are not determined, by work of Steel and Woodin
it is consistent that {\em variable length games\/}
of reasonable complexity are determined, assuming 
large cardinals. 
The expression \lq\lq variable length" means that
the length of the game varies with the play:
for example we can stipulate that the game
is over when we reach a position $p$ of length $\nu$, 
where $\nu$ is the least admissible in $p$ larger than
$\go$. A stronger game is obtained by letting
$\nu$ be the second admissible in $p$.
Another family of long games are the 
{\em continuously coded\/} ones:
at stage $\ga$, $I$ plays a real $x_{\ga}$ {\em and\/} a 
natural number $n_{\ga}$ such that 
$n_{\ga}\notin\{n_{\gb}\mid\gb<\ga\}$ and
the game is over when $I$ runs out of integers.
Continuously coded games are stronger than games ending
at the first admissible, but weaker than the 
ones ending at the second admissible in the play.

Steel and Woodin proved that if there is a 
supercompact cardinal, then it is consistent that 
{\em all continuously coded closed-$\Pi^1_1$ real games
are determined\/} (see \cite{Longgames} for a proof of this
and other basic facts about long games).
In \S 5 a new iteration game  ${\cal G}^+(M)$ is introduced.
It is a continuously coded closed-$\Pi^1_1$ real game. 
In \S 6 and \S 7 it is shown (Theorem \ref{main}) that $I$ does
not win ${\cal G}^+(M)$ for countable $M$ elementarily embeddable
in some $V_{\eta}$. Hence, modulo supercompact cardinals,
$II$ wins ${\cal G}^+(M)$.

The game ${\cal G}^+(M)$ should yield enough iterability
to give a comparison theorem for inner models with 
many strong cardinals overlapping Woodin cardinals, 
but it is still too weak for hypotheses like a 
Woodin limit of Woodins. In order to get a comparison 
theorem for inner models with large cardinals that powerful, 
we believe that progress must be made in two distinct areas.
For one, Theorem \ref{main} must be strengthened
to non-continuously coded games: unfortunately
our proof seems to use continuity in an essential way.
The second area that needs to be further developed
is more descriptive set theoretic in nature, as
we need more powerful and sharper results
concerning the determinacy of long games.

%% file: preliminaries.tex
%

In this section we define {\em pseudo-iteration trees\/},
which are a generalization of 
iteration trees (\cite{Iterationtrees}, \cite{Projdeter}).
Besides of being of independent interest, 
pseudo-iteration trees will be a key ingredient 
in the main part of the present paper \S 5, \S 6 and \S 7,
where an iterability result about ordinary iteration
trees is proved. 
Several basic facts about iteration trees hold also 
in this more general set-up, so we preferred to give 
a unified treatment to the subject, rather than repeating 
the arguments twice, first for ordinary iteration trees 
and then for their \lq\lq pseudo" siblings. 
Pseudo-iteration trees will make no appearance until
\S 5 so the reader only interested in \S 4 may skip
some of the material in the present section.
The reader should keep in mind, though, that the
notions of {\em support\/} and {\em  chunk\/}, and 
in particular Lemma \ref{chunks} will be used in \S 4. 

By a {\em coarse premouse\/}, or simply a premouse,
we mean a transitive set or class $M$ with a distinguished 
ordinal $\gd=\gd(M)\in M$ such that $M$ is power admissible,
satisfies choice, comprehension and the collection schema
for domains $\subseteq V_{\gd}$. Whenever a 
$(\gk,\gl)$-extender $E$ is applied to a premouse 
$M$, it will always be assumed that $\gk<\gd(M)$, 
so that \L os' theorem holds for $\ult(M,E)$
and the embedding $i_E^M$ is fully elementary.
An ordinal $\gg$, $\gd(M)<\gg<M\cap\On$
is a {\em cut-off point\/} of $M$ iff 
$M\cap V_{\gg}$ is still a premouse with
$\gd(M\cap V_{\gg})=\gd (M)$.
We say that two transitive sets or classes
$M$ and $N$ agree through an ordinal $\gr$
iff $M\cap V_{\gr}=N\cap V_{\gr}$.
A {\em tree ordering on\/ $\gth$ with\/ $\gl +1$
roots\/}, $\gl<\gth$, is a transitive,
irreflexive, wellfounded relation $<_T$ on $\gth$
such that

\begin{enumerate}

\item 
$\forall\ga,\gb <\gth (\ga<_T \gb\implies\ga<\gb)$ 
and for all $\gb<\gth$ the set 
$\{\ga<\gth\mid\ga<_T\gb\}$ is linearly
ordered by $<_T$.

\item 
$\forall\ga,\gb\leq\gl (\ga\neq\gb\implies
\ga,\gb\hbox{ are $<_T$-incomparable })$ and
$\forall\gb(\gl<\gb<\gth\implies\exists !\ga\leq\gl
(\ga<_T\gb))$.

The ordinals $\leq\gl$ are called roots and
$\rt_T(\gb)$ is the unique $\ga\leq\gl$
such that $\ga<_T\gb$ or $\ga=\gb$.

\item
$\forall\ga(\gl<\ga<\gth)$

$\ga$ is a successor $\iff$ $\ga$ is a $<_T$-successor,

$\ga$ is a limit $\implies$ $\{\gg\mid\gg <_T\ga\}$ is
cofinal in $\ga$.

\end{enumerate}
$\ga\leq_T\gb$ stands for $\ga<_T\gb\vee\ga=\gb$, 
and $[\ga,\gb]_T=
\{\gg\mid\ga\leq_T\gg\leq_T\gb\}$.
Similarly we define $[\ga,\gb)_T$, $(\ga,\gb)_T$, etc.
If $b$ is a {\em branch\/}, i.e. a
maximal $<_T$-linearly ordered subset of $\gth$,
$\rt_T(b)$ is the least $\ga\in b$.
If $\ga+1>\gl$, then $<_T$-pred$(\ga+1)$ 
is the least $\gb$ such that $\gb<_T\ga+1$.

\begin{definition}
A plus-$n$ pseudo-iteration tree of length $(\gth,\gl)$, 
with $\gl<\gth$, is a pair $(\cal T, \cal B)$ where 
\begin{enumerate}

\item
${\cal B}=\seq{B_{\ga}\mid\ga\leq\gl}$ 
is a sequence of premice,  
called base models, together with
a sequence of increasing ordinals 
$\gr_{\ga}$, for $\ga<\gl$ such that 
$B_{\ga}$ and $B_{\gb}$ agree 
through $\gr_{\ga}+n$, that is 
$B_{\ga}\cap V_{\gr_{\ga}+n}=B_{\gb}\cap V_{\gr_{\ga}+n}$, 
for $\ga<\gb\leq\gl$;

\item
$\cal T$ is a tree ordering $<_T$ on $\gth$ with
$\gl+1$ roots, together with a sequence
$$
\seq{(E_{\ga},\gr_{\ga})\mid \gl<\ga +1<\gth}
$$
of extenders and ordinals obeying the usual restrictions
for iteration trees, that is: there are premice 
$M_{\ga}^{(\cal T,\cal B)}$ and elementary embeddings
$i_{\ga,\gb}^{(\cal T,\cal B)}:
M_{\ga}^{(\cal T,\cal B)}\to M_{\gb}^{(\cal T,\cal B)}$,
$\gd(M_{\gb}^{({\cal T}, {\cal B})})=
i_{\ga,\gb}^{({\cal T},{\cal B})}(\gd(M_{\ga}^{({\cal T}, {\cal B})}))$,
for $\ga<_T\gb$, and such that 
\begin{enumerate}

      \item
      the sequence $\seq{\gr_{\ga}\mid\ga+1<\gth}$
      is increasing and;

      \item 
      $M_{\ga}^{(\cal T,\cal B)}=B_{\ga}$, for $\ga\leq\gl$;

      \item
      if $\gl<\ga+1<\gth$, then $M_{\ga}^{(\cal T,\cal B)}
      \models\lq\lq E_{\ga}$ is an extender 
      $\gr_{\ga}+n$ strong", $E_{\ga}\in V_{\gd(M_{\ga})}$,
      and letting $\gb=<_T \mbox{-pred}(\ga+1)$, then $\gb$ is least 
      such that $\gr_{\gb}+n>\mbox{crit}(E_{\ga})$,
      $$
      M_{\ga+1}^{(\cal T,\cal B)}=\ult 
      (M_{\gb}^{(\cal T,\cal B)},E_{\ga}),
      $$
      $i_{\gb,\ga+1}^{(\cal T,\cal B)}$
      is the canonical ultrapower embedding 
      $i_{E_{\ga}}^{M_{\gb}}$, and 
      $i_{\gb,\ga+1}^{({\cal T,\cal B})} 
      \circ i_{\gg,\gb}^{({\cal T,\cal B})}=
      i_{\gg,\ga+1}^{({\cal T,\cal B})}$, for $\gg<_T\gb<_T\ga+1$;

      \item
      if $\gl<\ga<\gth$ is limit, then $M_{\ga}^{(\cal T,\cal B)}$
      is the direct limit of $M_{\gb}^{(\cal T,\cal B)}$
      for $\gb<_T\ga$ and the $i_{\gb,\ga}^{({\cal T,\cal B})}$
      are the direct limit maps.

\end{enumerate}
\end{enumerate}
\end{definition}

\noindent {\bf Remarks.}
\begin{enumerate}

\item 
For $\ga<\gb<\gth$, 
$M_{\ga}^{(\cal T,\cal B)}$ and 
$M_{\gb}^{(\cal T,\cal B)}$ agree through 
$\gr_{\ga}+n$.
When there is no danger of confusion 
the superscript will be dropped from 
the $M$'s as well as from the embeddings 
$i_{\ga,\gb}:M_{\ga}\to M_{\gb}$.

\item
Iteration trees are pseudo-iteration trees $({\cal T},{\cal B})$
of length $(\gth,0)$, that is ${\cal B}=\seq{B_0}$ is a single
premouse. In this case it is customary to denote its length 
by $\gth$, rather than $(\gth,0)$.
On the other hand, 
any iteration tree $\cal T$ of length $\gth$
on a model $M$ can be
construed as a pseudo-iteration tree
of length $(\gth,\gl)$, any $\gl<\gth$. (Just forget
about the tree structure below $\gl$ and take
$B_{\ga}= M_{\ga}^{\cal T}$.)

\item
Implicit in 2.(c) of the above definition, is that
$M_{\ga}$ and $M_{\gb}$ agree through $\gr_{\gb}+n$,
when $\gb<\ga$. This is proved by induction on $\ga$.

\item
Note that plus-$n$ implies plus-$m$, for  $n>m$.
In this paper we will be mainly concerned with
plus-1 and plus-2 trees.

\item
The above definition, when restricted to ordinary iteration trees,
is less general than the one in \cite{Iterationtrees} 
as it covers only {\em normal}\/ iteration trees.
The reason we chose to eschew non-normal pseudo-iteration trees
was to avoid awkward notation.
On the other hand, the comparison process for models of the form
$L[\vec{\cal E}]$ entails normal trees only, so our 
present definition is not too restrictive.

\end{enumerate}

In order to prove a few basic results about pseudo-iteration 
trees we must restrict our definition a bit. 

\begin{definition}
Let $({\cal T,\cal B})$ be a plus-$n$ pseudo-iteration tree 
of length $(\gth,\gl)$.

{\sl (a)} $({\cal T,\cal B})$ is non-overlapping
if $lh(E_{\ga})<\hbox{\rm crit}(E_{\gb})$,
whenever $\ga+1=<_T\mbox{-pred}(\gb+1)$
and $\gb+1<\gth$.

{\sl (b)} $({\cal T,\cal B})$ is internal
if $\gth\in B_0$, $\seq{B_{\ga}\mid\ga<\gl}\in B_{\gl}$ 
and $B_{\gl}\models\lq\lq | B_{\ga}|=|V_{\gr_{\ga}+n}|$
and $B_{\ga}$ is $2^{\aleph_0}$-closed".

{\sl (c)}
If all the extenders $E_{\ga}$ are $2^{\aleph_0}$-closed
in the model they appear, i.e. $M_{\ga}\models\lq\lq 
\ult(V,E_{\ga})$ is $2^{\aleph_0}$-closed", then $({\cal T,\cal B)}$
is said to be $2^{\aleph_0}$-closed.

\end{definition}
Notice that if $({\cal T,\cal B})$ is internal
plus-$n$, $n\geq 1$, then 
$\seq{(B_{\ga},\gr_{\ga})\mid\ga<\gb}\in B_{\gb}$,
for any $\gb\leq\gl$, as
such sequence can be coded as a subset of 
$V_{\gr_{\gb}+(n-1)}\cap B_{\gl}$ and $B_{\gl}$ and $B_{\gb}$
agree up to $\gr_{\gb}+n$.

In the next two lemmata we derive some easy consequences
of $\cal T$ being non-overlapping or $2^{\aleph_0}$-closed.

\begin{lemma}
\label{successor}
Suppose we are given a countable, 
internal, non-overlapping, plus-1
pseudo-iteration tree $(\cal T, \cal B)$  of length 
$(\gth+1,\gl)$.
Assume also that $E\in M_{\gth}$ is an extender that can be 
applied to some earlier model $M_{\gn}$ in a non-overlapping way. 
Then $\ult(M_{\gn}, E)$
is wellfounded.
\end{lemma}

\begin{proof} The proof is an obvious modification of
Theorem 1.2 of \cite{Wfdd}. Let $\ga=\rt
(\gn)$, let $M_{\gth+1}=\ult(M_{\gn},E)$ and
let $i_{\ga,\gth+1}=i^{M_{\gn}}_E\circ i_{\ga,\gn}$.
As the pseudo-iteration tree is non-overlapping, every element in  
$\ult(M_{\gn}, E)$ is of the form $i_{\ga ,\gth +1}(f)(a)$, for
some $a\in [\gb]^{<\go}$, where $\gb= lh(E)$. 
[This follows from a straightforward induction on $\nu$:
the only place where the \lq\lq non-overlapping" 
condition is used is when $\nu$ is limit.]
Suppose, towards a contradiction, that $M_{\gth +1}$
is illfounded. 
As $M_{\gth+1}=\ult(M_{\nu}, E)$ agrees with 
$\ult(M_{\gth},E)$ through $i_E(\gk)+1$, then
$V_{i_E(\gk)+1}^{M_{\gth+1}}\in\mbox{\rm WFP}(M_{\gth+1})$.
By absoluteness $B_{\gl}\models\lq\lq M_{\gth+1}$ is
illfounded", hence there is a sequence of functions
$\seq{f_n\mid n\in\go}\in B_{\gl}$, with
each $f_n\in B_{\ga}$, and $a_n\in [\gb]^{<\go}$ 
such that $\seq{i_{\ga,\gth+1}(f_n)(a_n)\mid n\in\go}$
forms an infinite descending chain in $M_{\gth+1}$. 
As $B_{\ga}$ is $\go$-closed inside $B_{\gl}$,
$\seq{f_n\mid n\in\go}\in B_{\ga}$, hence the set 
$Y=\{i_{\ga,\gth+1}(f_n)(b)\mid n\in\go, 
b\in [\gb]^{<\go}\}\in M_{\gth+1}$.
Working in $M_{\gth+1}$, observe that
$|Y|\leq\gb$ so $Z$, its transitive collapse,
belongs to $V_{i_E(\gk)+1}^{M_{\gth+1}}$. 
But $Y$ is really illfounded (in $V$), and so
must be $Z$. Thus $V_{i_E(\gk)+1}^{M_{\gth+1}}$ cannot
be wellfounded: a contradiction.
\end{proof}

By inspecting the proof above we see that for the first
$\omega$ models, the non-overlapping condition is not needed.

\begin{corollary}
\label{corsuccessor}
If $({\cal T},{\cal B})$ is countable, internal, plus-1,
of length $(\gl+n+1,\gl)$, and such that
$E\in M_{\gl+n}$ can be applied to some previous
$M_{\ga}$, then $\ult(M_{\ga},E)$ is wellfounded.
\end{corollary}

\begin{lemma}
\label{finitesup}
If $(\cal T,\cal B)$ is plus-1, $2^{\aleph_0}$-closed, 
internal and of length $(\gth,\gl)$,
then every model $M_{\ga}^{(\cal T,\cal B)}$ 
is $2^{\aleph_0}$-closed in $B_{\gl}$, for
$\ga<\min (\gth,\gl+\go)$.
\end{lemma}

\noindent{\bf Remark:}
In general, $M_{\ga}^{\cal T}$ fails 
to be $\go$-closed for $\ga\geq\gl+\go$, 
so the lemma cannot be improved.

\medskip

\begin{proof}
By induction on $\ga$. 
We may assume $\ga=\gb+1>\gl$ as when 
$\ga\leq\gl$ the result follows at once.
Let $M_{\ga}=\ult(M_{\gg},E)$, where 
$E=E_{\gb}$ is a $(\gk,\nu)$-extender,
and $\gg=<_T\mbox{-pred}(\ga)$.
Given $\seq{(a_{\xi},f_{\xi})\mid\xi<2^{\aleph_0}}\in B_{\gl}$,
with $[a_{\xi},f_{\xi}]^{M_{\gg}}_E\in M_{\ga}$, we
want to show that
$$
\seq{i_E^{M_{\gg}}(f_{\xi})(a_{\xi})\mid\xi<2^{\aleph_0}}=
\seq{[a_{\xi},f_{\xi}]^{M_{\gg}}_E\mid\xi<2^{\aleph_0}}
\in M_{\ga}
$$
[Here, and in the rest of this proof, $2^{\aleph_0}$
means $(2^{\aleph_0})^{B_{\gl}}$.]

First notice that 
$\seq{a_{\xi}\mid\xi<2^{\aleph_0}}\in M_{\ga}$: by
the inductive hypothesis applied to $\gb$ and 
$2^{\aleph_0}$-closure of $E$,
$\seq{a_{\xi}\mid\xi<2^{\aleph_0}}$ belongs to
$\ult (M_{\gb}, E)$, which agrees with
$\ult(M_{\gg},E)=M_{\ga}$ through $i_E(\gk)+1$. 
Hence $\seq{a_{\xi}\mid\xi<2^{\aleph_0}}\in M_{\ga}$.

As each $f_{\xi}\in M_{\gg}$ and $M_{\gg}$ is 
$2^{\aleph_0}$-closed inside $B_{\gl}$, then
$\seq{f_{\xi}\mid\xi<2^{\aleph_0}}\in M_{\gg}$,
hence $F\in M_{\gg}$ where we set
$$
F\big(\seq{b_{\xi}\mid\xi<2^{\aleph_0}}\big)(\eta)=
f_{\eta}(b_{\eta})
$$
for all sequences $\seq{b_{\xi}\mid\xi<2^{\aleph_0}}\in M_{\gg}$
with $b_{\xi}\in [\nu]^{|a_{\xi}|}$. Thus
$$
i_E^{M_{\gg}}\big(F\big)(\seq{a_{\xi}\mid\xi<2^{\aleph_0}})=
\seq{i_E^{M_{\gg}}(f_{\xi})(a_{\xi})\mid\xi <2^{\aleph_0}}
\in M_{\ga}
$$ 
and this is what we had to prove.
\end{proof}

If $\cal T$ is a pseudo-iteration tree
on ${\cal B}=\seq{B_{\ga}\mid\ga\leq\gl}$
and ${\cal C}=\seq{C_{\ga}\mid\ga\leq\gl}$
are premice such that $B_{\ga}\subset C_{\ga}$
and $\gd(B_{\ga})=\gd(C_{\ga})$,
then $\cal T$ need not to be a pseudo-iteration
tree on $\cal C$: it is quite possible that for some 
$\gl<\gg< lh({\cal T})$, the $\gg$th model 
$M_{\gg}^{(\cal T,\cal C)}$ 
is illfounded, while  the corresponding model
on the $\cal B$-side is wellfounded, as required
by our definition.
Similarly, if $C_{\ga}\subset B_{\ga}$ and 
$\gd(C_{\ga})=\gd(B_{\ga})$, then again
$({\cal T,\cal C})$ can fail to be a pseudo-iteration tree,
as at some stage $\gg>\gl$, $E_{\gg}^{\cal T}$ 
might not belong to $M_{\gg}^{(\cal T,\cal C)}$.
In order to find sufficient conditions 
on $\cal C$ for $(\cal T,\cal C)$
to be a pseudo-iteration tree we introduce
the notion of embedding.

\begin{definition}
\label{embedding}
Suppose
$(\cal T,\cal B)$ and $(\cal S,\cal C)$ 
are pseudo-iteration trees of length 
$(\gth,\gl)$, $(\gth,\nu)$, respectively, 
and $\gl\leq\gn$.
A family of maps $\Pi=\seq{\gp_{\ga}\mid\ga<\gth}$ 
is an {\em embedding of pseudo-iteration trees\/},
$\Pi:({\cal T},{\cal B})\to ({\cal S},{\cal C})$, 
if there are ordinals $\eta_{\ga}\leq M_{\ga}^{\cal S}
\cap\On$ such that
\begin{enumerate}

\item
each $\gp_{\ga}:M_{\ga}^{\cal T}\to 
M_{\ga}^{\cal S}\cap V_{\eta_{\ga}}$ 
is an elementary embedding,
$\pi_{\ga}(\gd(M_{\ga}^{\cal T}))=\gd(M_{\ga}^{\cal S})$,
$\pi_{\ga}(\gr_{\ga}^{\cal T})=\gr_{\ga}^{\cal S}$
and $\pi_{\ga}(E_{\ga}^{\cal T})=E_{\ga}^{\cal S}$;

\item
for $\ga,\gb\geq\nu$,
$\ga<_T\gb\iff\ga<_S\gb$ and, 
for $\ga\leq\nu<\gb$,

$\ga<_S\gb\iff(\ga<_T\gb$ and 
$\neg\exists\ga'(\ga<\ga'\leq\nu\wedge\ga'<_T\gb))$;

\item
if $\ga<_S\gb$ then
$\eta_{\ga}\in M_{\ga}^{\cal S}\iff
\eta_{\gb}\in M_{\gb}^{\cal S}
\iff i^{\cal S}_{\ga,\gb}(\eta_{\ga})=\eta_{\gb}$
and the diagram
$$
\begin{diagram}
\node{M^{\cal T}_{\gb}} \arrow[2]{e,t}{\pi_{\gb}}
\node[2]{M_{\gb}^{\cal S}\cap V_{\eta_{\gb}}} \\

\node{M_{\ga}^{\cal T}} \arrow{n,l}{i^{\cal T}_{\ga,\gb}}
\arrow[2]{e,b}{\pi_{\ga}}
\node[2]{M_{\ga}^{\cal S}\cap V_{\eta_{\ga}}} 
\arrow{n,r}{i_{\ga,\gb}^{\cal S}}
\end{diagram}
$$
commutes.
\end{enumerate}

If for each $\ga<\gth$, $\eta_{\ga}=M^{\cal S}_{\ga}\cap\On$,
then $\Pi$ is an {\em elementary\/} embedding.

If for each $\ga<\gth$, $\eta_{\ga}\in M^{\cal S}_{\ga}$,
then $\Pi$ is a {\em bounded\/} embedding. 
Any sequence $\seq{\eta_{\ga}'\mid\ga<\gth}$ 
with $\eta_{\ga}'\geq\eta_{\ga}$, is
called a bound for $\Pi$.

\end{definition}
Note that an embedding can be both bounded and elementary.
Also if $\Pi:({\cal T},{\cal B})\to ({\cal S},{\cal C})$ is 
an embedding, $({\cal T},{\cal B})$ is plus-$n$
($2^{\aleph_0}$-closed, non-overlapping) 
iff $({\cal S},{\cal C})$ is plus-$n$
($2^{\aleph_0}$-closed, non-overlapping).

A particular kind of embedding is
obtained via the {\em copy construction}
(see \cite{Iterationtrees}).
Given $(\cal T,\cal B)$,
a plus-$n$ pseudo-iteration tree of length $(\gth,\gl)$, 
and a family of premice ${\cal C}=
\seq{C_{\ga}\mid\ga\leq\gl}$ and embeddings
$\Pi=\seq{\gp_{\ga}\mid\ga\leq\gl}$,
$\pi_{\ga}:B_{\ga}\to C_{\ga}\cap V_{\eta_{\ga}}$
such that $\gp_{\ga}\re V_{\gr_{\ga}^{\cal T}+n}=
\gp_{\gb}\re V_{\gr_{\ga}^{\cal T}+n}$, 
for $\ga\leq\gb\leq\gl$, we define the copied tree
$\Pi {\cal T}={\cal S}$ by boot-strapping the
definition of the $\pi_{\ga}$'s for $\ga>\gl$:
For any $\gl\leq\nu\leq\gth$ we want 
$\seq{\pi_{\ga}\mid\ga\leq\nu}$ to be an 
embedding of $({\cal T}\re\nu ,{\cal B})$
into $(\Pi{\cal T}\re\nu ,{\cal C})$
such that for $\eta\leq\xi<\nu$,
$\gp_{\eta}\re M_{\eta}^{(\cal T,\cal B)}
\cap V_{\gr_{\eta}^{\cal T}+n}=
\gp_{\xi}\re M_{\xi}^{(\cal T,\cal B)}
\cap V_{\gr_{\xi}^{\cal T}+n}$.
Thus if $\nu=\xi+1$ and $\gg=<_T\mbox{-pred}
(\nu)$ we let 
$M_{\nu}^{\Pi\cal T}=\ult\Big(
M_{\gg}^{\Pi\cal T},\pi_{\xi}(E_{\xi}^{\cal T})\Big)$,
if it is wellfounded and let 
$\pi_{\nu}:M_{\nu}^{\cal T}
\to M_{\nu}^{\Pi\cal T}$ be defined by
$$
\pi_{\nu}\Big([a,f]^{M}_{E}\Big)=
[\pi_{\xi}(a),\pi_{\gg}(f)]^N_F
$$
where $M=M_{\gg}^{\cal T}$,
$E=E^{\cal T}_{\xi}$, $N=M_{\gg}^{\Pi\cal T}$
and $F=\pi_{\xi}(E_{\xi}^{\cal T})=E_{\xi}^{\Pi\cal T}$.
If $\nu$ is limit, let $M_{\nu}^{\Pi\cal T}=
\lim_{\gg<_T\nu} M_{\gg}^{\Pi\cal T}$, 
if such direct limit is wellfounded, and 
$\pi_{\nu}$ is the limit map.
If at some stage $\nu<\gth$ we hit an illfounded
model $M_{\nu}^{\Pi{\cal T}}$, then we stop the 
construction and declare the length of 
$\Pi{\cal T}$ to be $(\nu,\gl)$.

If $lh(\Pi{\cal T})=lh({\cal T})$, 
then we say that $\cal T$ can be copied
on $\cal C$ via $\Pi$. Also, by a slight abuse
of notation, the system of maps 
$\seq{\pi_{\ga}\mid\ga<\gth}$
is still denoted by $\Pi$. Observe also that if
$\Pi:({\cal T},{\cal B})\to
(\cal S,\cal C)$ is obtained from copying 
via $\Pi$ and is a bounded embedding, 
then it is enough to specify the  
bounds on $\cal C$, i.e. it is enough to give
$\seq{\eta_{\ga}\mid\ga<lh(\cal C)}$.

In the case $({\cal T},{\cal B})$ 
is internal and $\Pi$ is elementary and 
$lh({\cal B})=\gl$, then $\Pi$ and $\cal C$
can be retrieved from $\pi_{\gl}$ and $C_{\gl}$ as
$C_{\ga}=\gp_{\gl}(B_{\ga})$ and
$\gp_{\ga}=\gp_{\gl}\re B_{\ga}$.

We should also notice that in order to run the copy 
construction the $\pi_{\ga}$'s need not 
to be fully elementary. If, for example, 
$B_{\ga}\subseteq C_{\ga}$, 
$\gd(B_{\ga})=\gd(C_{\ga})=\gd_{\ga}$ and 
$B_{\ga}$ and $C_{\ga}$
agree through $\gd_{\ga}$, 
then we can still try to copy $\cal T$ on $\cal C$
via the inclusion maps $\pi_{\ga}:
B_{\ga}\hookrightarrow C_{\ga}$.
[Of course $lh(\Pi{\cal T})<lh({\cal T})$ is possible.]

\begin{lemma}	
\label{remark}
Suppose $(\cal T,\cal B)$ is a 
pseudo-iteration tree of length $(\gth,\gl)$ and
let $\gd_{\nu}=\gd(M_{\nu}^{({\cal T,\cal B})})$.
Suppose also ${\cal C}=\seq{C_{\ga}\mid\ga\leq\gl}$
are premice with $\gd(C_{\ga})=\gd_{\ga}$
and $\Pi=\seq{\pi_{\ga}\mid\ga\leq\gl}$
are embeddings such that, for all $\ga\leq\gl$,
$$
C_{\ga}\cap V_{\gd_{\ga}}=B_{\ga}\cap V_{\gd_{\ga}}
\quad\pi_{\ga}:C_{\ga}\to B_{\ga}\cap V_{\eta_{\ga}}
\quad\mbox{and}\quad\pi_{\ga}\re V_{\gd_{\ga}}
\subseteq\mbox{id}
$$
where $\eta_{\ga}\leq B_{\ga}\cap\On$. 
Then $\cal T$ can be construed as a 
pseudo-iteration tree on $\cal C$ 
and $\Pi$ copies $({\cal T},{\cal C})$ 
to $({\cal T},{\cal B})$. Moreover for $\nu<\gth$,
$$
M_{\nu}^{({\cal T},{\cal C})}\cap V_{\gd_{\nu}}=
M_{\nu}^{({\cal T},{\cal B})}\cap V_{\gd_{\nu}}
\qquad\mbox{and}\qquad
\pi_{\nu}\re V_{\gd_{\nu}}\subseteq\mbox{id}
$$
\end{lemma}

\begin{proof}
We verify by induction on $\nu$ that 
$M_{\nu}^{({\cal T},{\cal C})}$ is 
wellfounded, that it agrees with 
$M_{\nu}^{({\cal T},{\cal B})}$ 
through $\gd_{\nu}$ and that the 
copy map $\pi_{\nu}$ is the identity
on $V_{\gd_{\nu}}$.

Suppose $\gl<\nu+1<\gth$ and let $\xi=
<_T\mbox{-pred}(\nu+1)$. By the agreement 
between $M_{\nu}^{({\cal T},{\cal C})}$ 
and $M_{\nu}^{({\cal T},{\cal B})}$, 
$E=E^{\cal T}_{\nu}\in M_{\nu}^{({\cal T},{\cal C})}$. 
Also $M_{\xi}^{({\cal T},{\cal C})}$ and 
$M_{\nu}^{({\cal T},{\cal C})}$
agree (at least) through $\gr_{\xi}+1$,
hence $E$ can be applied to $M_{\xi}^{({\cal T},{\cal C})}$. 
Let $\pi_{\nu+1}:M_{\nu +1}^{({\cal T},{\cal C})}\to
M_{\nu +1}^{({\cal T},{\cal B})}\cap V_{\eta_{\nu+1}}$
be given by
$$
\begin{array}{lll}
\pi_{\nu+1}([a,f]^{M_{\xi}^{({\cal T},{\cal C})}}_E) & = &
[\pi_{\nu}(a),\pi_{\xi}(f)]^{M_{\xi}^{({\cal T},{\cal B})}}_{\pi_{\nu}(E)} \\
 & & \\
 & = & [a,\pi_{\xi}(f)]^{M_{\xi}^{({\cal T},{\cal B})}}_E
\end{array}
$$
$\pi_{\nu+1}$ is well-defined and elementary, 
as $\pi_{\nu}$ is the identity on $V_{\gd_{\nu}}$ 
and $E\in V_{\gd_{\nu}}$.
Hence $M_{\nu +1}^{({\cal T},{\cal C})}$
is wellfounded. As 
$M_{\xi}^{({\cal T},{\cal C})}$ and
$M_{\xi}^{({\cal T},{\cal B})}$ agree through $\gd_{\xi}$,
$M_{\nu +1}^{({\cal T},{\cal C})}$ and
$M_{\nu +1}^{({\cal T},{\cal B})}$ agree
through $i_{\xi,\nu+1}^{\cal S}(\gd_{\xi})=
i_{\xi,\nu+1}^{\cal T}(\gd_{\xi})=\gd_{\nu+1}$.
Similarly $\pi_{\nu+1}\re V_{\gd_{\nu+1}}$ 
is shown to be the identity.

The case when $\gl<\nu<\gth$ is limit is left to the reader.
\end{proof}

The very same argument shows that if 
$V_{\gd_{\ga}}\cap B_{\ga}\subseteq 
C_{\ga}\subseteq B_{\ga}$ and $\cal T$ 
is a pseudo-iteration tree on $\cal B$, 
then $\cal T$ can be construed on $\cal C$ and 
$M_{\nu}^{({\cal T},{\cal C})}\cap V_{\gd_{\nu}}=
M_{\nu}^{({\cal T},{\cal B})}\cap V_{\gd_{\nu}}$, 
all $\nu<\gth$ and $\gd_{\nu}=\gd(M_{\nu}^{(\cal T,\cal B)})$. 
In fact this is almost 
a corollary of the preceding lemma, 
except for the fact that the inclusion maps 
$\pi_{\ga}:C_{\ga}\hookrightarrow B_{\ga}$
do not form an embedding in our official sense.
[See also the remarks after the proof of the next lemma.]

The next result shows that we can truncate a $B_{\ga}$
at a rank without affecting the illfoundedness
of a given branch. 

\begin{lemma}
Let $(\cal T, \cal B)$ be 
internal, plus-1 pseudo-iteration tree 
of length $(\gth,\gl)$, $\gth<\go_1$, and let
$\gd_{\ga}=\gd(B_{\ga})$, for $\ga\leq\gl$.

\noindent {\sl (a)\/}.
Suppose $b$ is an illfounded branch
with root $\ga$ and suppose  
$B_{\ga}\models |V_{\gd_{\ga}}|<\gd^{\ast}$. 
Then $b$ is illfounded below $\gd^{\ast}$,
that is: the least ordinal of $B_{\ga}$
sent by $i_{\ga, b}$ into the illfounded
part of $M_b^{({\cal T},{\cal B})}$ 
is $<\gd^{\ast}$.

\smallskip

\noindent {\sl (b)\/}.
Suppose $\gth=\nu+1$, $\ga=\rt(\gb)$, $\gb<\nu$
and $\gd^{\ast}$ is such that 
$B_{\ga}\models |V_{\gd_{\ga}}|<\gd^{\ast}$. 
Suppose also that $M_{\nu}^{({\cal T},{\cal B})}\models\lq\lq E$ 
is an extender with critical point $\leq\gr_{\gb} "$ 
and that $\ult(M_{\gb}^{({\cal T},{\cal B})}, E)$ is illfounded.
Then the least ordinal sent by $i_E^{({\cal T},{\cal B})}\circ
i_{\ga,\gb}$ into the illfounded part of the ultrapower is
$<\gd^{\ast}$.
\end{lemma}

\begin{proof}
{\sl (a).\/} Working inside $B_{\gl}$ choose 
a cofinal sequence $\gb_n\in b$,
with $\gb_0=\ga$, and ordinals
$\xi_n$ such that 
$i_{\gb_n,\gb_{n+1}}(\xi_n)>\xi_{n+1}$,
witnessing the illfoundedness of 
$M^{(\cal T,\cal B)}_b$.
Pick $\gz>\gd_{\gl}$ large enough
so that all the relevant stuff is in $V_{\gz}$.
We must consider whether or not $\ga=\gl$.

Suppose $\ga=\gl$. Let $C_{\gl}$ be the transitive
collapse of the Skolem hull,
computed inside $B_{\gl}$,
$$
C_{\gl}\cong \hull^{V_{\gz}}\big(V_{{\gd}_{\gl}}\cup
\{\seq{(\xi_n,\gb_n)\mid n\in\go}\}\big)
$$
and let $\pi_{\gl}$ be the inverse of 
the transitive collapse,
$C_{\gb}=B_{\gb}$ and 
$\pi_{\gb}=\mbox{id}\re C_{\gb}$, for $\gb<\gl$.
Lemma \ref{remark} implies that 
$(\cal T,\cal C)$ is a pseudo-iteration tree 
that copies to $(\cal T,\cal B)$.
Moreover $M_b^{({\cal T},{\cal C})}$ 
is illfounded via the ordinals
$\pi^{-1}_{\gl}(\xi_n)=\bar{\xi}_n$.
As $\bar{\xi}_0\in C_{\gl}$, then 
$\bar{\xi}_0<|V_{\gd_{\gl}}|^+ \leq\gd^{\ast}$.
As the Skolem hull above was
computed inside $B_{\gl}$, then
$C_{\gl}\subset B_{\gl}$, hence $\gF$ copies 
$({\cal T},{\cal C})$ to $({\cal T},{\cal B})$, 
where $\gf_{\gb}=\pi_{\gb}$, $\gb\leq\gl$, 
are the identity maps.  By commutativity
of the copy maps and the iteration embeddings
$$
i^{({\cal T},{\cal B})}_{\gb_n,\gb_{n+1}}
(\gf_{\gb_n}(\bar{\xi}_n))>\gf_{\gb_{n+1}}
(\bar{\xi}_{n+1})
$$
and $\gf_{\gb_0}(\bar{\xi}_0)=
\gf_{\gl}(\bar{\xi}_0)=\bar{\xi}_0$.
Thus the least ordinal mapped by
$i_{\gl,b}^{({\cal T},{\cal B})}$ into 
the illfounded part is 
$\leq\bar{\xi}_0<\gd^{\ast}$.
This completes the proof in the case when $\ga=\gl$.

\medskip

Suppose now $\ga<\gl$. 
We cannot simply repeat word-by-word the argument 
above, as the sequence $\seq{\xi_n\mid n\in\go}$
cannot be taken to be in $B_{\ga}$.
The plan is to get a countable copy 
$(\bar{\cal T}, \bar{\cal B})$ of the
tree belonging to $B_{\ga}$ and then
internalize the construction in $B_{\ga}$.
Let $\bar{B}_{\gl}$ be the collapse of the
countable Skolem hull, computed inside $B_{\gl}$,
$$
\bar{B}_{\gl}\cong\hull^{V_{\gz}}\big(\gth +1\cup
\{{\cal T},\seq{(\xi_n,\gb_n)\mid n\in\go}\}\big)
$$
and let $\bar{\pi}_{\gl}$ be the inverse of 
the collapsing function, 
$\bar{\pi}_{\gl}(\bar{\cal T})=\cal T$. 
Set also 
$\bar{B}_{\gb}=\bar{\pi}_{\gl}^{-1}(B_{\gb})$
and $\bar{\pi}_{\gb}=\bar{\pi}_{\gl}\re \bar{B}_{\gb}
:\bar{B}_{\gb}\to B_{\gb}$ for all $\gb<\gl$.
By elementarity of $\bar{\pi}_{\gl}$, 
$(\bar{\cal T},\bar{\cal B})$ 
is a pseudo-iteration tree and $b$ is illfounded
via the ordinals
$\bar{\xi}_n=\bar{\pi}^{-1}_{\gl}(\xi_n)$.
As the Skolem hull was taken inside $B_{\gl}$,
then $\bar{\pi}_{\ga}\in B_{\gl}$, it is countable and
$\bar{\pi}_{\ga}\subseteq B_{\ga}$,
hence $\bar{\pi}_{\ga}\in B_{\ga}$. Similarly 
$\seq{(\gb_n,\bar{\xi}_n)\mid n<\go},
(\bar{\cal T},\bar{\cal B})\in B_{\ga}$.
Let $\gg$ be large enough so that 
$\bar{\pi}_{\ga}, (\bar{\cal T},\bar{\cal B})
\in B_{\ga}\cap V_{\gg}$ and let
$C_{\ga}$ be transitive collapse of the following
Skolem hull, computed inside $B_{\ga}$
$$
C_{\ga}\cong\hull^{V_{\gg}}\big(V_{\gd_{\ga}}\cup
\{\bar{\gp}_{\ga}\}\big)
$$
and let $h$ the collapsing map. 
Set $C_{\gb}=B_{\gb}$, $\gp_{\gb}=\bar{\gp}_{\gb}$,
for $\gb\neq\ga$ and $\pi_{\ga}=h(\bar{\pi}_{\ga})$.
Then $\Pi$ copies $(\bar{\cal T},\bar{\cal B})$
to $(\cal T,\cal C)$, $M_b^{({\cal T},{\cal C})}$ 
is illfounded as witnessed by
$\seq{\pi_{\ga}(\bar{\xi}_n)\mid n<\go}$.
We now argue as in the case when 
$\ga=\gl$. Letting 
$\gf_{\gb}:C_{\gb}\hookrightarrow B_{\gb}$, 
$\gb\leq\gl$, be the inclusion maps, then
$\gF$ copies $({\cal T},{\cal C})$ to
$({\cal T},{\cal B})$ and 
$$
i_{\gb_n,\gb_{n+1}}^{({\cal T},{\cal B})}
\Big(\gf_{\gb_n}(\gp_{\ga}(\bar{\xi}_n))\Big)
>\gf_{\gb_{n+1}}(\gp_{\ga}(\bar{\xi}_{n+1})).
$$
Thus the least element mapped by 
$i_{\ga,b}^{({\cal T},{\cal B})}$ is 
$\leq\gf_{\gb_0}(\pi_{\ga}(\bar{\xi}_0))=
\pi_{\ga}(\bar{\xi}_0)<|C_{\ga}|^+<\gd^{\ast}$.
This concludes the proof of part {\sl (a)}.

\smallskip
{\sl (b).\/} 
The proof of this case is very similar to the 
one of {\sl (a)\/}, so we only indicate the main changes, 
leaving the details to the reader.
Let $[a_n,f_n]_E^{M_{\nu}^{({\cal T},{\cal B})}}$
witness the illfoundedness of the last ultrapower.
By absoluteness the $f_n$, $a_n$ can be 
taken to be inside $B_{\gl}$.
By replacing $\gb_n$, $\xi_n$ with $a_n$, $f_n$,
the proof adapts {\em verbatim}.
\end{proof}

The careful reader might question a few 
steps in the proof above: the collapses 
of those hulls are not, in general, premice, 
so we should not be allowed to build 
pseudo-iteration trees on them. 
One way to fix this problem would 
be to start with base premice 
$B_{\ga}$'s with arbitrarily large 
cut-off points. The other way, which 
we implicitly followed, is to relax 
a bit our official definition of 
pseudo-iteration tree, so that 
$({\cal T},{\cal C})$ makes sense 
even if the $C_{\ga}$'s don't satisfy 
replacement for domains of bounded rank.
The only difference is that the tree 
embeddings $i_{\ga,\gb}$ are only $\gS_0$-elementary
which is enough, anyway, to show that the branch $b$ 
is illfounded via the (images of the) $\bar{\xi}_n$'s.

\begin{corollary}
\label{truncation}
Suppose $(\cal T,\cal B)$ is 
internal, plus-1, of length
$(\gth,\gl)$, $\gth<\go_1$.
Suppose also that, for all $\ga\leq\gl$,
the $\gg_{\ga}$'s are cut-off points of the $B_{\ga}$'s,
and let $C_{\ga}=B_{\ga}\cap V_{\gg_{\ga}}$.
Then $(\cal T,\cal C)$ is a plus-1 pseudo-iteration 
tree and for any $\nu<\gth$ with root $\ga$
$$
M_{\nu}^{(\cal T,\cal C)}=
M_{\nu}^{(\cal T,\cal B)}\cap V_{i_{\ga,\gn}(\gg_{\ga})}.
$$
Moreover if $b$ is a branch 
$$
M_b^{(\cal T,\cal B)}\mbox{ is wellfounded }
\quad\iff\quad M_b^{(\cal T,\cal C)}\mbox{ is wellfounded,}
$$
and if $\gth=\gt +1$, $E$ is an extender in
$M_{\gt}^{({\cal T},{\cal B})}$ with critical 
point $\leq\gr_{\ga}$, then
$$
\ult(M_{\ga}^{(\cal T,\cal B)},E)
\mbox{ is wellfounded }\quad\iff\quad 
\ult(M_{\ga}^{(\cal T,\cal C)},E)\mbox{ is wellfounded.}
$$
\end{corollary}

\begin{proof}
The result follows from the last two lemmata 
and the fact that for any premouse $M$ and 
any cut-off point $\gg$, $M\models |V_{\gd(M)}|<\gg$.
\end{proof}

So far we only really used that 
the $B_{\ga}$'s are $\go$-closed 
inside $B_{\gl}$, rather than 
$2^{\aleph_0}$-closed.
The reason for requiring the stronger 
closure property in the definition 
of \lq\lq internal" will be clear 
from the proof of Lemma \ref{chunks}.
In order to get to it we must first 
introduce the notion of {\em support\/} 
for pseudo-iteration trees. 
This is the generalization to our 
present set up of the notion 
defined in \cite{Manywoodins}.

\begin{definition}
Let $T$ be a tree ordering on $\gth$ with 
$\gl +1$ roots. A set $X\subseteq\gth$ is 
$T$-compatible iff
\begin{enumerate}

\item
$X\cap (\gl+1)\neq\emptyset$.

\item
If $\ga+1\in X$ and $\ga+1>\gl$, then
$\ga,<_T\mbox{-pred}(\ga+1)\in X$.

\item
Suppose $\gg\in X$ is limit, $\gg>\gl$.
Then $\rt_T(\gg)\in X$.

If there is a largest $\ga\in X$ 
such that $\ga<_T\gg$, then $\gb\in X$,
where $\gb+1$ is
least such that $\ga<_T\gb+1<_T\gg$.

\end{enumerate}
\end{definition}
It is easy to see that if $X\subseteq\gl$ 
or $X\in\gth$, then $X$ is $T$-compatible. 
Clause (2) implies that if $\gb+n\in X$ 
and $\gb\geq\gl$, then $\gb+i\in X$, for all $i\leq n$.
Clause (3) implies that, letting 
$Y=(X\cap(\gb+1))\cup\{\gb+1\}$
and $Z=X\cap (\gg+1)$, 
$(Y, T\cap Y\times Y)\cong(Z, T\cap Z\times Z)$.

\begin{definition}
\label{definitionsupport}
Let $({\cal T},{\cal B})$ be of length $(\gth,\gl)$. 
A set $X\subseteq\gth$ is a support iff $X$
is $T$-compatible and there are elementary substructures
$\big(M_{\ga}^{(\cal T,\cal B)}\big)_X\prec 
M_{\ga}^{(\cal T,\cal B)}$, for $\ga\in X$, such that

\begin{enumerate}

\item
if $\ga\in X\cap(\gl+1)$ then $\big(M_{\ga}\big)_X=M_{\ga}=B_{\ga}$;

\item
if $\ga+1\in X$ and $\ga+1>\gl$ and
$\gb=<_T\mbox{-pred}(\ga+1)$, then 
$E_{\ga},\gr_{\ga}\in\big(M_{\ga}\big)_X$ and 
$$
\big(M_{\ga+1}\big)_X=
\{[a,f]_{E_{\ga}}^{M_{\gb}}\mid f\in (M_{\gb})_X
\wedge a\in (M_{\ga})_X\};
$$

\item
if $\ga\in X$ then $Y=X\cap (\ga+1)$ is
a support and for all $\gb\in Y$, 
$\big(M_{\gb}\big)_X=\big(M_{\gb}\big)_Y$;

\item
suppose $\gg\in X$ is limit and $\gg>\gl$ 
and let $A=\{\nu\in X\mid \gn<_T\gg\}$:

	\begin{enumerate}

	\item

	if $A$ has limit order type, let
	$$
	\big(M_{\gg}\big)_X=
	\bigcup_{\nu\in A}i_{\gn,\gg}^{\prime\prime}
	\big(M_{\gn}\big)_X\; ;
	$$

	\item
	if $A$ has a largest element $\ga$, let
	$\gb+1$ be least such that $\ga<_T\gb+1<_T\gg$.
	Then $Y=(X\cap(\gb+1))\cup\{\gb+1\}$ is a support and
	$$ 
	\big(M_{\gg}\big)_X=
	i^{\prime\prime}_{\gb+1,\gg} (M_{\gb+1})_Y.
	$$
\end{enumerate}

\end{enumerate}
\end{definition}

There is no suggestion that the $(M_{\ga})_X$'s 
should be transitive: in fact, in general, they
are not. The next lemma lists a few basic results 
about supports. The proof (a tedious but
straightforward induction on $\gth$) is left to the reader. 

\begin{lemma}
Fix $({\cal T},{\cal B})$ of length $(\gth,\gl)$. 
\begin{enumerate}

\item
If $X\setminus (\gl+1)$ is non-empty, $X$ a support,
then $\gl\in X$.

\item
If $\ga\in X\subseteq Y$ and $X$, $Y$
are supports then $\big(M_{\ga}\big)_X\prec
\big(M_{\ga}\big)_Y$.

\item
For any $Y\subseteq\gth$ there is 
a smallest support $X\supseteq Y$,
called the support generated by $Y$.
Moreover if $Y$ is finite, $X$ is finite too.

\item
For any $y\in M_{\ga}$ there 
is a finite support $X\supseteq\{\ga\}$ such that 
$y\in \big(M_{\ga}\big)_X$.

\end{enumerate}
\end{lemma}

Given a pseudo-iteration tree $(\cal T,\cal B)$
of length $(\gth,\gl)$ and a support $X$, a new
pseudo-iteration tree 
$\big({\cal T,\cal B}\big)_X=
\big({\cal T}_X,{\cal B}_X\big)$
is defined as follows. Let $h:\gth_X\to X$ be the 
enumerating function and let $\gl_X=$o.t.$(X\cap\gl)$.
\begin{enumerate}
\item
${\cal B}_X=\seq{B_{h(\ga)}\mid\ga\leq \gl_X}$;

\item
the tree ordering $<_{T_X}$ on $\gth_X$ is isomorphic
to $<_T\re X$ via $h$;

\item
let 
$$
j^{-1}_{\ga,X}:\big(M_{h(\ga)}^{(\cal T,\cal B)}\big)_X
\to M_{\ga}^{{(\cal T,\cal B)}_X}
$$
be the transitive collapse and set 
$\gr_{\ga}^X=j^{-1}_{\ga,X}(\gr_{h(\ga)})$ and
$E_{\ga}^{{\cal T}_X}=j^{-1}_{\ga,X}(E_{\ga}^{\cal T})$.

\end{enumerate}
It is immediate to verify that $({\cal T, 
\cal B})_X$ is a pseudo-iteration tree and that,
for $\ga<\gb<\gth_X$, if $\ga<_{T_X}\gb$, then
$h(\ga)<_T h(\gb)$ and
$$
\begin{diagram}
\node{M_{\gb}^{({\cal T,\cal B})_X}}
\arrow[2]{e,t}{j_{\gb,X}}
\node[2]{M_{h(\gb)}^{(\cal T,\cal B)}} \\
\node{M_{\ga}^{({\cal T,\cal B})_X}}
\arrow{n,l}{i_{\ga,\gb}^{({\cal T,\cal B})_X}}
\arrow[2]{e,b}{j_{\ga,X}}
\node[2]{M_{h(\ga)}^{(\cal T,\cal B)}}
\arrow{n,r}{i_{h(\ga),h(\gb)}^{(\cal T,\cal B)}}
\end{diagram}
$$
commutes. We call $J_{X,\gth}=
\seq{j_{\ga,X}\mid\ga\leq\gth_X}$
an immersion of $({\cal T,\cal B})_X$ 
in $(\cal T,\cal B)$.
Note that $J_{X,\gth}$ is {\em not\/} 
an embedding in the sense of 
Definition \ref{embedding}, 
unless $X=\gth$, in which case
it is the identity.
Also if $(\cal T,\cal B)$ is internal (plus-$n$, 
non-overlapping, $2^{\aleph_0}$-closed), 
so is $({\cal T,\cal B})_X$.

If $X\subseteq Y$ are supports 
for $(\cal T,\cal B)$ and $\pi$ is the collapse of
$Y$ and $W=\pi^{\prime\prime} X$, then, by a tedious
but straightforward verification, it can be shown that
$$
W \mbox{ is a support for }
({\cal T,\cal B})_Y \mbox{ and }
\Big(\,({\cal T,\cal B})_Y\Big)_W=
({\cal T,\cal B})_X.
$$
There is also an immersion 
$J_{X,Y}:({\cal T,\cal B})_X\to({\cal T,\cal B})_Y$, 
such that $J_{X,\gth}=J_{Y,\gth}\circ J_{X,Y}$.
Hence for supports $X\subseteq Y\subseteq Z$
$J_{X,Z}=J_{Y,Z}\circ J_{X,Y}$.
Summarizing: any pseudo-iteration tree 
$({\cal T,\cal B})$ of length $(\gth,\gl)$
is the direct limit of the system
$\seq{({\cal T,\cal B})_X,J_{X,Y}\mid
X\subseteq Y\subseteq\gth}$, with $J_{X,\gth}$ the limit maps. 

Suppose $\Pi:({\cal T,\cal B})\to(\cal S,\cal C)$
is an embedding of pseudo-iteration 
trees of length $(\gth,\gl)$,
$(\gth,\gn)$ respectively, and
suppose $X$ is a support for $(\cal T,\cal B)$. 
Then $X$ need not be a support for 
$({\cal S},{\cal C})$: in fact
$X\setminus (\nu+1)$ could be non-empty
and yet $\nu\notin X$.
Thus we set $\Pi X$ to be the support for
$(\cal S,\cal C)$ generated by $X$,
and for $\ga\in X$ we have
$$
\pi_{\ga}^{\prime\prime}\big(M_{\ga}^{\cal T}\big)_X
\prec\big(M_{\ga}^{\cal S}\big)_{\Pi X}. 
$$
On the other hand, if $\gl=\nu$ then $\Pi X=X$. 
In particular, if $\Pi:({\cal T},M)\to({\cal S},N)$ 
is an embedding of ordinary iteration trees then
a support for $\cal T$ is also a support for $\cal S$.

\begin{lemma}
\label{closure}
Suppose $(\cal T,\cal B)$ of length $(\gth,\gl)$, $\gth<\go_1$,
is plus-1,  internal and $2^{\aleph_0}$-closed.
Let $\seq{S_n\mid n<\go}$ be an increasing
sequence of finite supports such that
$\bigcup_n S_n=\gth$.
Then, for any $\ga<\gth$,
there is $n_0=n_0(\ga)$ such that
$\ga\in S_{n_0}$ and
$$
M_{\ga}=\bigcup_{n\geq n_0}\big(M_{\ga}\big)_{S_n}
$$
and each $\big(M_{\ga}\big)_{S_n}$
is $2^{\aleph_0}$-closed inside $B_{\gl}$.
\end{lemma}

\begin{proof}
As any element of $M_{\ga}$
belongs to $\big(M_{\ga}\big)_X$,
for some finite support $X$ containing $\ga$, 
and as $X\subseteq S_n$,
for $n$ sufficiently large,
$M_{\ga}$ is the increasing union 
of the $\big(M_{\ga}\big)_{S_n}$.
As for $2^{\aleph_0}$-closure, note that
$\big(M_{\ga}\big)_{S_n}$ is isomorphic
(via the transitive collapse) to a model
of ${(\cal T,\cal B)}_{S_n}$. As 
${(\cal T,\cal B)}_{S_n}$ is
internal, plus-1, $2^{\aleph_0}$-closed
and $B_{\gl}$ is the last model of 
${\cal B}_{S_n}$, the result follows
easily from Lemma \ref{finitesup}.
\end{proof}

The submodels $\big(M_{\ga}\big)_{S_n}$ will be called
sometimes {\em chunks\/} of $M_{\ga}$. 
The next result will be a key ingredient in the
main proofs of this paper.

\begin{lemma}
\label{chunks}
Suppose $\Pi:({\cal T,\cal B})\to({\cal S,\cal C})$ is
a bounded embedding with bounds $\seq{\eta_{\ga}\mid\ga<\gth}$
and that $({\cal T,\cal B})$ and
$({\cal S,\cal C})$ are internal, plus-1,
$2^{\aleph_0}$-closed of length $(\gth,\gn)$
and $(\gth,\gl)$, respectively, and $\gth<\go_1$.
Suppose also that 
$({\cal T,\cal B})\in C_{\nu}$ and
$C_{\nu}\models\forall\ga\leq\gl
\big(|B_{\ga}|\leq 2^{\aleph_0}\big)$.
Let $\seq{S_n\mid n<\go}\in C_{\nu}$ be an 
increasing sequence of finite supports
for $(\cal T,\cal B)$, such that 
$\bigcup_n S_n=\gth$.
Then, for any $\ga<\gth$, and $n$ such that $\ga\in S_n$,

\begin{enumerate}

\item
$\pi_{\ga}\re\big(M_{\ga}^{\cal T}\big)_{S_n}
\in M_{\ga}^{\cal S}$ and

\item
there is an elementary embedding
$$
\gf_{\ga}:M_{\ga}^{\cal T}\to
M_{\ga}^{\cal S}\cap V_{\eta_{\ga}}
$$
such that 
$\gf_{\ga}\re \big(M_{\ga}^{\cal T}
\big)_{S_n}=\pi_{\ga}\re\big(M_{\ga}^{\cal T}
\big)_{S_n}$ and $\gF$, the system of embeddings 
obtained from $\Pi$ by changing $\pi_{\ga}$ 
to $\gf_{\ga}$, is an embedding of pseudo-iteration trees
$\gF:{(\cal T,\cal B)}\to {(\cal S,\cal C)}$ with
the same bounds $\seq{\eta_{\ga}\mid\ga<\gth}$.
\end{enumerate}

\end{lemma}

\begin{proof}
As $({\cal T,\cal B})\in C_{\nu}$, every model
$M_{\ga}^{\cal T}$ belongs to every 
$M_{\gb}^{\cal S}$ and is of size $\leq 2^{\aleph_0}$.
Thus, for fixed $\ga<\gth$, and $n$ such that $\ga\in S_n$,
$$
M_{\ga}^{\cal S}\models \big(
M_{\ga}^{\cal T}\big)_{S_n}\mbox{ is
of size }\leq 2^{\aleph_0}.
$$
As $\pi^n=\pi_{\ga}\re\big(M_{\ga}^{\cal T}
\big)_{S_n}$ is an elementary embedding of 
$\big(M_{\ga}^{\cal T}\big)_{S_n}$
into $\big( M_{\ga}^{\cal S}\big)_{\Pi S_n}\cap V_{\eta_{\ga}}$
and $C_{\nu}\models |\pi^n|\leq 2^{\aleph_0}$, by 
Lemma \ref{closure} $\pi^n\in\big(M_{\ga}^{\cal S}
\big)_{\Pi S_n}\subseteq M_{\ga}^{\cal S}$,
proving thus part (1).

Now for (2). Fix $n<\go$ such that 
$\ga\in S_n$ and let 
${\cal V}\in M_{\ga}^{\cal S}$ 
be the tree of attempts to find a sequence
like $\seq{\pi^n,\pi^{n+1},\ldots}$. That is,
working inside $M_{\ga}^{\cal S}$, let
$$
\seq{\gt_0,\ldots,\gt_k}\in{\cal V}\iff
\gt_0=\pi^n, \gt_0\subseteq\ldots\subseteq\gt_k
$$ 
where $\gt_i:\big(M_{\ga}^{\cal T}\big)_{S_{n+i}}
\to V_{\eta_{\ga}}$ is an elementary embedding.
By part (1), $\pi^m\in M_{\ga}^{\cal S}$, 
for any $m\geq n$, so 
$\seq{\pi^i\mid n\leq i<\go}$ 
is a branch of $\cal V$ in $V$. 
By absoluteness there is a branch
$\seq{\gt_i\mid i<\go}\in M_{\ga}^{\cal S}$
and let $\gf_{\ga}=\bigcup_{i<\go}\gt_i$.
\end{proof}

%% file: strongpast.tex
%

In this section all iteration trees will be 
countable, $2^{\aleph_0}$-closed and
non-overlapping. Suppose we are given an 
iteration tree $({\cal T}, M)$ and a sequence
$\seq{({\cal T}_{\ga},P_{\ga})\mid\ga<\nu}$
of iteration trees ${\cal T}_{\ga}$ on $P_{\ga}$
of successor length $\gth_{\ga}+1$, 
together with a last extender
$E_{\gth_{\ga}}^{{\cal T}_{\ga}}\in M_{\gth_{\ga}}^{{\cal T}_{\ga}}$
and such that ${\cal T}_0={\cal T}\re\gth_0+1$
and $M=P_0$. We define, by induction on 
$0<\nu\leq lh({\cal T})$, what it means for
$\seq{({\cal T}_{\ga},P_{\ga})\mid\ga<\nu}$
to be a decomposition of $({\cal T},M)$.

\begin{itemize}
\item
If $\nu=1$, then $\seq{({\cal T},M)}$ is the only possible
decomposition of $({\cal T},M)$, hence $lh({\cal T})=\gth_0+1$.

\item
If $\nu$ is limit, then for every $\xi<\nu$
there is an ordinal $\gth<lh({\cal T})$, 
such that $\seq{({\cal T}_{\ga},P_{\ga})\mid\ga<\xi}$
is a decomposition of $({\cal T}\re\gth,M)$.

\item
If $\nu=\xi+1>1$, then there is $\gth<lh({\cal T})$
such that $\seq{({\cal T}_{\ga},P_{\ga})\mid\ga<\xi}$
is a decomposition of $({\cal T}\re\gth,M)$,
$lh({\cal T})=\gth+\gth_{\xi}+1$ and 
for $\ga<\gth_{\xi}$
$$
M_{\ga}^{{\cal T}_{\xi}}=M_{\gth+\ga}^{\cal T}
\qquad\mbox{and}\qquad
E_{\ga}^{{\cal T}_{\xi}}=E_{\gth+\ga}^{\cal T}.
$$
Hence in particular $P_{\xi}=M_{\gth}^{\cal T}$.
Moreover

\begin{enumerate}
\item
if $\xi$ is limit, then $P_{\xi}$ is the 
direct limit of a (cofinal in $\xi$)
sequence of $P_{\ga}$'s,

\item
if $\xi=\eta+1$, then $\gth$ is a successor and either
\begin{enumerate}

	\item
	$P_{\xi}=M_{\gth_{\eta}}^{{\cal T}_{\eta}}=$
	the last model of $({\cal T}_{\eta},P_{\eta})$,
	in which case we write $({\cal T}_{\xi},P_{\xi})
	\|({\cal T}_{\eta},P_{\eta})$, or else

	\item
	$P_{\xi}=\ult(M_{\gg}^{{\cal T}_{\gb}},
	E_{\gth-1}^{\cal T})$, where $\gb\leq\eta=\xi-1$ 
	and $\gg\leq\gth_{\gb}$, if $\gb<\eta$, or
        $\gg<\gth_{\gb})$, if $\gb=\eta$. In this case we write
	$({\cal T}_{\xi},P_{\xi})\bot ({\cal T}_{\eta},P_{\eta})$.
\end{enumerate}
(When there is no danger of confusion  we simply drop 
the $P_{\ga}$'s and write ${\cal T}_{\ga+1}\|{\cal T}_{\ga}$ 
or ${\cal T}_{\ga+1}\bot {\cal T}_{\ga}$.)
\end{enumerate}
\end{itemize}

\noindent
The idea here is that $({\cal T},M)$ can be written as 
a tree of trees $({\cal T}_{\ga},P_{\ga})$: 
every model $M_{\ga}^{\cal T}$ and extender $E_{\ga}^{\cal T}$
are of the form $M_{\gg}^{\cal T_{\gb}}$, 
$E_{\gg}^{{\cal T}_{\gb}}$, for some $\gg\leq\gth_{\gb}$,
and $\gb<\nu$. The pair $(\gb,\gg)$ is unique except when 
${\cal T}_{\xi+1}\|{\cal T}_{\xi}$.
In this case $(\xi,\gth_{\xi})$ and $(\xi+1,0)$ yield
the same model $M_{\gth_{\xi}}^{{\cal T}_{\xi}}=P_{\xi+1}=
M_0^{{\cal T}_{\xi+1}}$. Vice versa, for any pair
$(\gb,\gg)$ with $\gg\leq\gth_{\gb}$ and $\gb<\nu$,
there is a unique $\ga$, such that
$M_{\ga}^{\cal T}=M_{\gg}^{{\cal T}_{\gb}}$.
Let $\gb\otimes\gg$ be such $\ga$, with 
$$
\otimes :\{(\gb,\gg)\mid\gg
\leq\gth_{\gb}\mbox{ and }\gb<\nu\}\to lh({\cal T}).
$$

\begin{definition}
Let $\cal S=
\seq{({\cal T}_{\ga},P_{\ga})\mid\ga<\nu}$
be a decomposition of $({\cal T},M)$ and let
${\gb\otimes\gg}=\gk< lh({\cal T})$. By induction
on $\nu$, we will define:

\begin{enumerate}
\item
what it means for $\cal S$ to be {\em quasi-linear\/};

\item
the ordinal $\mbox{\bf B}(\gk)$, the
{\em back-up point\/} of $\gk$, 
which depends on the tree ordering $<_T$ only;

\item
the set $\mbox{\bf F}_{\nu}\subseteq lh({\cal T})$
of {\em forbidden nodes.}
\end{enumerate}

\begin{itemize}
\item
If $\nu=1$, then $\cal S=\seq{({\cal T},M)}$ is quasi-linear,
$\mbox{\bf B}(\gk)=0$ and $\mbox{\bf F}_1=\emptyset$.

\item
Suppose $\nu$ is limit. $\cal S$ is quasi-linear iff
$\forall\xi<\nu({\cal S}\re\xi\mbox{ is quasi-linear})$, 
$\mbox{\bf F}_{\nu}=\bigcup_{\xi<\nu}\mbox{\bf F}_{\xi}$,
and $\mbox{\bf B}(\gk)$ is the back-up point of $\gk$
as computed in ${\cal S}\re\gb+1$.

\item
Suppose $\nu=\xi+1$ and $\xi$ limit.
Set $\mbox{\bf F}_{\nu}=\mbox{\bf F}_{\xi}$.

If $\gb<\xi$, then $\mbox{\bf B}(\gk)$
has already been defined so we may assume $\gb=\xi$.

$\cal S$ is quasi-linear iff ${\cal S}\re\xi$ is
quasi-linear and either
\begin{itemize}
	\item
	there is $\xi_0<\xi$ such that for $\xi_0\leq\eta<\xi$, 
	${\cal T}_{\eta+1}\| {\cal T}_{\eta}$ 
	and $P_{\xi}$ is the direct limit of such $P_{\eta}$'s,
        and $\xi_0$ is least such.
	Set $\mbox{\bf B}(\gk)=\mbox{\bf B}({\xi_0\otimes 0})$. 
	Or else

	\item
	there is an increasing sequence $\xi_n\to\xi$ such that
	${\cal T}_{\xi_n+1}\bot {\cal T}_{\xi_n}$, 
	$(\xi_n+1)\otimes 0<_T \xi\otimes 0$. Set $\mbox{\bf B}(\gk)
	={\xi\otimes 0}$.
\end{itemize}
\item
Suppose $\nu=\xi+2$. If $\gb<\xi+1$, then 
$\mbox{\bf B}(\gk)$ has already been defined,
so assume $\gb=\xi+1$.

$\cal S$ is quasi-linear iff ${\cal S}\re\xi+1$
is quasi-linear and either
\begin{itemize}
	\item
	${\cal T}_{\xi+1}\|{\cal T}_{\xi}$,
	$\mbox{\bf B}(\gk)=\mbox{\bf B}\big({\xi\otimes 0}\big)$
	and $\mbox{\bf F}_{\nu}=\mbox{\bf F}_{\xi+1}$, or else

	\item
	${\cal T}_{\xi+1}\bot{\cal T}_{\xi}$, that is
	$P_{\xi+1}=\ult(M_{\eta}^{{\cal T}_{\gz}},
	E_{\gth_\xi}^{{\cal T}_{\xi}})$, and 
	${\eta\otimes \gz}\notin
	\mbox{\bf F}_{\xi+1}$; set 
	$\mbox{\bf F}_{\nu}=\mbox{\bf F}_{\xi+1}\cup
	[\gz\otimes\eta,\xi\otimes\gth_{\xi}]$ and 
	$\mbox{\bf B}(\gk)={\eta\otimes\gz}$.
\end{itemize}
\end{itemize}
\end{definition}
{\bf Remarks.}
\begin{enumerate}
\item
The sets $\mbox{\bf F}_{\ga}$'s are increasing, i.e.
$\ga<\gb\implies\mbox{\bf F}_{\ga}\subseteq\mbox{\bf F}_{\gb}$.
If $\ga\in\mbox{\bf F}_{\gb}$, then we are not allowed 
to visit the models $M_{\ga}$ past round $\gb$. In other words:
there is no $\gg\geq\gb$ such that $P_{\gg+1}=\ult(M_{\ga},E)$.
This is the content of Lemma \ref{next}
to be proved shortly.

\item
$\mbox{\bf B}(\gk)<\gk$, unless $\gk=0$ or 
$\gk={\xi\otimes 0}$, where $\xi$ is 
limit of an increasing sequence $\xi_n$, with 
${\cal T}_{\xi_n+1}\bot {\cal T}_{\xi_n}$. 
In this case $\mbox{\bf B}(\gk)=\gk$.

\item
If ${\cal T}_{\ga+1}\bot{\cal T}_{\ga}$, then 
$\mbox{\bf B}\big((\ga+1)\otimes 0\big)$ is the immediate 
$<_T$-predecessor of $(\ga+1)\otimes 0$.

\item
To make  the notation a bit simpler we shall write
$\mbox{\bf B}(\gb,\gg)$ rather than $\mbox{\bf B}({\gb\otimes\gg})$.
\end{enumerate}

A few words on the motivations 
behind the notion of quasi-linearity are in order here.
In the case of ${\cal W\cal G}(M,\nu)$, 
the weak iteration game of length $\nu$,
a tower of trees 
$\seq{({\cal T}_{\ga},P_{\ga})\mid\ga<\gn}$
is built, with $P_0=M$, 
${\cal T}_{\ga+1}\|{\cal T}_{\ga}$
and for any limit $\gl<\nu$, $P_{\gl}$ is the direct
limit of the $P_{\ga}$'s, $\ga<\gl$. 
Thus the resulting iteration tree can be
construed as a {\em linear iteration of 
iteration trees\/}.
In order to consider a wider spectrum of trees
for which we can still prove an iterability result,
we relax the \lq\lq linearity"
condition a bit. The trees ${\cal T}_{\ga}$
are arranged themselves in a tree ordering,
but this tree ordering is not too removed
from a linear ordering: whenever we \lq\lq go back"
and take the ultrapower of the $\gg$th model
of the $\gb$th iteration tree to start a new 
${\cal T}_{\ga+1}$, then we give up the right to
visit, from this point on, any model with index 
$(\gb',\gg')$, with 
$\gb\otimes\gg<\gb'\otimes\gg'<(\ga+1)\otimes 0$.
(For example: an \lq\lq alternating chain of iteration trees"
is not quasi-linear.) 

\begin{lemma}
\label{next}
Suppose $\seq{({\cal T}_{\ga},P_{\ga})\mid\ga<\nu}$
is a quasi-linear decomposition of $({\cal T},M)$.
If ${\cal T}_{\ga+1}\bot {\cal T}_{\ga}$ and 
$\mbox{\bf B}\big(\ga+1,0\big)={\gb\otimes\gg}$,
then 
$\forall\gk\forall\xi(\gb\otimes\gg\leq\gk\leq\ga\otimes\gth_{\ga}$
and $(\ga+1)\otimes 0\leq\xi<lh({\cal T})$ implies $\gk\not <_T\xi)$.
\end{lemma}

\begin{proof}
Deny. Choose counter-examples $\gk$ and $\xi$,
first minimizing $\xi$, and then taking $\gk$ 
as large as possible (relative to this $\xi$). 
Then $\gk$ is the immediate $<_T$-predecessor of $\xi$:
otherwise, if $\gk<_T\gz<_T\xi$, then by maximality of $\gk$,
$\gz\not\leq{\ga\otimes\gth_{\ga}}$ hence $\gz\geq
(\ga\otimes\gth_{\ga})+1=(\ga+1)\otimes 0$, contradicting 
the minimality of $\xi$.
As $M_{\xi}^{\cal T}$ is a model of some
tree ${\cal T}_{\eta}$ with $\eta>\ga$,
while $M_{\gk}^{\cal T}$ appears in some
${\cal T}_{\gz}$ with $\gz\leq\ga$,
then it must be the case that $\xi={\ga'+1\otimes 0}$,
for some $\ga'>\ga$,
${\cal T}_{\ga'+1}\bot{\cal T}_{\ga'}$
and $\gk=\mbox{\bf B}(\xi)$. So $\gk={\gb'\otimes\gg'}$ and 
$P_{\ga'+1}=\ult(M_{\gg'}^{{\cal T}_{\gb'}},
E_{\gth_{\ga'}}^{{\cal T}_{\ga'}})$.
By quasi-linearity $\gk\notin\mbox{\bf F}_{\ga'}\supseteq
\mbox{\bf F}_{\ga+1}$.
But $\mbox{\bf F}_{\ga+1}\supseteq 
[{\gb\otimes\gg},{\ga\otimes\gth_{\ga}}]$,
so $\gk\in\mbox{\bf F}_{\ga+1}$: a contradiction.
\end{proof}

Given a quasi-linear decomposition 
$\seq{({\cal T}_{\ga}, P_{\ga})\mid\ga<\nu}$ of
$({\cal T}, M)$, let $A=\{\ga+1<\nu\mid {\cal T}_{\ga+1}
\bot {\cal T}_{\ga}\}$. 
For $\ga+1, \gb+1\in A$ with $\ga<\gb$, then
$\mbox{\bf B}\big(\ga+1,0\big)\in
\mbox{\bf F}_{\ga+1}\subseteq\mbox{\bf F}_{\gb}$
and $\mbox{\bf B}\big(\gb+1,0\big)\notin\mbox{\bf F}_{\gb}$,
hence $\mbox{\bf B}\big(\ga+1,0\big)\neq
\mbox{\bf B}\big(\gb+1,0\big)$. 
In other words, the
function $A\ni\ga+1\mapsto\mbox{\bf B}\big(\ga+1,0\big)$
is injective. We want to thin-down $A$ so that 
this function is also increasing.
Let $B=\{\mbox{\bf B}\big(\ga+1,0\big)\mid\ga+1\in A\}$
and define $f:B'\to A$, $B'\subseteq B$ be 
such that $f$ is increasing and 
$\mbox{\bf B}\big(f(\gb),0\big)=\gb$.
For $\gb\in B$ let
$$
f(\gb)= \left\{\begin{array}{ll}
\ga+1\in A &\quad \mbox{if }\gb=\mbox{\bf B}(\ga+1,0)\mbox{ and }
\big(\forall\gg\in B\cap\gb\big)\, \ga>f(\gg),\\
 & \\
\mbox{undefined} &\quad \mbox{otherwise.}
\end{array}\right.
$$
Let $B'=\mbox{dom}(f)$. 
Let $A'=\mbox{ran}(f)\subseteq A$ and let
$\seq{\ga_i+1\mid i<\gl}$ be the increasing enumeration
of $A'$, $\gl=\mbox{o.t.}(A')=\mbox{o.t.}(B')$.
Note that $\mbox{\bf B}\big(\ga_0+1,0\big)=
\min B'=\min B$ and $i<j<\gl$ implies
$\mbox{\bf B}\big(\ga_i+1,0\big)<
\mbox{\bf B}\big(\ga_j+1,0\big)$.
$\seq{\ga_i+1\mid i<\gl}$ is the {\em basic sequence\/}
of $\seq{({\cal T}_{\ga}, P_{\ga})\mid\ga<\nu}$.

Note that $\seq{\ga_i(\xi)+1\mid i<\gl(\xi)}$,
the basic sequence of $\seq{({\cal T}_{\ga},P_{\ga})\mid\ga<\xi}$,
$\xi<\nu$, is not, in general, an initial segment
of $\seq{\ga_i+1\mid i<\gl}$. In fact it is not even
true that $\gl(\xi)\leq\gl$.
Let us list a few facts whose proof is left to the reader.

\begin{enumerate}
\item
If $\xi\leq\eta+1<\nu$ implies 
${\cal T}_{\eta+1}\|{\cal T}_{\eta}$, 
then $\gl(\xi)=\gl$ and
$$
\seq{\ga_i(\xi)+1\mid i<\gl(\xi)}=
\seq{\ga_i+1\mid i<\gl}.
$$

\item
Suppose $\xi=\eta+1<\nu$ 
and ${\cal T}_{\eta+1}\bot {\cal T}_{\eta}$. Let
$$
\gl'= \left\{\begin{array}{ll}
\gl(\xi) &\quad \mbox{if }\mbox{\bf B}\big(\xi,0\big)\geq
\sup\{\mbox{\bf B}\big(\ga_i(\xi)+1,0\big)
\mid i<\gl(\xi)\},\\

{} & {} \\

j &\quad \mbox{otherwise, where }j<\gl(\xi)\mbox{ is least
such that }\\

{} &\quad \mbox{\bf B}\big(\xi,0\big)<
\mbox{\bf B}\big(\ga_j(\xi)+1,0\big).
\end{array}\right.
$$
Then $\gl(\xi+1)=\gl'+1$ and
$$
\seq{\ga_i(\xi+1)+1\mid i<\gl(\xi+1)}=
\seq{\ga_i(\xi)+1\mid i<\gl'}^{\frown}
\seq{\xi}.
$$

\item
If $\xi_n\to\xi$ and ${\cal T}_{\xi_n+1}\bot{\cal T}_{\xi_n}$, then
$$
\lim_{n\to\infty}\seq{\ga_i(\xi_n+1)+1\mid i<\gl(\xi_n+1)}
=\seq{\ga_i(\xi)+1\mid i<\gl(\xi)}
$$
meaning that $\gl(\xi_n+1)\to\gl(\xi)$ and for every 
$i<\gl(\xi)$, $\ga_i(\xi)=\lim_{n\to\infty}
\ga_i(\xi_n+1)$. 
If we choose the $\xi_n$'s to be in
$\{\ga_i(\xi)+1\mid i<\gl(\xi)\}$, then
$$
\seq{\ga_i(\xi)+1\mid i<\gl(\xi)}=
\bigcup_{n<\go}\seq{\ga_i(\xi_n+2)+1\mid i<\gl(\xi_n+2)}.
$$
\end{enumerate}

\begin{lemma}
\label{uniqueness}
Suppose $\seq{({\cal T}_{\ga},P_{\ga})\mid\ga<\nu}$ is 
a quasi-linear decomposition of $({\cal T}, M)$ and
let $\seq{\ga_i+1\mid i<\gl}$ be its basic sequence.

\begin{enumerate}
\item
For $i<\gl$ and $\gk<\gth=lh({\cal T})$, 
if $(\ga_i+1)\otimes 0<\gk$ then $(\ga_i+1)\otimes 0<_T\gk$.
In particular, $i<j$ implies $(\ga_i+1)\otimes 0<_T
(\ga_j+1)\otimes 0$.

\item
${\cal T}$ has exactly one cofinal branch $b$ and
$(\ga_i+1)\otimes 0\in b$, for every $i<\gl$. 

\end{enumerate}
\end{lemma}

\begin{proof}
(1). 
By induction on $\gk$, for fixed $i<\gl$.
Let $\gb=\mbox{\bf B}(\gk)$.
If $\gb=\gk$, then 
$\gk={\xi\otimes 0}$ with $\xi$ limit and for 
some increasing sequence $\xi_n\to\xi$, 
${\cal T}_{\xi_n+1}\bot {\cal T}_{\xi_n}$. 
By inductive hypothesis, for $n$ sufficiently large, 
$(\ga_i+1)\otimes 0<_T (\xi_n+1)\otimes 0$
and by quasi-linearity $(\xi_n+1)\otimes 0<_T
{\xi\otimes 0}=\gk$,
so by transitivity $(\ga_i+1)\otimes 0<_T\gk$.
Assume now that $\gb<\gk$.
By quasi-linearity, $\gb\notin
\mbox{\bf F}_{\ga_i}\supseteq [
\mbox{\bf B}\big(\ga_i+1,0\big),
{\ga_i\otimes\gth_{\ga_i}}]$,
so either $\gb<\mbox{\bf B}\big(\ga_i+1,0\big)$ 
or $\gb\geq{\ga_i\otimes\gth_{\ga_i}}+1=
(\ga_i+1)\otimes 0$.
By definition of basic sequence, 
$\gb<\mbox{\bf B}\big(\ga_i+1,0\big)$ is impossible,
so $\gb\geq(\ga_i+1)\otimes 0$.
As $\gb<_T\gk$, then, by inductive hypothesis, 
$(\ga_i+1)\otimes 0\leq_T\gb<_T\gk$.

\medskip\noindent
(2). If the $(\ga_i+1)\otimes 0$'s are bounded in $\gth=lh({\cal T})$, then
for $\sup_{i<\gl}\ga_i<\eta\leq\gk<\nu$, 
${\eta\otimes 0}<_T {\gk\otimes 0}$, so they yield a 
cofinal branch $b$. 
If, otherwise, the $\ga_i$'s are unbounded in
$\gth$, let $b=\{\ga<\gth\mid (\exists i<\gl)\,\ga<_T
(\ga_i+1)\otimes 0\}$. By part (1) this is a branch.
In both cases it is immediate to verify that
$b$ is the only cofinal branch of $\cal T$.
\end{proof}

Given a quasi-linear decomposition
$\seq{({\cal T}_{\ga},P_{\ga})\mid\ga<\nu}$
of $({\cal T},M)$ and $\gb<\nu$, it is not 
the case, in general, that
$\seq{({\cal T}_{\ga}, P_{\ga})\mid\gb\leq\ga<\nu}$
yields an iteration tree: it is true if either we 
are dealing with a linear decomposition,
i.e. $(\forall\ga+1<\nu)\: {\cal T}_{\ga+1}
\|{\cal T}_{\ga}$, or else if $\gb$ is chosen carefully.

\begin{definition}
Given $({\cal T},M)$ of length $\gth$, 
with quasi-linear decomposition
$\seq{({\cal T}_{\ga},P_{\ga})\mid\ga<\nu}$
and basic sequence
$\seq{\ga_i+1\mid i<\gl}$, and given
$\gb\leq\go_1$, let
$$
\gG(\gb)=\sup\{(\ga_i+1)\otimes 0\mid i<\gl\mbox{ and }
(\ga_i+1)\otimes 0\leq\gb\}.
$$
\end{definition}
{\bf Remarks.}
\begin{enumerate}

\item
If $\bar{\ga}=\sup\{\ga_i+1\mid i<\gl
\mbox{ and }(\ga_i+1)\otimes 0 \leq\gb\}$,
then $\gG(\gb)=\bar{\ga}\otimes 0$. 
In other words, $\gG(\gb)$ is always of 
the form $\ga\otimes 0$, for some $\ga$.
In particular, if 
$\seq{({\cal T}_{\ga},P_{\ga})\mid\ga<\gn}$
is linear, then $\gG(\gb)=0\otimes 0=0$, for all $\gb$.

\item
Note that $\gG(\gb)\leq\min(\gb,\gth)$ and 
it is defined even for $\gb\geq\gth$:
if $\gb\geq\ga'=\sup\{(\ga_i+1)\otimes 0\mid i<\gl\}$, then
$\gG(\gb)=\ga'$. 

\item
By part (1) of \ref{uniqueness}, $\ga\leq\gb$
implies $\gG(\ga)\leq_T\gG(\gb)\leq_T\gb$.

\item
We will write $\gG^{\xi}(\gb)$ when  the basic 
sequence $\seq{\ga_i(\xi)+1\mid i<\gl(\xi)}$,
$\xi\leq\nu$, is used.
\end{enumerate}

\begin{lemma}
Let $\seq{({\cal T}_{\ga},P_{\ga})\mid\ga<\nu}$ be
a quasi-linear decomposition of $({\cal T},M)$. 
For $\gb\leq\go_1$
$$
{\cal T}\lfloor \gG(\gb)=
\seq{(E^{\cal T}_{\ga},\gr_{\ga}^{\cal T})\mid
\gG(\gb)<\ga+1<\gth}
$$
is an iteration tree on 
$M_{\gG(\gb)}^{\cal T}=P_{\bar{\ga}}$, where
$\gG(\gb)={\bar{\ga}\otimes 0}$. Its tree ordering
is $<_T\re [\gG(\gb),\gth)$ and 
$\seq{({\cal T}_{\ga},P_{\ga})\mid\bar{\ga}\leq
\ga<\nu}$ is its quasi-linear decomposition.
\end{lemma}

\begin{proof}
We only need to check that, for a fixed $\gb$,
$\gG(\gb)<\gk$ implies $\gG(\gb)<_T\gk$.
By part (1) of Lemma \ref{uniqueness}, for any $i<\gl$
such that $(\ga_i+1)\otimes 0\leq\gb$,
$(\ga_i+1)\otimes 0<_T\gk$, so taking the supremum, we have that
$\gG(\gb)<_T\gk$.
\end{proof}

\noindent {\bf Remarks:}
\begin{enumerate}
\item
The careful reader might object that, formally,
${\cal T}\lfloor\gG(\gb)$ is not an iteration tree
as the field of its tree ordering is not an ordinal.
On the other hand, by a trivial re-indexing it
can be construed as an iteration tree of length 
$$
\gth-\gG(\gb)=\mbox{o.t.}\,\big[\gG(\gb),\gth\big).
$$
It is convenient though, for our purposes, to 
think of the models and embeddings 
of ${\cal T}\lfloor\gG(\gb)$ 
as indexed by ordinals in the interval
$\big[\gG(\gb),\gth\big)$.

\item
If $\gg$ is an ordinal such that 
$<_T\cap [\gg,\gth)\times [\gg,\gth)$
is a tree ordering on $\xi$, then 
${\cal T}\lfloor\gg$ makes sense.
In a way, $\gG(\gb)$ is the optimal 
such $\gg\leq\gb$: say that $\gg$ is
{\em minimal \/} if there is no 
$\ga<\gg$ such that for all $\gb$ $(\ga\leq\gb+1<\gg\implies
{\cal T}_{\gb+1}\|{\cal T}_{\gb})$. Then
$$
\gG(\gb)=\sup\{\gg\leq\gb\mid\gg\mbox{ is minimal and }
{\cal T}\lfloor\gg \mbox{ is defined }\}.
$$
\end{enumerate}

\bigskip\noindent
{\bf The game ${\cal G}(M,\go_1+1)$.}
This game lasts $\go_1+1$ rounds, with $I$ playing only at 
successor rounds and $II$ playing at every round,
in which $I$ and $II$ cooperatively construct
an iteration tree on $M$. 
There are several constraints on the moves that $I$ and $II$
are allowed: the first player to violate these constraints
loses. 
If $II$ has not lost by $\go_1$, then he wins. 
At round $\gn+1<\go_1$, 
a certain iteration tree $\cal T$ is
given and player $I$ 
plays an iteration tree of successor length
${\cal S}_{\gn}$ with the intent of extending
$\cal T$ in a quasi-linear way. $II$ can either accept 
it, or reject it and play a maximal wellfounded branch $b$ 
of ${\cal S}_{\gn}$: the tree ${\cal S}_{\gn}$
is truncated at $\gl=\sup (b)$ and then extended by $b$.
The resulting tree ${\cal T}_{\gn}$ is then used to extend
$\cal T$. At limit stages $II$ must play a cofinal wellfounded
branch $b$ of the current $\cal T$. Formally:

\bigskip\noindent
At round $\nu<\go_1$ a quasi-linear decomposition
$\seq{({\cal T}_{\ga},P_{\ga})\mid\ga<\nu}$ of
some $({\cal T},M)$ of length $\gth<\go_1$
has been played.
\begin{itemize}
\item
If $\nu$ is limit, it is $II$'s turn to move:
he must play a cofinal wellfounded branch $b$ of
$({\cal T},M)$. By Lemma \ref{uniqueness},
$II$ does not have much choice, so if the 
only cofinal branch $b$ is illfounded,
$II$ loses.
If $II$ does not lose at round $\nu$,
set $P_{\nu}=M_b^{\cal T}$ and ${\cal T}_{\nu}=$
the empty tree on $P_{\nu}$. Thus 
$\seq{({\cal T}_{\ga},P_{\ga})\mid\ga<\nu+1}$
is a quasi-linear decomposition of $\cal T$
extended by $b$.

\item
Suppose $\nu=\xi+1$. Player $I$ has two choices.

\begin{enumerate}
\item
$I$ plays an iteration tree $({\cal S}_{\nu}, P_{\nu})$,
where $P_{\nu}=M_{\gth}^{\cal T}$,
such that, extending $({\cal T},M)$ via ${\cal S}_{\nu}$,
we still have a non-overlapping tree.

\item
$I$ plays $({\cal S}_{\nu},P_{\nu},E_{\gth},\gb,\gg)$ such that 

\begin{enumerate}
	\item
	$M_{\gth}^{\cal T}\models\mbox{\lq\lq}E_{\gth}$
	is a $2^{\aleph_0}$-closed extender",
	$\gg\leq\gth_{\gb}$, $\gb\leq\xi$, 
	the least model $E_{\gth}$ can be applied to,
	is $M_{\gb\otimes\gg}^{\cal T}$ and if $\gb=\xi$
	then $\gg<\gth_{\xi}$,

	\item
	${\cal S}_{\nu}$ is an iteration tree 
	on $P_{\nu}=\ult\big(
	M_{\gb\otimes\gg}^{\cal T},E_{\gth}\big)$,

	\item 
	the tree resulting from extending $\cal T$ by 
	${\cal S}_{\nu}$ is still non-overlapping
	and 

        $\seq{({\cal T}_{\ga},P_{\ga})\mid\ga<\nu}^{\frown}
	\seq{({\cal S}_{\nu},P_{\nu})}$ satisfies the
	definition of quasi-linearity, except, possibly,
	for not having a last model.
\end{enumerate}
\end{enumerate}
Player $II$ responds by playing either
\begin{enumerate}
\item
({\bf accept}). Then set ${\cal T}_{\nu}={\cal S}_{\nu}$ and
extend $\cal T$ via ${\cal T}_{\nu}$ (and $E_{\gth}$,
if $I$ played as in case 2).
This move is legal only if ${\cal S}_{\nu}$ 
is of successor length, in which case we set 
$\gth_{\nu}=lh({\cal T}_{\nu})$.

\item
$(\mbox{\bf accept},b)$, where $b$ is a cofinal wellfounded
branch of ${\cal S}_{\nu}$. Let ${\cal T}_{\nu}$ be
${\cal S}_{\nu}$ extended by $b$.
Extend $\cal T$ via ${\cal T}_{\nu}$ (and $E_{\gth}$,
if $I$ played as in case 2).

\item
$(\mbox{\bf reject},b)$, where $b$ is a maximal 
wellfounded branch of ${\cal S}_{\nu}$, such that 
$\sup(b)<lh({\cal S}_{\nu})$, i.e. $b$ is non-cofinal. 
Let ${\cal T}_{\nu}$ be 
${\cal S}_{\nu}\re\sup(b)$ extended by $b$.
Extend $\cal T$ via ${\cal T}_{\nu}$ (and $E_{\gth}$, 
if $I$ played as in case 2).
\end{enumerate}
\end{itemize}
This concludes the definition of ${\cal G}(M,\go_1+1)$.

\begin{theorem}
\label{little}
Suppose $M$ is a countable premouse elementarily
embeddable in some $V_{\eta}$, $\pi:M\to V_{\eta}$.
Then $II$ has a winning strategy for ${\cal G}(M,\go_1+1)$.
\end{theorem}

\begin{proof}
Choose $\gk>\eta$ large enough so that
$V_{\gk}$ is a premouse with $\gd(V_{\gk})=\eta$ and with 
$\go_1+1$ cut-off points above $\eta$. 
We first construct premice 
$N_{\ga}$, for $1\leq\ga\leq\go_1$,
ordinals $\eta(\ga,\gb)$, $\eta(\ga,\infty)$,
and maps $\pi_{\ga,\gb}$ and $\pi_{\ga,\infty}$
for $0\leq\ga\leq\gb\leq\go_1$ such that 
\begin{enumerate}

\item
$N_0=M$ and, for $\ga>0$, each 
$N_{\ga}$ contains $HC$ and
$\seq{(N_{\gb},\eta(\gg,\gb),\pi_{\gg,\gb})
\mid 0\leq\gg\leq\gb<\ga}$, 
$\gd(N_{\ga})=\eta(0,\ga)$,
and for each $0<\gb<\ga$,
$N_{\ga}\models |N_{\gb}|=2^{\aleph_0}$.

\item
For $0\leq\ga<\gb\leq\gg\leq\go_1$,
$\eta(\ga,\gg)<\eta(\gb,\gg)$, 
$N_{\gg}\cap\On=\eta(\gg,\gg)$.
Moreover $\eta(0,\infty)=\eta$ and, for $\ga>0$,
$\eta(\ga,\infty)$ is the $\ga$th cut-off 
point of $V_{\gk}$ above $\eta$.

\item
For $\ga\leq\gb\leq\gg\leq\go_1$,
$\pi_{\ga,\gb}:N_{\ga}\to N_{\gb}\cap V_{\eta(\ga,\gb)}$
is an elementary embedding such that,
$\pi_{\ga,\gg}=\pi_{\gb,\gg}\circ\pi_{\ga,\gb}$ 
and $\pi_{\ga,\ga}=\mbox{id}\re N_{\ga}$. 

\item
For $\ga\leq\gb\leq\go_1$,
$\gp_{\ga,\infty}:N_{\ga}\to V_{\eta(\ga,\infty)}$
is an elementary embedding such that
$\pi_{\gb,\infty}\circ\pi_{\ga,\gb}=\pi_{\ga,\infty}$ 
and $\pi=\pi_{0,\infty}$.

\end{enumerate}
To see this suppose $N_{\gb}$, $\eta(\gg,\gb)$, 
$\pi_{\gg,\gb}$, $\pi_{\gg,\infty}$ have
been defined for $\gg\leq\gb<\ga$,
and let 
$$
H=\hull ^{V_{\eta(\ga,\infty)}}\Big(HC\cup\bigcup_{\gb<\ga}
\mbox{ran}(\gp_{\gb,\infty})\cup
\seq{\eta(\gb,\infty)\mid\gb<\ga}
\Big)
$$
let $\gp_{\ga,\infty}:N_{\ga}\to V_{\eta(\ga,\infty)}$ be the 
inverse of the transitive collapse and let 
$\pi_{\gb,\ga}=(\gp_{\ga,\infty})^{-1}\circ \gp_{\gb,\infty}$.

In order to show that $II$ has a winning strategy 
$\gS$ in ${\cal G}(M)$, we shall define a system of maps
$\gf_{\ga,\gb}^{\nu}$, $\gs_{\ga}^{\nu}$ and trees
${\cal U}_{\ga}^{\nu}$ at round $\nu$.
Suppose we are at a position of length $\nu$, 
according to $\gS$, and suppose we have built 
so far a tree $({\cal T},M)$ of length $\gth$ 
with quasi-linear decomposition 
$\seq{({\cal T}_{\ga},P_{\ga})\mid\ga<\nu}$
and basic sequence $\seq{\ga_i+1\mid i<\gl}$.
The maps $\gf^{\nu}_{\ga,\gb}$ should be 
thought of as being the \lq\lq stage $\nu$" 
versions of the $\pi_{\ga,\gb}$:
if ${\cal T}_{\xi +1}\|{\cal T}_{\xi}$
for all $\xi+1<\nu$, then 
$\gf_{\ga,\gb}^{\nu}=\pi_{\ga,\gb}$
and ${\cal U}_{\ga}^{\nu}=\pi_{\ga,\gb}{\cal T}$.
If, otherwise, there is $\xi+1<\nu$ such that
${\cal T}_{\xi+1}\bot{\cal T}_{\xi}$, and 
$\xi$ is least such, then the $\gf$'s change:
say that $P_{\xi+1}=\ult(M_{\gg}^{{\cal T}_{\gb}}, E)$
and let $\ga=\gb\otimes\gg$. Then 
$\gf^{\xi+1}_{0,\ga+1}:P_{\xi+1}\to N_{\ga+1}$.
On the other hand, the $\gs^{\nu}$'s 
guarantee that, if $\nu$ is limit and for 
all sufficiently large $\xi+1<\nu$ 
${\cal T}_{\xi+1}\|{\cal T}_{\xi}$, then 
the direct limit model is wellfounded.

Here is our official definition.
Suppose we are given:
\begin{enumerate}
\item 
elementary embeddings 
$\gf_{0,\ga}:M_{\gG(\ga)}^{\cal T}\to
N_{\ga}\cap V_{\eta(0,\ga)}$ such that, 
for $0\leq\ga\leq\gb\leq\go_1$
$$
\gf_{0,\gb}\circ i^{\cal T}_{\gG(\ga),\gG(\gb)}
=\pi_{\ga,\gb}\circ\gf_{0,\ga}.
$$
Note that $\gG(\ga)=\gG(\gb)$ with
$\ga<\gb$ is possible:
in this case $M_{\gG(\ga)}^{\cal T}=M_{\gG(\gb)}^{\cal T}$
and $\gf_{\ga,\gb}=\pi_{\ga,\gb}\circ\gf_{0,\ga}$.

\item
${\cal U}_{\ga}=\gf_{0,\ga}\big({\cal T}\lfloor\gG(\ga)\big)$,
the iteration tree ${\cal T}\lfloor\gG(\ga)$ copied
on $N_{\ga}$ via $\gf_{0,\ga}$.
So $\gf_{0,\ga}:({\cal T}\lfloor \gG(\ga), M_{\gG(\ga)}^{\cal T})
\to ({\cal U}_{\ga}, N_{\ga})$ is a bounded embedding with bound
$\eta(0,\ga)$.
Following our convention above, 
all models and embeddings
of ${\cal U}_{\ga}$ are indexed by ordinals in 
$[\gG(\ga),\gth)$.

\item
For $0\leq\ga\leq\gb\leq\go_1$, let
$\gf_{\ga,\gb}:M_{\gG(\gb)}^{{\cal U}_{\ga}}
\to N_{\gb}\cap V_{\eta(\ga,\gb)}$ 
be an elementary map such that $\gf_{\ga,\gb}\circ
i_{\gG(\ga),\gG(\gb)}^{{\cal U}_{\ga}}=\pi_{\ga,\gb}$.
[Recall that $M_{\gG(\ga)}^{{\cal U}_{\ga}}$ is 
the base model of ${\cal U}_{\ga}$, that is $N_{\ga}$.]
Moreover, letting $\gf^{\ast}:M_{\gG(\gb)}^{\cal T}
\to M_{\gG(\gb)}^{{\cal U}_{\ga}}$ be the copy map
induced by $\gf_{0,\ga}$, then 
$\gf_{\ga,\gb}\circ\gf^{\ast}=\gf_{0,\gb}$.
$$
\begin{diagram}
\node{M_{\gG(\gb)}^{\cal T}}
\arrow[3]{se,t}{\gf_{0,\gb}}
\arrow{sse,r}{\gf^{\ast}} \\

\\

\node{M_{\gG(\ga)}^{\cal T}}
\arrow[2]{n,l}{i_{\gG(\ga),\gG(\gb)}^{\cal T}}
\arrow{se,b}{\gf_{0,\ga}} 

\node{M_{\gG(\gb)}^{{\cal U}_{\ga}}}
\arrow{ese,t}{\gf_{\ga,\gb}} \\

\node{M}
\arrow{n,l}{}
\arrow{e,b}{}

\node{N_{\ga}=M_{\gG(\ga)}^{{\cal U}_{\ga}}}
\arrow{n,r}{i_{\gG(\ga),\gG(\gb)}^{{\cal U}_{\ga}}}
\arrow[2]{e,b}{\gp_{\ga,\gb}}

\node[2]{N_{\gb}}

\end{diagram}
$$

\item
For $\sup_{i<\gl}\ga_i+1\leq\ga\leq\gb<\nu$ 
the following is a commutative diagram of elementary maps
$$
\begin{diagram}
   \node{M_{\gb\otimes\gth_{\gb}}^{{\cal U}_{\go_1}}}
   \arrow[2]{se,t}{\gs_{\gb+1}}\\
   \node{M_{\ga\otimes 0}^{{\cal U}_{\go_1}}}
   \arrow{n,l}{} \arrow{ese,t}{\gs_{\ga}}\\
   \node{N_{\go_1}=M_{\gG(\go_1)}^{{\cal U}_{\go_1}}}
   \arrow{n,l}{} \arrow[2]{e,b}{\gp_{\go_1,\infty}}
   \node[2]{V_{\eta(\go_1,\infty)}}
\end{diagram}
$$
In particular, $\gs_{\nu}$ is defined when
$\nu=\gb+1$.
Note that, for $\gb\geq\sup_{i<\gl}\ga_i+1$,\,\ 
${\cal T}_{\gb+1}\|{\cal T}_{\gb}$, so
$\gb\otimes\gth_{\gb}=(\gb+1)\otimes 0$.
\end{enumerate}

As the objects above change, as the play goes on, 
we should really write:
$$
{\cal U}_{\ga}^{\nu},\qquad
\gf^{\nu}_{\ga,\gb},\qquad
\gs^{\nu}_{\ga},\qquad
\seq{\ga_i(\nu)+1\mid i<\gl(\nu)}\qquad
\mbox{and}\qquad
\gG^{\nu}(\gb).
$$
For the sake of legibility, we will drop the 
superscript $\nu$ whenever possible.
The objects defined at round $\nu$ are related 
to the ones defined at round $\xi<\nu$ as follows.

\begin{itemize}
\item
Suppose $\exists\xi_0<\nu$ such that for all
$\xi_0<\xi+1<\nu$, ${\cal T}_{\xi+1}\|{\cal T}_{\xi}$.
[In this case $\nu$ could be limit or successor.]
Then, for such $\xi_0$ and $\xi$,
$\gl(\nu)=\gl(\xi_0)=\gl(\xi)$, and the 
basic sequences are the same,
$\seq{\ga_i(\nu)+1\mid i<\gl(\nu)}=
\seq{\ga_i(\xi)+1\mid i<\gl(\xi)}$.
Also $\gf_{\ga,\gb}^{\nu}=\gf_{\ga,\gb}^{\xi}$ and 
${\cal U}^{\nu}_{\ga}$ extends ${\cal U}^{\xi}_{\ga}$ 
via $\gf_{0,\ga}^{\nu}{\cal T}_{\xi}$, 
for $0\leq\ga\leq\gb\leq\go_1$.
Finally, if $\nu$ is limit, then
$\gs^{\nu}_{\ga}=\gs^{\xi}_{\ga}$ for all $\ga$ 
such that
$\sup\{\ga_i(\nu)+1\mid i<\gl(\nu)\}\leq\ga
<\xi+1\leq\nu$,
and if $\nu=\xi+1$, then $\gs_{\ga}^{\nu}=\gs_{\ga}^{\xi}$
for all $\ga\leq\xi$.

\item
Suppose $\nu$ is limit and there is an increasing 
sequence $\xi_n\to\nu$ such that 
${\cal T}_{\xi_n+1}\bot{\cal T}_{\xi_n}$.
Without loss of generality we may assume that
each $\xi_n+1\in\{\ga_i(\nu)+1\mid i<\gl(\nu)\}$.
Then, for $\gg\leq\nu$, $\gb\leq\go_1$ 
and ${\nu\otimes 0}\leq\ga\leq\go_1$
$$
\seq{\ga_i(\nu)+1\mid i<\gl(\nu)}=
\bigcup_{n<\go}\seq{\ga_i(\xi_n+2)+1\mid i<\gl(\xi_n+2)},
$$
$$
\gG^{\nu}(\gb)=\sup_n\gG^{\xi_n+2}(\gb),\qquad
\gG^{\nu}(\ga)={\nu\otimes 0} \qquad
\mbox{and}\qquad\gs_{\gg}^{\nu}=\pi_{\go_1,\infty}.
$$
[The \lq\lq +2" in the indices above is because we want 
to consider quasi-linear systems where ${\cal T}_{\xi_n+1}$
is the last tree.]
Note that for $\gb<{\nu\otimes 0}$ the sequence $\gG^{\xi_n+2}(\gb)$
is eventually $=\gG^{\nu}(\gb)$. 
Thus, for $n$ larger than some fixed $m$,
$$
\gf_{\ga,\gb}^{\nu}=\gf_{\ga,\gb}^{\xi_n+2}
\qquad\mbox{ and }\qquad
{\cal U}_{\gb}^{\nu}=\gf^{\nu}_{0,\gb}
({\cal T}\lfloor\gG^{\nu}(\gb))=
\bigcup_{m<n}{\cal U}_{\gb}^{\xi_n+2}.
$$
If, instead, $\gb\geq{\nu\otimes 0}$, then the ordinals
$\gG^{\xi_n+2}(\gb)$ are stricly increasing and
converge to $\gG^{\nu}(\gb)={\nu\otimes 0}$. In this case,
for $\ga<\gb=\nu\otimes 0$, 
$\gf_{\ga,\gb}^{\nu}$ is the direct limit map 
of the commutative system of embeddings
$$
\begin{diagram}
\node{M_{\gG^{\nu}(\gb)}^{\cal T}}
\arrow{ese,t}{}\\

\node{M_{\gG^{\xi_n+2}(\gb)}^{\cal T}}
\arrow{n,l,..}{}
\arrow{ese,t}{}

\node[2]{M_{\gG^{\nu}(\gb)}^{{\cal U}^{\nu}_{\ga}}}
\arrow[2]{se,t}{\gf^{\nu}_{\ga,\gb}}\\

\node{M_{\gG^{\xi_n+2}(\ga)}^{\cal T}}
\arrow{n,l}{}
\arrow{ese,t}{\gf_{0,\ga}^{\xi_n+2}}

\node[2]{M_{\gG^{\xi_n+2}(\gb)}^{{\cal U}^{\xi_n+2}_{\ga}}}
\arrow{n,l,..}{}
\arrow{ese,b}{\gf^{\xi_n+2}_{\ga,\gb}} \\

\node{M}
\arrow{n,l}{}
\arrow[2]{e,b}{\pi_{0,\ga}}

\node[2]{N_{\ga}}
\arrow{n,l}{}
\arrow[2]{e,b}{\gp_{\ga,\gb}}

\node[2]{N_{\gb}}

\end{diagram}
$$
For 
$\nu\otimes 0\leq\gb\leq\gg$, 
$\gf_{\gb,\gg}^{\nu}=\pi_{\gb,\gg}$
as ${\cal U}^{\nu}_{\gb}$ and ${\cal U}^{\nu}_{\gg}$ 
are the empty iteration trees
on $N_{\gb}$ and $N_{\gg}$, respectively.
Finally, if $\ga<{\nu\otimes 0}\leq\gb$,
we set 
$\gf_{\ga,\gb}^{\nu}=\pi_{\nu\otimes 0,\gb}\circ
\gf_{\ga,\nu\otimes 0}^{\nu}$.

\item
Suppose $\nu=\xi+1$ and $\xi$ is limit.
By definition of the game $\cal G$, 
the quasi-linear decomposition 
$\seq{({\cal T}_{\ga},P_{\ga})\mid\ga<\xi+1}$
of $({\cal T},M)$ is such that $P_{\xi}$ 
is the direct limit of the $P_{\xi_n}$'s,
for some cofinal sequence $\xi_n\to\xi$,
and ${\cal T}_{\xi}$ is the empty tree on $P_{\xi}$.
Then 
$$
\seq{\ga_i(\nu)+1\mid i<\gl(\nu)}=
\seq{\ga_i(\xi)+1\mid i<\gl(\xi)}
$$
and for $0\leq\ga\leq\gb\leq\go_1$ and $\gg<\xi$
$$
\gG^{\nu}(\gb)=\gG^{\xi}(\gb),\quad
\gf^{\nu}_{\ga,\gb}=\gf^{\xi}_{\ga,\gb},\quad
\gs_{\gg}^{\gn}=\gs_{\gg}^{\xi},\quad
{\cal U}_{\gb}^{\nu}\mbox{ extends }
{\cal U}_{\gb}^{\xi}
$$
and 
$\gs_{\xi}^{\nu}:M_{\xi\otimes 0}^{{\cal U}_{\go_1}}\to V_{\gk}$ 
is the direct limit map induced by the embeddings
$\gs_{\gg}^{\nu}$.
\end{itemize}

\begin{claim}
Suppose that $\nu\leq\go_1$ 
is limit and that $II$ has not 
lost by round $\nu$, and that
there are $\gf^{\nu}$'s, 
${\cal U}^{\nu}$'s and 
$\gs^{\nu}$'s as above, then
$II$  does not lose at round $\nu$. 
In particular if $II$ has not lost by
round $\go_1$, then the $\gf^{\go_1}$'s, 
${\cal U}^{\go_1}$'s and $\gs^{\go_1}$'s 
witness that he wins ${\cal G}(M,\go_1+1)$.
\end{claim}

\begin{proof}
As $\nu$ is limit, it is $II$'s turn to move:
he has to verify that  the only cofinal branch $b$ 
of the iteration tree $\cal T$ built so far, is wellfounded.
By hypothesis we are given $\gf_{\ga,\gb}=\gf_{\ga,\gb}^{\nu}$,
${\cal U}_{\ga}={\cal U}_{\ga}^{\nu}$ and
$\gs_{\ga}=\gs_{\ga}^{\nu}$ satisfying the 
conditions above.

\medskip

\noindent {\bf Case 1:} 
$(\exists\xi_0<\nu)(\forall\xi)\xi_0<\xi+1<\nu\implies
{\cal T}_{\xi+1}\|{\cal T}_{\xi}$.

\smallskip

\noindent
Choosing $\xi_0$ least such, $\gG(\go_1)=\xi_0\otimes 0$
and $\gf_{0,\go_1}$ copies ${\cal T}\lfloor\gG(\go_1)$
to ${\cal U}_{\go_1}$. Let 
$\xi_n\to\nu$ be increasing and let 
$$
f_n:M_{\xi_n\otimes 0}^{\cal T}\to
M_{\xi_n\otimes 0}^{{\cal U}_{\go_1}}
$$
be the copy map induced by $\gf_{0,\go_1}$.
Then $\gs_{\xi_n}\circ f_n$ witness 
$M_b^{\cal T}$ is wellfounded.

\medskip

\noindent{\bf Case 2:} There is an increasing sequence 
$\xi_n\to\nu$ such that
${\cal T}_{\xi_n +1}\bot{\cal T}_{\xi_n}$.

\smallskip

\noindent
Without loss of generality we can assume that 
for all $n<\go$, $\ga_n=\gG(\xi_n +1)$ belongs to
the basic sequence, hence to $b$.
Let 
$$
f_n=\pi_{\xi_n+1,\go_1}\circ\gf_{0,\xi_n+1}:
M_{\ga_n}^{\cal T}\to N_{\go_1}.
$$
Thus, for $n<m$, $f_n=f_m\circ i_{\ga_n,\ga_m}$
so that there is a limit map 
$f_{\go}:M_b^ {\cal T}\to N_{\go_1}$
witnessing wellfoundedness.
\end{proof}
So it is enough to show that $\gf^{\nu}$'s,
${\cal U}^{\nu}$'s and $\gs^{\nu}$'s as above 
exist for every $\nu\leq\go_1$.

\bigskip

If $\nu=0$, then $\gf_{\ga,\gb}=\pi_{\ga,\gb}$,
${\cal U}_{\ga}$ is the empty tree on $N_{\ga}$ and 
$\gs_0=\pi_{\go_1,\infty}:N_{\go_1}\to V_{\gk}$.

\bigskip

If $\nu$ is limit, then the remarks before the claim 
show how to define the
$\gf^{\nu}$, ${\cal U}^{\nu}$ and $\gs^{\nu}$ 
from the 
$\gf^{\xi}$, ${\cal U}^{\xi}$ and $\gs^{\xi}$,
for $\xi<\nu$. 

\bigskip

We now take care of round $\nu+1$.
Suppose we are given $\gf^{\nu}$'s, $\gs^{\nu}$'s 
and ${\cal U}^{\nu}$'s, we must describe how $II$ answers
to $I$'s moves and how to build $\gf^{\nu+1}$'s, $\gs^{\nu+1}$'s
and ${\cal U}^{\nu+1}$'s. In order to avoid making
our formul\ae\ exceedingly ornate, we shall drop the
suffix $\nu$ for the objects at stage $\nu$ and use
$\gG'$, $\psi_{\ga,\gb}$, $\gt_{\ga}$ and ${\cal V}_{\ga}$ for 
$\gG^{\nu+1}$, $\gf^{\nu+1}_{\ga,\gb}$, $\gs_{\ga}^{\nu+1}$ 
and ${\cal U}_{\ga}^{\nu+1}$.

\bigskip\noindent
$\bullet$ {\bf Suppose $I$ plays $({\cal S}_{\nu},P_{\nu})$.}

\noindent
Then $P_{\nu}=M_{\gth}^{\cal T}$ and $II$ tries to copy 
${\cal S}_{\nu}$ on $V_{\gk}$. If he succeeds 
to do so and finds a branch $b$, then $II$ 
plays $(\mbox{\bf accept},b)$. Otherwise
he rejects ${\cal S}_{\nu}$.
To be more precise.

\noindent
As $\gth+1=lh({\cal T})$, there is
$\gs_{\nu}:M_{\gth}^{{\cal U}_{\go_1}}\to 
V_{\eta(\go_1,\infty)}\subset V_{\gk}$ 
and let 
$\gf:M_{\gth}^{\cal T}\to 
M_{\gth}^{{\cal U}_{\go_1}}$ 
be the copy map induced by 
$\gf_{0,\go_1}$.
Let $\gth'\leq lh({\cal S}_{\nu})$
be largest such that ${\cal S}_{\nu}\re\gth'$
can be copied on $V_{\gk}$
via $\gs_{\nu}\circ\gf$
and there are no wellfounded maximal 
branches cofinal in some $\gg<\gth'$.
If $\gth'=lh({\cal S}_{\nu})$ is a successor
ordinal, let $II$ play $(\mbox{\bf accept})$
and set ${\cal T}_{\nu}={\cal S}_{\nu}$ 
and $\gth_{\nu}=\gth'-1$.
By Lemma \ref{successor}, if $\gth'< lh({\cal S}_{\nu})$,
then it must be a limit ordinal, so we may assume that
$\gth'\leq lh({\cal S}_{\nu})$ is limit. 
By $2^{\aleph_0}$-closure of the extenders
and Theorem 5.6 in \cite{Iterationtrees}.
$\gs_{\nu}\circ\gf ({\cal S}_{\nu}\re\gth')$
must have a cofinal wellfounded branch $b$. 
Let $II$ play $(\mbox{\bf accept},b)$,
if $\gth_{\nu}=lh({\cal S}_{\nu})$, 
or $(\mbox{\bf reject},b)$, if 
$\gth_{\nu}<lh({\cal S}_{\nu})$, and in either case
let ${\cal T}_{\nu}$ be ${\cal S}_{\nu}\re\gth_{\nu}$
extended by $b$. 

In order to keep the induction going, 
the $\psi$'s, $\gt$'s and $\cal V$'s must be defined
and shown to satisfy the inductive hypothesis.
The tree $\cal T$ extended by ${\cal T}_{\nu}$ will 
be denoted by $\cal S$. 
As ${\cal T}_{\nu}\|{\cal T}_{\xi}$, then
$\seq{({\cal T}_{\ga},P_{\ga})\mid\ga<\nu+1}$ 
is  quasi-linear, 
its basic sequence is still $\seq{\ga_i+1\mid i<\gl}$,
hence $\gG'=\gG$.
Set $\psi_{\ga,\gb}=\gf_{\ga,\gb}$ and
extend ${\cal U}_{\ga}$ to ${\cal V}_{\ga}$, 
by tagging (an isomorphic copy of) ${\cal T}_{\nu}$
on top: ${\cal V}_{\ga}=\psi_{0,\ga}\big(
{\cal S}\lfloor\gG(\ga)\big)$.
For $\xi\leq\nu$, $\gt_{\xi}$ can be taken to be
$\gs_{\xi}$, so we are only left to define 
$\gt_{\nu+1}:M_{\nu\otimes\gth_{\nu}}^{{\cal V}_{\go_1}}
\to V_{\eta(\go_1,\infty)}$, which will be obtained
using Lemma \ref{chunks}.
[The argument will be presented in some detail and will serve
as a template for other proofs in this paper.]

First observe that  ${\cal V}_{\go_1}\lfloor\gth=
\psi_{0,\go_1}({\cal S}\lfloor\gth)$ is (isomorphic to)
${\cal T}_{\nu}$ and, by construction,
${\cal V}_{\go_1}$ can be copied on $V_{\gk}$ 
via $\gt_{\nu}$. 
Let us denote by $W_{\ga}$ and $Z_{\ga}$
the $\ga$th models of the tree 
${\cal T}_{\nu}$ copied on $V_{\gk}$
and on $V_{\eta(\go_1,\infty)}$, respectively,
via $\gt\circ\psi_{0,\go_1}$. 
In other words, $Z_{\ga}=W_{\ga}\cap V_{\eta'}$
where $\eta'={i_{0,\ga}
(\eta(\go_1,\infty))}$.
Let $\gt^{\ast}:M_{\nu\otimes\gth_{\nu}}^{{\cal V}_{\go_1}}
\to Z_{\gth_{\nu}}\subset W_{\gth_{\nu}}$ be 
the copy map induced by $\gt_{\nu}$.
The commutative diagram below may 
help to follow the argument.
(The horizontal arrows come from the copy
construction, while the vertical ones are the tree
embeddings.)
$$
\begin{diagram}
\node{M_{\nu\otimes \gth_{\nu}}^{\cal S}}
\arrow{e,t}{}
\node{M_{\nu\otimes\gth_{\nu}}^{{\cal V}_{\go_1}}}
\arrow{se,t}{\gt^{\ast}} \\

\node{M_{\gth}^{\cal S}}
\arrow{n,l}{}
\arrow{e,t}{}
\node{M_{\gth}^{{\cal V}_{\go_1}}}
\arrow{se,t}{\gt_{\nu}}
\arrow{n,l}{}
\node{Z_{\gth_{\nu}}} 
\arrow{e,t}{\mbox{id}}
\node{W_{\gth_{\nu}}} \\

\node{M_{\gG(\go_1)}^{\cal S}}
\arrow{n,l}{}
\arrow{e,b}{\psi_{0,\go_1}}
\node{N_{\go_1}}
\arrow{n,l}{}
\arrow{e,b}{\gp_{\go_1,\infty}}
\node{V_{\eta(\go_1,\infty)}}
\arrow{n,r}{}
\arrow{e,b}{\mbox{id}}
\node{V_{\gk}}
\arrow{n,l}{}
\end{diagram}
$$
Note that $|N_{\go_1}|=2^{\aleph_0}$, ${\cal S}$ is  countable and
$2^{\aleph_0}$-closed, so the hypotheses of Lemma
\ref{chunks} hold. Thus $\gt^{\ast}$ can be taken 
to be an element of $W_{\gth_{\nu}}$, and
by elementarity of the tree-embedding from $V_{\gk}$ 
to $W_{\gth_n}$, there is $\gt_{\nu+1}\in V_{\gk}$ such that 
$$
\begin{diagram}
\node{M^{\cal S}_{\nu\otimes\gth_{\nu}}}
\arrow{e,t}{}
   \node{M_{\gn\otimes\gth_{\gn}}^{{\cal V}_{\go_1}}}
   \arrow{sse,t}{\gt_{\gn+1}}\\
\node{M_{\gth}^{\cal S}}
\arrow{n,l}{}
\arrow{e,t}{}
   \node{M_{\gth}^{{\cal V}_{\go_1}}}
   \arrow{n,l}{} \arrow{se,b}{\gt_{\nu}}\\
\node{M_{\gG(\go_1)}^{\cal S}}
\arrow{n,l}{}
\arrow{e,t}{}
   \node{N_{\go_1}}
   \arrow{n,l}{} \arrow{e,b}{\gp_{\go_1,\infty}}
   \node{V_{\eta(\go_1,\infty)}}
\end{diagram}
$$
commutes, and this is what we had to prove.

\bigskip
\noindent$\bullet$ {\bf Suppose $I$ plays $({\cal S}_{\nu},P_{\nu},
E_{\gth},\gb_0,\gg_0)$.}

\noindent
Let $\ga=\gb_0\otimes\gg_0$. As $I$'s move is 
legal, then $\ga$ cannot belong to $\mbox{\bf F}_{\nu}$,
the set of forbidden nodes at stage $\nu$.
The next Claim is crucial for the present construction.

\begin{claim}
$\gG(\ga+1)\leq\ga$.
\end{claim}

\begin{proof}
Deny. As $\gG(\ga+1)\leq\ga+1$, it follows that
$\ga+1=\gG(\ga+1)=(\ga_i+1)\otimes 0$, for some 
$i<\gl(\nu)$. The definition of basic sequence implies that
${\cal T}_{\ga_i+1}\bot{\cal T}_{\ga_i}$,
and $\ga=\ga_i\otimes\gth_{\ga_i}$.
Thus $\ga\in\mbox{\bf F}_{\ga_i+1}\subseteq
\mbox{\bf F}_{\nu}$: a contradiction.
\end{proof}
As $\gG(\ga+1)\leq\ga$, the embedding
$\gf_{\ga,\ga+1}: M_{\gG(\ga+1)}^{{\cal U}_{\ga}}\to
N_{\ga+1}\cap V_{\eta(\ga,\ga+1)}$
copies to bounded embeddings
$$
\gf:M_{\ga}^{{\cal U}_{\ga}}\to 
M_{\ga}^{{\cal U}_{\ga+1}}
\qquad\mbox{and}\qquad
\gf^{\ast}:M_{\gth}^{{\cal U}_{\ga}}\to M_{\gth}^{{\cal U}_{\ga+1}}.
$$
Let also 
$\gf':M_{\gth}^{\cal T}\to M_{\gth}^{{\cal U}_{\ga}}$
be the copy map induced by 
$\gf_{0,\ga}:M_{\gG(\ga)}^{\cal T}\to N_{\ga}$. 
By Lemma \ref{successor}, 
$$
P_{\nu}=\ult (M_{\ga}^{\cal T},E_{\gth}), \qquad  
P=\ult (M_{\ga}^{{\cal U}_{\ga}},\gf'(E_{\gth})), 
\qquad\mbox{and}\qquad
Q=\ult(M_{\ga}^{{\cal U}_{\ga+1}},\gf^{\ast}\circ\gf' (E_{\gth}))
$$
are wellfounded, so the copy constructions 
yield bounded embeddings $\psi^{\ast}:P\to Q$ and 
$\psi':P_{\nu}\to P$.
By Lemma \ref{chunks}, 
we get 
$\psi_{\ga,\ga+1}:P \to N_{\ga+1}\cap V_{\eta(\ga,\ga+1)}$.
Set
$$
\psi_{\gb,\gg}=\left\{
\begin{array}{ll}
\gf_{\gb,\gg}&\quad \mbox{ if } \gb\leq\gg\leq\ga,\\

\pi_{\ga,\gb}&\quad \mbox{ if } \ga+1\leq\gb\leq\gg, \\

\pi_{\ga+1,\gg}\circ\psi_{\ga,\ga+1}\circ\gf_{\gb,\ga}&
\quad \mbox{ if } \gb\leq\ga<\gg,
\end{array}\right.
$$
and $\gt_{\xi}=\pi_{\go_1,\infty}$, for all $\xi\leq\nu$.
Now the argument proceeds as before.
Let $\gth'\leq lh({\cal S}_{\nu})$
be largest such that ${\cal S}_{\nu}\re\gth'$
can be copied on $V_{\gk}$
via $\gt_{\nu}\circ\psi_{0,\go_1}$
and there are no wellfounded maximal
branches cofinal in some $\gg<\gth'$.
If $\gth'=lh({\cal S}_{\nu})$ is a successor
ordinal, let $II$ play $(\mbox{\bf accept})$
and set ${\cal T}_{\nu}={\cal S}_{\nu}$.
If, otherwise, $\gth'\leq lh({\cal S}_{\nu})$
is limit, then let $b$ be a wellfounded
branch of 
$\gt_{\nu}\circ\psi_{0,\ga}({\cal S}_{\nu}\re\gth')$
and let $II$ play $(\mbox{\bf accept},b)$, or
$(\mbox{\bf reject},b)$, depending on
whether $\gth'=lh({\cal S}_{\nu})$ or not.
In either case set ${\cal T}_{\nu}={\cal S}_{\nu}$.
By Lemma  \ref{chunks} there is $\gt_{\nu+1}$ such that
the diagram
$$
\begin{diagram}
\setlength{\dgARROWLENGTH}{2.8em}
\node{M_b^{{\cal S}_{\nu}}}
\arrow[2]{e,t}{}

\node[2]{M_b^{\psi'{\cal S}_{\nu}}}
\arrow{se,t}{} 
\\

\node{P_{\nu}}
\arrow[2]{e,t}{\psi'}
\arrow{n,l}{}

\node[2]{P}
\arrow{se,t}{}
\arrow{n,l}{}

\node{M_b^{\psi_{0,\ga+1}{\cal S}_{\nu}}}
\arrow[2]{e,t}{}

\node[2]{M_b^{\psi_{0,\go_1}{\cal S}_{\nu}}}
\arrow{se,t}{\gt_{\nu+1}}
\\

\node{M_{\gG^{\nu}(\ga)}}
\arrow[2]{e,b}{\gf_{0,\ga}=\psi_{0,\ga}}
\arrow{n,l}{}

\node[2]{N_{\ga}}
\arrow{n,l}{}
\arrow{e,b}{\pi_{\ga,\ga+1}}

\node{N_{\ga+1}}
\arrow{n,l}{}
\arrow[2]{e,b}{\pi_{\ga+1,\go_1}=\psi_{\ga+1,\go_1}}

\node[2]{N_{\go_1}}
\arrow{e,b}{\pi_{\go_1,\infty}}
\arrow{n,l}{}

\node{V_{\gk}}
\end{diagram}
$$
commutes. Let
$$
\gG'(\gb)=\left\{
\begin{array}{ll}
\gG(\gb)&\qquad\mbox{for }\gb\leq\ga=\gb_0\otimes\gg_0,\\
 & \\
\ga+1&\qquad\mbox{for }\ga<\gb,
\end{array}
\right.
$$
and ${\cal V}_{\gb}=\psi_{0,\gb}({\cal S}\lfloor\gG'(\gb))$.
We leave to the reader the verification that the 
$\gf^{\nu+1}$'s, $\gs^{\nu+1}$'s, ${\cal U}^{\nu+1}$'s so
defined satisfy the inductive hypothesis.

\bigskip
As we have taken care of all possible cases, the Theorem is proved.
\end{proof}

%% file: hownotlose.tex
%

Let $M$ be a countable coarse premouse.
We will consider the following iteration game on $M$, 
called ${\cal G}^+(M)$.
It is played like the ordinary full iteration game, 
with $II$ on the move at limits and $I$ on 
the move at the other rounds, except that
$I$ must also play (besides the extenders $E_{\ga}$'s)
distinct natural numbers $n_{\ga}$, with
$n_{\ga}\notin\{n_{\gb}\mid\gb<\ga\}$.
The game is over when $I$ runs out of integers.
We also require that the iteration tree 
that $I$ and $II$ construct is
$2^{\aleph_0}$-closed.

The length of the game thus depends on the play,
but is always $<\go_1$. On the other hand, it
is easy to see that ${\cal G}^+(M)$
is stronger than the full iteration game (for
$2^{\aleph_0}$-closed extenders) of 
{\em fixed\/} countable length.

The rest of this paper is devoted to a proof of

\begin{theorem}
\label{main}
Player $I$ does not have a winning strategy for ${\cal G}^+(M)$,
for a countable premouse $M$ elementarily embeddable
in some $V_{\eta}$.
\end{theorem}
If $\gS$ is a strategy for $I$ then 
we have a continuous coding of ${\cal G}^+(M)$ in the sense
of \cite{Longgames}, i.e. a function {\bf c} form the set of all legal
positions of ${\cal G}^+(M)$ to $\go$ such that if $p$ and $q$ 
are positions and $q$ extends $p$, then 
$\hbox{\bf c}(p)\neq \hbox{\bf  c}(q)$: just take
$\mbox{\bf c}(p)=$ the natural number $n$ 
given by $\gS$ at position $p$.
A position $p$ of ${\cal G}^+(M)$ is, essentially, an iteration tree
of successor length. 
To avoid confusions, we denote with $lh (p)$ the length of $p$
as a position, and with $ht(p)$ the length (or height) of the associated
iteration tree, so that $ht(p)=lh(p)+1$.

\noindent For the reader's benefit, here is a brief 
description the plan of the proof.
In all previous iterability theorems (see 
\cite{Iterationtrees},\cite{Wfdd}), one argues
by contradiction: given a \lq\lq bad" 
tree, i.e. a tree without cofinal wellfounded branches, 
using the ordinals witnessing such 
badness, an infinite descending $\in$-chain is
constructed.
This contradiction forces us 
to conclude that such bad tree cannot exist
in the first place, hence 
the theorem follows. To be more specific.
Ordinals are assigned to the models of the bad tree
$\cal T$, witnessing continuous illfoundedness.
These ordinals are then used to build Skolem hulls
of the models of the bad tree.
Exercising proper care in the construction of such
hulls, it can be shown that they form an inconsistent
enlargement, that is a system of models resembling
the original iteration tree, but containing an infinite
descending $\in$-chain.

Back to our proof.
If we try to argue by contradiction following the pattern 
above, we are immediately faced with the problem that
we are not given a bad tree, but rather a bad strategy
$\gS$, i.e. a winning strategy for $I$, hence we cannot 
{\em first\/} fix a tree and {\em then\/} get the ordinals
for the construction. In other words, the ordinals
should be given \lq\lq continuously in the tree".
The cure for this is to construct positions
$p_1\subset\ldots\subset p_n\ldots$ of the game ${\cal G}^+(M)$
together with enlargements 
${\cal P}^n=\seq{P_{\ga}\mid\ga<
\gth_n}^{\frown}\seq{P^{\ast}_n}$
of $p_n$ (here $\gth_n=lh(p_n)$) such that 
$P^{\ast}_{n+1}\in P^{\ast}_n$, obtaining thus a contradiction.
The enlargement ${\cal P}^n$ will be 
constructed by taking hulls
of the models of the (pseudo-)iteration tree obtained
by copying the position $p_n$ on ${\cal P}^{n-1}$.
The ordinals needed for this construction are ranks of certain
nodes on a wellfounded tree $\cal U$ on some $V_{\gk}$,
searching for a defeat of $\gS$.
Let's take a closer look at $\cal U$.

Fix towards a contradiction, $\gp:M\to V_{\eta}$ and a 
winning strategy $\gS$ for $I$ in ${\cal G}^+(M)$. 
Let $\gk>\eta$ be such that $V_{\gk}$ is a premouse
with $\gd(V_{\gk})=\eta$.
For any position $p$, let
$$
I_p=\big\{n\in\go \mid \forall\gb\leq lh(p)
(n\neq \hbox{\bf c}(p\re \gb))\big\}
$$
that is, the set of all natural numbers not yet
played by position $p$.
$\cal U$ searches, among other things, for a sequence of positions 
according to $I$'s strategy $\gS$,
$p_1\subset p_2\subset\ldots\subset\bigcup_n p_n=\cal T$, a 
cofinal wellfounded branch $b$ of $\cal T$ and
families of premice ${\cal C}^1\subset\ldots\subset{\cal C}^n$,
with $lh({\cal C}^n)=lh(p_n)$. 
The position $p_n$ can be copied, as a pseudo-iteration tree, 
on ${\cal C}^{n-1\,\frown}\seq{V_{\gk}}$. 
Also if $m$ is the least element of $I_{p_n}$,
we make sure that either $m\notin I_{p_{n+1}}$,
or there is no $q\in\cal U$ extending $p_n$ according to $\gS$
that can be copied onto $V_{\gk}$ and that $m\notin I_q$.
We then add to the node of the tree $\cal U$ 
a family of models ${\cal B}^n$ of length $=ht(p_n)$, 
each model of size at most $2^{\aleph_0}$,
and resembling enough to 
${\cal C}^{n\,\frown}\seq{V_{\gk}}$, 
and witnessing that no such $q$ can be copied on ${\cal B}^n$.

Suppose $\cal U$ had a branch, namely positions
$p_1\subset p_2\subset\ldots\subset\bigcup_n p_n=\cal T$, 
and a cofinal wellfounded branch $b$ of $\cal T$. 
This would determine a new position $p=({\cal T}, b)$, and
let {\bf c}$(p)=m$. As the integer $m$ was considered
at some stage $n$, while choosing position $p_{n+1}$, 
and as $m\in I_{p_{n+1}}$, it follows that no position
$q$ with $m\notin I_q$ could have been copied
onto ${\cal B}^n$. But, as it turns out, $p$ is either
such a position or a defeat for $\gS$.
This shows that $\cal U$ is wellfounded.

The proof now proceeds as follows.
Positions $p_n$ and families of premice ${\cal P}^n$
are built inductively so that:

\begin{enumerate}

\item 
the $p_n$'s are according to $\gS$, extend each other and 
$\seq{p_n\mid n<\go}$ is a complete play of ${\cal G}^+(M)$;

\item
${\cal P}^n=\seq{P_{\ga}\mid\ga<\gth_n}^{\frown}\seq{P_n^{\ast}}$, 
${\cal P}^n\re\gth_n={\cal P}^{n+1}\re\gth_n$,
where $\gth_n +1=ht(p_n)$, and
$p_n$ together with ${\cal C^n}={\cal P}^n\re\gth_n$ 
make up for part of a node of $\cal U$ of length $n+1$;

\item 
$P^{\ast}_{n+1}\in P^{\ast}_n$.

\end{enumerate}
The idea is to choose $p_{n+1}$ first, and then copy 
the iteration tree on the current enlargement ${\cal P}^n$.
For any $\gth_n\leq\ga<\gth_{n+1}$, the $\ga$th model 
on the pseudo-iteration tree on ${\cal P}^n$ is replaced
by the transitive collapse of a Skolem hull. 
Call these models $P_{\ga}$. By a tree
argument the $P_{\ga}$'s can be taken to be in the last model 
(the one with index $\gth_{n+1}$) of the pseudo-iteration 
tree on ${\cal P}^n$. By taking another hull and calling it
$P_{n+1}^{\ast}$ we obtain ${\cal P}^{n+1}$. 
The ordinal needed to take the hull in the $\ga$th model
is the rank of the node given by $p_{n+1}$, ${\cal P}^n$
and an initial segment of a cofinal branch passing through $\ga$.

This concludes our brief description of the structure of the
proof.

%% file: bigtree.tex
%

From andretta@math.ucla.edu Mon Oct  4 18:52:28 1993
Received: from julia.math.ucla.edu by pianeta (4.1/SMI-4.1)
	id AA15440; Mon, 4 Oct 93 18:52:05 +0100
Received: from sonia.math.ucla.edu by julia.math.ucla.edu 
	(Sendmail 4.1/1.07) id AA09075; Mon, 4 Oct 93 10:55:52 PDT
Return-Path: <andretta@math.ucla.edu>
Received: by sonia.math.ucla.edu 
	(Sendmail 4.1/1.07) id AA07559; Mon, 4 Oct 93 10:55:51 PDT
Date: Mon, 4 Oct 93 10:55:51 PDT
From: Alessandro Andretta <andretta@math.ucla.edu>
Message-Id: <9310041755.AA07559@sonia.math.ucla.edu>
To: andretta@di.unito.it
Subject: bigtree.tex
Status: RO

Let's first introduce some handy notation.
Suppose $({\cal T},{\cal B})$ is a pseudo-iteration tree
of length $(\gth+1,\gl)$ and suppose that
${\cal C}=\seq{C_{\ga}\mid\ga\leq\gth}$ is a family
of premice and $\Pi=\seq{\pi_{\ga}\mid\ga\leq\gth}$
$$
\Pi:({\cal T},{\cal B})\to (\emptyset,{\cal C})
$$
is a bounded embedding. If $\cal C$ is internal, we call
the pair $({\cal C},\Pi)$ an {\em enlargement\/}.

\begin{definition}
\label{definitionchunk}
Let $p$ be a position of ${\cal G}^+(M)$,
$\gth=ht(p)$, let $\ga_0,\ldots\ga_k,<\gth$
and let $\Pi:(p,M)\to (\Pi p,\cal B)$ be 
an embedding of pseudo-iteration trees.
\begin{itemize}

\item
$e(p):\gth\to\go$ is the 1--1 function 
defined by $e(p)(\gb)=\mbox{\bf c}(p\re\gb)$.

\item
Let $X$ be the support for $(p,M)$ generated by
$\{\ga_0,\ldots,\ga_k\}\cup e(p)^{-1\prime\prime} (k+1)$.
$$
\chunk{M_{\ga_k}^{(\Pi p,\cal B)};\ga_0,\ldots,\ga_k}
=\big( M_{\ga_k}^{(\Pi p,\cal B)}\big)_{\Pi X}.
$$
\end{itemize}
\end{definition}
{\bf Remark.}
If $q$ extends $p$, then $e(q)\supseteq e(p)$ so that,
if $\ga_0,\ldots,\ga_k<ht(p)$ 
$$
\chunk{M_{\ga_k}^{(\Pi p,\cal B)};\ga_0,\ldots,\ga_k}
\prec
\chunk{M_{\ga_k}^{(\Pi q,\cal B)};\ga_0,\ldots,\ga_k}.
$$

\bigskip
We are now ready to define the tree $\cal U$,
a set of finite sequences from $V_{\gk}$ closed under
initial segment.
A node \br, of length $n+1$, of $\cal U$ is of 
the form
$$
\br=\langle(p_1,{\cal B}^0,\Pi^0,\nu_0,\gs_0),
(p_2,{\cal B}^1,\Pi^1,\Pi^{0,1},{\cal C}^1,
\mbox{H}^1,\gF^1,\Psi^1,\nu_1,\gs_1),\ldots
$$
$$
\ldots,(p_{n+1},{\cal B}^n,\Pi^n,\Pi^{n-1,n},
{\cal C}^n,\mbox{H}^n,\gF^n,
\Psi^n,\nu_n,\gs_n)\rangle
$$
such that the following 6 clauses must hold.
\begin{enumerate}

\item
$p_1\subset p_2\subset\ldots\subset p_{n+1}$ 
are non-empty positions 
of the game ${\cal G}^+(M)$, according to $\gS$. 
Let $\gth_i+1=ht(p_i)$ and, for notational convenience,
let's agree that $p_0=\emptyset=$
the empty position, hence $I_{p_0}=\go$ and 
$\gth_0=lh(p_0)=0$.

\item \begin{enumerate}
      \item
      ${\cal C}^1\subset\ldots\subset{\cal C}^n$,
      $$
      {\cal C}^n=\seq{C_{\ga}\mid\ga<\gth_n}
      $$
      and, for each $\ga<\gth_n$, $C_{\ga}$ is a premouse and 
      $H_{(2^{\aleph_0})^+}\in C_{\ga}$.

      \item
      $\gF^1\subset\ldots\subset\gF^n$, 
      $\mbox{H}^1\subset\ldots\subset\mbox{H}^n$,
      $$
      \gF^n=\seq{\gf_{\ga}\mid\ga<\gth_n}\qquad
      \mbox{H}^n=\seq{\eta_{\ga}\mid\ga<\gth_n}
      $$ 
      such that, for each $\ga<\gth_n$,
      $\gf_{\ga}:M_{\ga}^{p_n}\to C_{\ga}\cap 
      V_{\eta_{\ga}}$ is an elementary embedding.

      \item
      There is an elementary embedding 
      $\gf^{\ast}_n:M_{\gth_n}^{p_n}
      \to V_{\eta}$ such that
      $$
      ({}^{\ast}{\cal C}^n,{}^{\ast}\gF^n)=
      ({\cal C}^{n\,\frown}\seq{V_{\gk}},\gF^{n\,\frown}
      \seq{\gf^{\ast}_n})
      $$
      is an enlargement for $(p_n, M)$ with bounds
      $\seq{\eta_{\ga}\mid\ga<\gth_n}^{\frown}
      \seq{\eta}$.
\end{enumerate}

Let's agree to define $({}^{\ast}{\cal C}^0,{}^{\ast}\gF^0)$
to be the pair $(V_{\gk},\pi)$.
There is a further constraint on the sequence of the positions $p_i$'s.

\begin{claim}
For all $n\geq 0$ there is an $k\in I_{p_n}$ and a 
$q\supset p_n$ according to $\gS$ 
such that $k\notin I_q$ and $q$ can 
be copied on $({}^{\ast}{\cal C}^n,{}^{\ast}\gF^n)$.
\end{claim}

\begin{proof}
$I_{p_n}\neq\emptyset$ as {\bf c}$(p_{n+1})$ belongs
to this set. By Corollary \ref{corsuccessor},  any $q$ extending 
$p_n$ of height $\gth_n+2=ht(p_n)+1$ can be copied on 
$({}^{\ast}{\cal C}^n,{}^{\ast}\gF^n)$,
so we can take $(q,k)$ to be the result 
of $\gS$ applied to $p_n$, hence 
$k=\mbox{\bf c}(q)\in I_{p_n}$.
\end{proof}

Define $N(p_n)$ to be the least integer 
$k$ as in the Claim above.
Thus $N(p_0)$ is the least $k\in\go$ such that
there is a non-empty position $q$ that can be copied
on $V_{\gk}$ via $\pi$ and such that $k=\mbox{\bf c}(q\re\gb)$
for some $\gb\leq lh(q)$.

\item
For all $n\geq 0$, $N(p_n)\notin I_{p_{n+1}}$ 
and $p_{n+1}$ can be copied on 
$({}^{\ast}{\cal C}^n,{}^{\ast}\gF^n)$.

\item
${\cal B}^0\subseteq\ldots\subseteq {\cal B}^n$,
$\Pi^0\subseteq\ldots\subseteq\Pi^n$ and
\begin{enumerate}
      \item
      ${\cal B}^n=\seq{B_{\ga}\mid\ga\leq\gth_m}$,
      for some $m\leq n$, and for each $\ga\leq\gth_m$,
      $B_{\ga}$ is a premouse of cardinality 
      $2^{\aleph_0}$ containing $HC$.
      Thus ${\cal B}^0=\seq{B_0}$ is a single premouse.
      
      \item
      For $n>0$, $N(p_n)\neq\min I_{p_n}\iff {\cal B}^{n-1}
      \neq{\cal B}^n$, and, if this is the case,
      $lh({\cal B}^n)=\gth_n+1$.

      \item
      $\Pi^n=\seq{\pi_{\ga}\mid\ga< lh({\cal B}^n)}$
      and there are ordinals $\ge_{\ga}\in B_{\ga}$ such that
      $\pi_{\ga}:M_{\ga}^{p_n}\to B_{\ga}\cap 
      V_{\ge_{\ga}}$ is an elementary embedding.
      The iteration tree $p_{n+1}$ can be copied on 
      $({\cal B}^n,\Pi^n)$ hence
      $\Pi^n:(p_{n+1},M)\to(\Pi^n p_{n+1},{\cal B}^n)$
      is a bounded embedding with bounds 
      $\seq{\ge_{\ga}\mid\ga<lh({\cal B}^n)}$.

      \item
      For $n>0$, 
      $\Pi^{n-1,n}=\seq{\pi^{n-1,n}_{\ga}\mid\ga\leq\gth_{n+1}}$,
      where the
      $$
      \pi^{n-1,n}_{\ga}:M_{\ga}^{(\Pi^{n-1}p_{n+1},{\cal B}^{n-1})}
      \to M_{\ga}^{(\Pi^n p_{n+1},{\cal B}^n)}
      $$
      are elementary, $\pi^{n-1,n}_{\ga}=\mbox{id}\re B_{\ga}$,
      for $\ga<lh({\cal B}^{n-1})$, and the diagram
      $$
      \setlength{\dgARROWLENGTH}{2.1em}
      \begin{diagram}
      \node{(p_{n+1},M)}
      \arrow[2]{e,t}{\gP^n} \arrow{se,b}{\gP^{n-1}}  
      \node[2]{(\Pi^n p_{n+1},{\cal B}^n)} \\
      \node[2]{ (\Pi^{n-1} p_{n+1},{\cal B}^{n-1})} 
      \arrow{ne,r}{\Pi^{n-1,n}}
      \end{diagram}
      $$
      commutes.

      \item
      For $n\geq 0$,
      $$
      \Psi^n=\seq{\psi^n_{\ga}\mid\ga<\gth_n}
      $$
      $\psi^n_{\ga}:M_{\ga}^{(\Pi^n p_n,{\cal B}^n)}
      \to C_{\ga}$ are elementary embeddings, and there is a
      $ \psi^{\ast}_n:M_{\gth_n}^{(\Pi^n p_n,
      {\cal B}^n)} \to V_{\gk}$
      such that
      $$
      {}^{\ast}\Psi^n=\Psi^{n\frown}
      \seq{\psi^{\ast}_n}:
      (\Pi^n p_n,{\cal B}^n)\to {}^{\ast}{\cal C}^n
      $$
      is an elementary embedding of pseudo-iteration trees. Moreover
      $$
      \setlength{\dgARROWLENGTH}{2.1em}
      \begin{diagram}
      \node{(p_n,M)}
      \arrow[2]{e,t}{{}^{\ast}\gF^n} \arrow{se,b}{\gP^n}  
      \node[2]{{}^{\ast}{\cal C}^n} \\
      \node[2]{(\Pi^n p_n,{\cal B}^n)} 
      \arrow{ne,r}{{}^{\ast}\Psi^n}

      \end{diagram}
      $$
      commutes and for $\ga<lh({\cal B}^n)$,
      $\psi^n_{\ga}(\ge_{\ga})=\eta_{\ga}$.

      \item 
      For $n\geq 0$, if $N(p_n)>\min(I_{p_n})$, then for every
      $k\in I_{p_n}\cap N(p_n)$ and every position
      $q\supset p_n$ according to $\gS$ with $k\notin I_q$,
      $q$ cannot be copied on ${\cal B}^n$. 
\end{enumerate}
Finally we take care of the $\nu$'s and $\gs$'s.

\item
$\nu_0<_T\ldots<_T\nu_n$, 
where $<_T$ is the tree ordering of the largest
iteration tree, $p_{n+1}$, and 
$\gth_n\leq\nu_n\leq\gth_{n+1}$.

\item
$\gs_0,\ldots,\gs_n$ are elementary embeddings 
$\gs_n: \chunk{M_{\nu_n}^{(\Pi^n p_{n+1},{\cal B}^n)};
\nu_0,\ldots,\nu_n}\to V_{\gk}$ such that,
for $n>0$, the diagram 
$$
\setlength{\dgARROWLENGTH}{2.5em}
\begin{diagram}
\node{\chunk{M_{\nu_{n-1}}^{(\Pi^n p_{n+1},
{\cal B}^n)};\nu_0,\ldots,\nu_{n-1}}} 
\arrow[3]{e,t}{\hat{\imath}}
\node[3]{\chunk{M_{\nu_n}^{(\Pi^n p_{n+1},{\cal B}^n)};
\nu_0,\ldots,\nu_n}} \arrow{s,r}{\gs_n}\\
\node{\chunk{M_{\nu_{n-1}}^{(\Pi^{n-1} p_n,{\cal B}^{n-1})};
\nu_0,\ldots,\nu_{n-1}}} 
\arrow[3]{e,b}{\gs_{n-1}} \arrow{n,l}{\hat{\pi}}
\node[3]{V_{\gk}} 
\end{diagram}
$$
commutes, where $\hat{\imath}$ and $\hat{\pi}$ are
the restrictions of 
$i_{\nu_{n-1},\nu_n}^{(\Pi^n p_{n+1},{\cal B}^n)}$ and 
$\pi^{n-1,n}_{\nu_{n-1}}$. Moreover 
$\pi=\gs_0\circ\pi^0_{\nu_0}\circ i_{0,\gn_0}^{p_1}$, 
where $\pi:M\to V_{\eta}$ is as in the hypothesis of our theorem.
In other words
$$
\setlength{\dgARROWLENGTH}{2.5em}
\begin{diagram}
\node{\chunk{M_{\nu_0}^{p_1};\nu_0}}
\arrow[2]{e,t}{\bar{\pi}}
\node[2]{\chunk{M_{\nu_0}^{(\Pi^0 p_1,{\cal B}^0)};\nu_0}}
\arrow{s,r}{\gs_0} \\
\node{M}
\arrow{n,l}{i_{0,\nu_0}^{p_1}}
\arrow[2]{e,b}{\pi}
\node[2]{V_{\gk}}
\end{diagram}
$$
commutes, where $\bar{\pi}$ is the restriction of $\pi_{\nu_0}^0$
to $\chunk{M_{\nu_0}^{p_1};\nu_0}$. 

\end{enumerate}

This concludes the definition of $\cal U$.
\bigskip

\noindent {\bf Remarks.}
The definition above has several awkward 
features that might appear unduly arbitrary.
The only reason we have chosen this particular
definition of $\cal U$, rather than more 
natural ones, is that it will greatly
simplify the construction in the next section.
The remarks that follow should help the reader
to understand some of the motivations behind the definition above.

\begin{enumerate}

\item 
The definition of $\cal U$ involves 
$M$, $\gS$, $\pi$, $\eta$ and $\gk$ as parameters.
We will relativize $\cal U$ to several models, all of which
contain $HC$, hence only $\pi$, $\eta$, $\gk$ will have to be
changed. 

\item
Every model of the pseudo-iteration tree 
$(\Pi^n p_n,{\cal B}^n)$ is of size $2^{\aleph_0}$, 
hence belongs to every $C_{\ga}$.

\item
For $0\leq n\leq m$ the elementary embedding 
$\Pi^{n,m}:(\Pi^n p_{m+1},{\cal B}^n)\to (\Pi^m p_{m+1},{\cal B}^m)$
mentioned in clause (6) is defined by induction:
$\Pi^{n,n}=$ the identity embedding, $\Pi^{n,m+1}=
\Pi^{m,m+1}\circ\Pi^{n,m}$

\item
Clauses (4.d) and (4.e) can be stated more concisely as
$$
\Pi^n=\Pi^{n-1,n}\circ\Pi^{n-1}\quad\mbox{and}\quad
{}^{\ast}\gF^n={}^{\ast}\Psi^n\circ\Pi^n
$$
(4.e) implies that for all $q\supset p_n$, 
if $q$ can be copied on ${}^{\ast}{\cal C}^n$ via ${}^{\ast}\Phi^n$
then $q$  can be copied on ${\cal B}^n$ via $\Pi^n$.
Using this and (4.d), every $p_m$ can be copied on 
any $({\cal B}^k,\Pi^k)$.

\item
Clause (4.f) says that $({\cal B}^n,\Pi^n)$ 
and $({}^{\ast}{\cal C}^n,{}^{\ast}\Phi^n)$
agree on the value of $N(p_n)$. To be more specific:
$N(p_n)=$ the least $k\in I_{p_n}$ such that there is a
$q\supset p_n$ according to $\gS$, such that
$k\notin I_q$ and $q$ can be copied on 
$({\cal B}^n,\Pi^n)$.
Clause (3) and (4.f) also imply that 
$N(p_{n+1})>N(p_n)$ and $I_{p_{n+1}}\subset I_{p_n}$.

\item
The reason why the $\gs_n$'s are defined on a chunks will 
become evident in the proof of Lemma \ref{B} in the next section.
The idea is that the $\gs_n$'s will be obtained from the
copying construction, and Lemma \ref{chunks} will 
be used.

\item
A few words  on the commutative diagrams of
clause (6). The embedding 
$$
i_{\nu_{n-1},\nu_n}^{(\Pi^n p_{n+1},{\cal B}^n)}:
M_{\nu_{n-1}}^{(\Pi^n p_{n+1},{\cal B}^n)} \to
M_{\nu_n}^{(\Pi^n p_{n+1},{\cal B}^n)} 
$$ 
is well defined as $\nu_{n-1}$ precedes 
$\nu_n$ in the tree ordering of $(\Pi^n p_{n+1},{\cal B}^n)$
and, by Definition \ref{definitionchunk}, it maps 
$\chunk{M_{\nu_{n-1}}^{(\Pi^n p_{n+1},{\cal B}^n)}; 
\nu_0,\ldots,\nu_n}$ elementarily into
$\chunk{M_{\nu_n}^{(\Pi^n p_{n+1},{\cal B}^n)}; 
\nu_0,\ldots,\nu_n}$, and 
$$
\chunk{M_{\nu_{n-1}}^{(\Pi^n p_{n+1},{\cal B}^n)};  
\nu_0,\ldots,\nu_{n-1}}\prec 
\chunk{M_{\nu_{n-1}}^{(\Pi^n p_{n+1},{\cal B}^n)};  
\nu_0,\ldots,\nu_n}.
$$
Thus $\hat{\imath}$ is well defined.
Similarly, as
$$
\chunk{M_{\nu_{n-1}}^{(\Pi^n p_n,{\cal B}^n)};  
\nu_0,\ldots,\nu_{n-1}}\prec
\chunk{M_{\nu_{n-1}}^{(\Pi^n p_{n+1},{\cal B}^n)};  
\nu_0,\ldots,\nu_{n-1}},
$$ 
$\hat{\gp}$ is well defined. 
Regarding the second commutative square, notice that
$0$ belongs to any support hence
$M=M_0^{p_1}=\chunk{M_0^{p_1};\nu_0}$.
\end{enumerate}

\begin{lemma}
$\cal U$ is wellfounded.
\end{lemma}

\begin{proof}
Deny. A branch of $\cal U$ is, essentially, a sequence
$$
\seq{(p_n,{\cal B}^n,\mbox{H}^n,\Pi^n,\Pi^{n,n+1},
{\cal C}^n,\gF^n,\Psi^n,\nu_n,\gs_n)
\mid n<\go}
$$
Then ${\cal T}=\bigcup_{n<\go}p_n$ 
is a countable iteration tree 
according to $\gS$, $\gth=\sup\{\gth_n+1\mid n<\go\}=lh(\cal T)$
and $b=\{\gb<\gth\mid\exists n(\gb<_T \gn_n)\}$
a cofinal branch of $\cal T$. 

\begin{claim}
$b$ is wellfounded branch of $\cal T$ and
letting ${\cal T}^+$ be the extension of 
$\cal T$ via $b$, ${\cal T}^+$ can be copied on 
$({\cal B}^n,\Pi^n)$, for all $n$.
\end{claim}

\begin{proof}
Fix $0\leq n<\go$. By remark (4) $\cal T$ can be copied
on $({\cal B}^n,\Pi^n)$. Also for $m\geq n$, the ordinals
$\nu_m$ are linearly ordered in the tree ordering
of $\Pi^n{\cal T}$ and determine a
cofinal branch of $\Pi^n {\cal T}$. By a minor abuse of 
notation, such a branch will still be denoted by $b$.
The direct limit of $(\Pi^n{\cal T},{\cal B}^n)$ 
along $b$ will be shown to be wellfounded, 
proving thus the claim.
By clause (6), for $m\geq n$, the diagram
$$
\setlength{\dgARROWLENGTH}{3.0em}
\begin{diagram}[\chunk{M_{\nu_m}^{(\Pi^n{\cal T},
{\cal B}^n)};\nu_0,\ldots,\nu_m}]
\node{M_b^{(\Pi^n{\cal T},{\cal B}^n)}}
\arrow{sse,t,1}{\gt_{\infty}}\\

\node{\chunk{M_{\nu_m}^{(\Pi^n p_{m+1},{\cal B}^n)};
\nu_0,\ldots\nu_m }}
\arrow{n,l,..}{i_{\nu_m,b}} 
\arrow{e,b,1}{\pi^{n,m}_{\nu_m}}
\arrow{se,b}{\gt_m} 
\node{\chunk{M_{\nu_m}^{(\Pi^m p_{m+1},{\cal B}^m)};
\nu_0,\ldots,\nu_m }}
\arrow{s,r}{\gs_m}\\

\node{\chunk{M_{\nu_n}^{(\Pi^n p_{n+1},{\cal B}^n)};
\nu_0,\ldots,\nu_n }}
\arrow{n,l}{i_{\nu_n,\nu_m}} \arrow{e,b}{\gs_n}
\node{V_{\gk}}
\end{diagram}
$$
commutes, where $\gt_m=\gs_m\circ\pi^{n,m}_{\nu_m}$, 
and $\gt_{\infty}$ is the limit map:
as every element $x\in M_b$ is of the form $i_{\gn_m,b}(y)$,
pick $k\geq m$ large enough so that 
$i_{\nu_m,\nu_k}(y)\in\chunk{M_{\gn_k}^{(\Pi^n p_{k+1},{\cal B}^n)};
\nu_0,\ldots,\nu_k}$,
and set $\gt_{\infty}(x)=\gt_k(i_{\nu_m,\nu_k}(y))$.

\end{proof}

Note that the game cannot be over once all the $p_n$'s have been
played, i.e. $\bigcap_n I_{p_n}\neq\emptyset$, 
as otherwise $II$ would win, $M_b$ being wellfounded, thus 
contradicting the assumption that $\gS$ 
is a winning strategy for $I$.

Thus extending $\cal T$ via $b$ yields a legal move of ${\cal G}^+(M)$,
call it $p$, that, by the claim, can be copied on any ${\cal B}^n$. 
Let $m=\hbox{\bf c}(p)$ and let $0\leq i\leq m+1$ be least such that 
$N(p_i)>m$. As $m\in\bigcap_n I_{p_n}$, then in particular
$m\in I_{p_i}$, hence $N(p_i)>m\geq\min I_{p_i}$.
Thus, by (3) in the definition of $\cal U$,
no extension $q$ of $p_i$ according to 
$\gS$ with $m\notin I_q$ could have 
been copied onto $({}^{\ast}{\cal C}^i,{}^{\ast}\gF^i)$, 
hence on $({\cal B}^i,\Pi^i)$.
In particular this should hold of $p$, but
the Claim shows that $p$ can be copied on ${\cal B}^i$:
a contradiction.
\end{proof}

%% file: enlargement.tex
%

From andretta@math.ucla.edu Mon Oct  4 18:53:13 1993
Received: from julia.math.ucla.edu by pianeta (4.1/SMI-4.1)
	id AA15446; Mon, 4 Oct 93 18:52:37 +0100
Received: from sonia.math.ucla.edu by julia.math.ucla.edu 
	(Sendmail 4.1/1.07) id AA09092; Mon, 4 Oct 93 10:56:22 PDT
Return-Path: <andretta@math.ucla.edu>
Received: by sonia.math.ucla.edu 
	(Sendmail 4.1/1.07) id AA07569; Mon, 4 Oct 93 10:56:21 PDT
Date: Mon, 4 Oct 93 10:56:21 PDT
From: Alessandro Andretta <andretta@math.ucla.edu>
Message-Id: <9310041756.AA07569@sonia.math.ucla.edu>
To: andretta@di.unito.it
Subject: enlargement.tex
Status: RO

Let $\br\in\cal U$ be a node:  
\brr\ is the finite sequence
obtained from \br\ by dropping the $\ga$'s and the $\gs$'s
and let ${\cal U}^-=\{\brr\mid\br\in{\cal U}\}$.
${\cal U}^-$ is still a tree on $V_{\gk}$, but
it is not wellfounded. 
In fact we will construct a branch through it.
Note that ${\cal U}^-$ can also be defined using 
clauses (1)---(4), without mentioning $\cal U$ at all.
Obviously, different $\br$'s may yield the same $\brr$,
so there is no way to retrieve e.g. the $\nu$'s from $\brr$.
Yet, for what we are going to do, we 
would like to be able to do this. To be more specific.
Suppose 
$\br\in{\cal U}^-$ is of length $n+1$ and that
$p_1\subset\ldots\subset p_{n+1}$ are its positions, 
and that $ht(p_n)\leq\ga<ht(p_{n+1})$:
can we find $\nu_0<\ldots<\nu_n=\ga$ so that
for some $\gs_0,\ldots,\gs_n$, 
$(\br,\vec{\nu},\vec{\gs})\in{\cal U}$?
The answer is no, as $\ga$ might not have $n+1$ 
predecessors in the $p_{n+1}$ tree ordering. 
In fact $<_T\mbox{-pred}(\ga)$ could be 0.
The next two definition address to this problem.

\begin{definition}
Suppose we are given positions 
$p_1\subset\ldots\subset p_{n+1}$, 
with $\gth_i+1=ht(p_i)$.
For every $\ga$ such that
$\gth_n<\ga\leq\gth_{n+1}$, if $n>0$,
or $0\leq\ga\leq\gth_1$, if $n=0$, 
the {\em backward sequence
of $\ga$ relative to\/} $p_1,\ldots ,p_{n+1}$, 
is the sequence 
$\seq{(\ga_0,m_1),\ldots,(\ga_k,m_{k+1})}$ defined as follows.

If $n=0$, then $\seq{(\ga_0,m_1)}=\seq{(\ga,1)}$ is the backward sequence
of $\ga$ relative to $p_1$.

Suppose $n>0$, then $k>0$. Let $<_T$ be the tree 
ordering associated to $p_{n+1}$.

\begin{enumerate}
\item
$0\leq\ga_0<_T\ldots<_T\ga_k=\ga$.

\item
$1=m_1<\ldots<m_{k+1}=n+1$, (hence $k\leq n$).

\item 
$\ga_{k-1}$ is the largest $\gb<_T\ga$ such that
$\gb\leq\gth_m$ for some $m<n+1$.
The least such $m$ is $m_k$.

\item
$\seq{(\ga_0, m_1),\ldots,(\ga_{k-1},m_k)}$ is the backward
sequence of $\ga_{k-1}$ relative to $p_1,\ldots,p_{m_k}$.

\end{enumerate}
\end{definition}

\noindent
{\bf Remarks.}

\begin{enumerate}

\item
Note that a largest $\gb$ as in (3) always exists as
the iteration trees $p_i$ have successor length $\gth_i +1$.

\item
For $0<i\leq k$, $\ga_i\in (\gth_{(m_{i+1})-1}, 
\gth_{m_{i+1}}]$.

\item
If $\bar{p}_i=p_{m_i}$,
the backward sequence of $\ga$ relative to
$\bar{p}_1,\ldots,\bar{p}_{k+1}$
is $\seq{(\ga_0,1),\ldots,(\ga_k, k+1)}$.
Notice that if $n>0$, then $k>0$.
\end{enumerate}

\begin{definition}
Let $\br\in{\cal U}^-$ be a node of length $n+1$
and let $p_1\subset\ldots\subset p_{n+1}$ be its positions.
Given $\gth_n\leq\ga<\gth_{n+1}$, let
$\seq{(\ga_0,m_1),\ldots,(\ga_k,m_{k+1})}$ be the
backward sequence of $\ga$ relative to
$\seq{p_1,\ldots,p_{n+1}}$.

The {\em contraction of\/} {\bf r} {\em relative to $\ga$\/}
is the node $\hbox{\bf t}\in\cal U^-$ of length
$k+1$ obtained by rearranging some of the stuff
contained in $\br$ in $k+1$ pieces:
$$
\begin{array}{rcl}
\mbox{\bf t}&=&\langle(p_{m_1},{\cal B}^0,\Pi^0),
\ldots (p_{m_{k+1}},{\cal B}^{m_k},\Pi^{m_k},\Pi^{m_{k-1},m_k},
{\cal C}^{m_k},\mbox{H}^{m_k}, \gF^{m_k},\Psi^{m_k})\rangle\\
 &=&\seq{(\bar{p}_1,\bar{\cal B}^0,\bar{\Pi}^0),
\ldots (\bar{p}_{k+1},\bar{\cal B}^k,\bar{\Pi}^k,
\bar{\Pi}^{k-1,k},\bar{\cal C}^k,\bar{\mbox{H}}^k,
\bar{\gF}^k,\bar{\Psi}^k)}
\end{array}
$$
where we set $m_0=0$ so that $\Pi^{m_0,m_1}=\Pi^{0,m_1}$.
\end{definition}

\noindent
{\bf Remarks.}
\begin{enumerate}
\item
In rearranging \br\ into {\bf t}, a few of
the $B$'s, $C$'s, $\gf$'s, $\gp$'s 
and $\psi$'s may be lost if $m_k<n$, 
but all the information coded by the $p$'s 
is preserved, as $m_{k+1}=n+1$.

\item
Suppose that for any $m\leq n$ and any
$\gth_m<\ga\leq\gth_{m+1}$ there is 
$\gth_{m-1}<\gb\leq\gth_m$
such that $\gb<_T\ga$. That amounts to say 
that in extending $p_m$ to $p_{m+1}$, 
we never visit any model in the iteration 
trees $p_l$, for $l<m$.
Then for $\gth_n<\ga\leq\gth_{n+1}$ the
backward sequence of $\ga$ relative to
$p_1,\ldots,p_{n+1}$ is of length $n+1$ and the 
contraction of $\br$ relative to $\ga$ is $\br$ itself.
Unfortunately, we cannot assume this holds in general, 
and this is why we had to introduce this further
complication. And, after all, life was not meant to be easy.

\end{enumerate}
Before we move on, we still must verify that

\begin{lemma}
$\mbox{\bf t}\in{\cal U}^-$
\end{lemma}

\begin{proof}
By induction on $k+1=$ the length of {\bf t} = the length of the 
backward sequence of $\ga$. 

Assume $k=0$. Then $\mbox{\bf t}=\seq{(p_{n+1},{\cal B}^0,\Pi^0)}
=\seq{(\bar{p}_1,\bar{\cal B}^0,\bar{\Pi}^0)}$.
By remark (4), $\bar{p}_1=p_{n+1}$ can be copied on 
$(\bar{\cal B}^0,\bar{\Pi}^0)=({\cal B}^0,\Pi^0)$,
so clause (3) holds. The other clauses are immediate or vacuous.

So we may assume $k>0$. The embedding 
$\bar{\gf}_k^{\ast}:M_{\bar{\gth}_k}^{\bar{p}_k}=
M_{\gth_{m_k}}^{p_{m_k}}\to V_{\eta}$ 
can be taken to be $\gf_{m_k}$, so clause (2.c) holds. 
As 
$$
m_{k}+1\leq m_{k+1}\implies 
I_{p_{m_k +1}}\supseteq I_{p_{m_{k+1}}}=I_{\bar{p}_{k+1}}
$$
and $N(\bar{p}_k)=N(p_{m_k})\notin I_{p_{m_k +1}}$, then 
$N(\bar{p}_k)\notin I_{\bar{p}_{k+1}}$. 
For any $q\supset p_{m_k}$, $q$ can be 
be copied on $({}^{\ast}{\cal C}^{m_k},{}^{\ast}\gF^{m_k})=
({}^{\ast}\bar{\cal C}^k,{}^{\ast}\bar{\gF}^k)$
if and only if $q$ can be copied on 
$({\cal B}^{m_k},\gP^{m_k})=(\bar{\cal B}^k,\bar{\gP}^k)$,
so by clause (4.f) and remark (4), 
and taking $q=\bar{p}_{k+1}=p_{m_{k+1}}$, clauses (3) and (4.f) hold.
The other clauses are left to the reader.
\end{proof}

\noindent
A branch  $\mbox{\bf b}$ through ${\cal U}^-$
$$
\mbox{\bf b}= \seq{(p_1,{\cal B}^0,\Pi^0),
\ldots,(p_{n+1},{\cal B}^n,\Pi^n,\Pi^{n-1,n},
{\cal C}^n,\mbox{H}^n,\gF^n,\Psi^n),\ldots}
$$
will be constructed inductively together with sequences 
$$
\seq{(P_{\ga},\gk_{\ga},
\br_{\ga})\mid\ga<\gth}
$$
$$
\seq{(P^{\ast}_n,\eta^{\ast}_n, \gk^{\ast}_n,
\gf^{\ast}_n,\psi^{\ast}_n,\br^{\ast}_n)\mid n<\go}
$$
where $\gth=\sup \gth_n=ht(\cal T)$, ${\cal T}=\bigcup_n p_n$.  
Let $\gr_{\ga}=\gr_{\ga}^{\cal T}$.

\bigskip\noindent
$\bullet$ {\bf The conditions.}
The following 9 conditions must hold for
every $0\leq n<\go$.

\begin{enumerate}
\item\begin{enumerate}
     \item
     $({\cal P}^n,{}^{\ast}\gF^n)=
     \seq{(P_{\ga},\gf_{\ga})\mid\ga<\gth_n}^{\frown}
     \seq{(P^{\ast}_n,\gf^{\ast}_n)}$
     is an enlargement for $p_n$ with bounds
     $\seq{\eta_{\ga}\mid\ga<\gth_n}^{\frown}
     \seq{\eta^{\ast}_n}$.

     \item
     $({\cal P}^n,{}^{\ast}\Psi^n)=
     \seq{(P_{\ga},\psi^n_{\ga})|\ga<\gth_n}^{\frown}
     \seq{(P^{\ast}_n,\psi^{\ast}_n)}$
     is an enlargement for $(\Pi^n p_n,{\cal B}^n)$ 
     with bounds $\seq{\gk_{\ga}\mid\ga<\gth_n}^{\frown}
     \seq{\gk^{\ast}_n}$. 
     
      \item
      $\gf_{\ga}=\psi_{\ga}^n\circ\pi^n_{\ga}$, 
      for $\ga<\gth_n$ and
      $\gf^{\ast}_n=\psi^{\ast}_n \circ\pi^n_{\gth_n}$.

      \item
      For $0\leq m<n$ and $\ga<\gth_m$, 
      $\psi^n_{\ga}=\psi^m_{\ga}\circ\pi^{m,n}_{\ga}$.
     \end{enumerate}

\item 
$\seq{(P_{\ga},\gf_{\ga})\mid\ga<\gth_n}\in P^{\ast}_n$
and for every $\ga<\gth_n$
$$
P^{\ast}_n\models P_{\ga}\mbox{ is $2^{\aleph_0}$-closed 
and of size } |V_{\gf^{\ast}_n(\gr_{\ga})+1}|
$$
and
$$
\gf^{\ast}_n\re V_{\gr_{\gth_n}+1}=
\gf_{\gth_n}\re V_{\gr_{\gth_n}+1}.
$$

\item
$p_{n+1}$ can be copied on $({\cal P}^n,{}^{\ast}\gF^n)$.

\item
Let $M(p_n)\in I_{p_n}$ be the least $k$ such that 
$\exists q\supset p_n\exists\gb( k=\mbox{\bf c}
(q\re\gb)$, $q$ is according to $\gS$ and can be copied on 
$({\cal P}^n,{}^{\ast}\gF^n))$.
Then $M(p_n)\notin I_{p_{n+1}}$.

\item\begin{enumerate}
     \item
     $\gd(P_{\ga})=\eta_{\ga}$ and
     $P_{\ga}\models\lq\lq |V_{\eta_{\ga}}|<\gk_{\ga}$ and
     $\gk_{\ga}$ is a cut-off point."

     \item
     $P_{\ga}\cap V_{\gk_{\ga}}
      =C_{\ga}\supset H_{(2^{\aleph_0})^+}$.
     \end{enumerate}

\item
For $\gth_n\leq\ga<\gth_{n+1}$, let ${\cal U}_{\ga}$ be the tree $\cal U$
relativized to $P_{\ga}$, where the ordinals $\eta$, $\gk$
are interpreted as $\eta_{\ga}$, $\gk_{\ga}$
and $\pi$ is replaced by $\gf_{\ga}\circ i_{0,\ga}^{p_{n+1}}$.
${\cal U}_{\ga}$ is wellfounded
and $||\ ||$ denotes its rank function.

Then $\br_{\ga}$ is a non-empty node of
${\cal U}_{\ga}$ and

     \begin{enumerate}

     \item
     the $\gn$'s of $\br_{\ga}$ are 
     $\seq{\gn_0,\ldots,\gn_k}$, where 
     $\seq{(\gn_0,m_1),\ldots,(\gn_k,m_{k+1})}$ is 
     the backward sequence of $\ga$ relative to 

     {\em (i)\/} $p_1,\ldots, p_{n+1}$, if $\ga\neq\gth_n$
                 or $n=0$;

     {\em (ii)\/} $p_1,\ldots, p_n$, if $\ga=\gth_n$ and $n>0$.

     \item
     $(\br_{\ga})^-$ is the contraction of 

     {\em (i)\/} $\mbox{\bf b}\re n+1$ relative to $\ga$,
                 if $\ga\neq\gth_n$ or $n=0$; 

     {\em (ii)\/} $\mbox{\bf b}\re n$ relative to $\ga$,
                 if $\ga=\gth_n$ and $n>0$. 

     In particular, ${\br_0}^-=\seq{p_1,{\cal B}^0,\Pi^0}$. 

     \item
     $\gs_k=\psi^{n+1}_{\ga}\circ\pi^{m,n+1}_{\ga}
     \re\chunk{M_{\ga}^{(\Pi^m p_{n+1},{\cal B}^m)}
     ;\gn_0,\ldots,\gn_k}$, where $m=m_k$ and $n+1=m_{k+1}$.
     Note that $\mbox{ran}(\gs_k)\subseteq C_{\ga}
     =P_{\ga}\cap V_{\gk_{\ga}}$, the relativization
     of $V_{\gk}$ to $P_{\ga}$.

     \item
     $P_{\ga}$ has at least $||\br_{\ga}||\cdot 2$ cut-off 
     points above $\gk_{\ga}$;
     \end{enumerate}

\item
\begin{enumerate}
     \item
     $\gd(P_n^{\ast})=\eta^{\ast}_n$  and
     $P_n^{\ast}\models\lq\lq |V_{\eta_n^{\ast}}|<\gk^{\ast}_n$
     and $\gk_n^{\ast}$ is a cut-off point." 

     \item
     $\gf^{\ast}_n:M_{\gth_n}^{p_n}\to 
     P_n^{\ast}\cap V_{\eta^{\ast}_n}$ and
     $\psi^{\ast}_n:M_{\gth_n}^{\Pi^n p_n}
     \to P_n^{\ast}\cap V_{\gk^{\ast}_n}$
     are elementary embeddings  such that
     $\gf^{\ast}_n=\pi^n_{\gth_n}\circ\psi^{\ast}_n$ and
     $$
     P^{\ast}_n\models\lq\lq
     ({\cal C}^{n\frown}\seq{V_{\gk^{\ast}_n}},
     \gF^{n\frown}\seq{\gf^{\ast}_n})
     \mbox{ is an enlargement of }(p_n, M)"
     $$
     and
     $$
     P^{\ast}_n\models\lq\lq
     ({\cal C}^{n\frown}\seq{V_{\gk^{\ast}_n}},
     \Psi^{n\frown}\seq{\psi^{\ast}_n}) 
     \mbox{ is an enlargement of }
     (\Pi^n p_n,{\cal B}^n)." 
     $$

     \item
     $P^{\ast}_n\cap V_{\gk_n^{\ast}}\supset H_{(2^{\aleph_0})^+}$.
\end{enumerate}

\item
Let ${\cal U}^{\ast}_n$ be the relativization of
$\cal U$ to $P^{\ast}_n$, with the ordinals $\eta$, $\gk$
interpreted as $\eta^{\ast}_n$, $\gk^{\ast}_n$
and $\pi$ interpreted as 
$\gf_n^{\ast}\circ i_{0,\gth_n}^{p_n}$.
${\cal U}^{\ast}_n$ is wellfounded
and $||\ ||$ denotes its rank function.

Then $\br_n^{\ast}$ is a 
node of ${\cal U}^{\ast}_n$ and 
\begin{enumerate}

     \item
     The $\gn$'s of $\br_n^{\ast}$ are 
     $\seq{\gn_0,\ldots,\gn_k}$, where 
     $\seq{(\gn_0,m_1),\ldots,(\gn_k,m_{k+1})}$ is 
     the backward sequence of $\gth_n$ relative to
     $p_1,\ldots, p_n$; 

      \item
      $(\br_n^{\ast})^-$ is the contraction of 
      $\mbox{\bf b}\re n$ relative to $\gth_n$.
      In particular $\br^{\ast}_0=\emptyset$.

     \item
     For $n>0$, $\gs_k=\psi_n^{\ast}\circ\pi^{m,n}_{\gth_n}
     \re\chunk{M_{\gth_n}^{(\Pi^m p_n,{\cal B}^m)};
     \gn_0,\ldots,\gn_k}$, where $m=m_k$ and $n=m_{k+1}$. 
     Note that $\mbox{ran}(\gs_k)\subseteq 
     P^{\ast}_n\cap V_{\gk_n^{\ast}}$, the interpretation
     of $V_{\gk}$ in $P_n^{\ast}$.

     \item
     $P_n^{\ast}$ has at least $||\br_n^{\ast}||\cdot 2+1$ 
     cut-off points above $\gk^{\ast}_n$.
     \end{enumerate}

\item
$P_{n+1}^{\ast}\in P^{\ast}_n$.
\end{enumerate}

\noindent
This contradiction will show that 
our assumption about $\gS$ being
a winning strategy for $I$ in ${\cal G}^+(M)$ is 
false, hence the theorem will be proved.

\bigskip\noindent
$\bullet$ {\bf Base step.}

\noindent
Let $\gz>\gk$ be large enough so that 
$V_{\gz}$ is a premouse with $\gd(V_{\gz})=\eta$
and such that there are $||{\cal U}||\cdot 2+1$ 
cut-off points above $\gk=\gk^{\ast}_0>\eta^{\ast}_0=\eta$. 
Let 
$$
P^{\ast}_0=V_{\gz},\qquad\br^{\ast}_0=\emptyset
\qquad\mbox{and}\qquad\gf^{\ast}_0=\pi.
$$ 
Now we have to define ${\cal B}^0=\seq{B_0}$.

\smallskip\noindent
{\bf Case 1:} $M(p_0)=\min I_{p_0}$.

\noindent 
Recall that $p_0=\emptyset$ and $I_{p_0}=\go$, 
so there is a position
$q$ according to $\gS$ such that $0=\mbox{\bf c}(q\re\gb)$
and $q$ can be copied on 
$({\cal P}^0,{}^{\ast}\gF^0)=(V_{\gk},\pi)$.
In this case, let $\psi^{\ast}_n:B_0\to V_{\gk}$,
such that $|B_0|=2^{\aleph_0}$,
and let $\psi^{-1}(\pi)=\pi_0$.

\smallskip\noindent
{\bf Case 2:}  $M(p_0)>\min I_{p_0}$.

\noindent
For every integer $k<M(p_0)$ and for every
position $q$ according to $\gS$ such that 
$\mbox{\bf c}(q\re\gb)=k$ for some $\gb\leq lh(q)$,
there are witnesses to the fact that such $q$ cannot be copied
on $({\cal P}^0,{}^{\ast}\gF^0)$. 
That is, for any such $q$, there is 
an ordinal $\gb$ such that $q\re\gb$ can be 
copied on $({\cal P}^0,{}^{\ast}\gF^0)$,
but $q\re\gb+1$ cannot. 
Let $\cal S$ be $q\re\gb$ copied 
on $({\cal P}^0,{}^{\ast}\gF^0)$.

If $\gb$ is a limit ordinal, then fix an increasing sequence
$\gb_m\to \gb$ and ordinals $\xi_m\in M^{\cal S}_{\gb_m}$
such that $i_{\gb_m,\gb_{m+1}}(\xi_m)>\xi_{m+1}$, 
and let $\mbox{w}_m=(\xi_m,\gb_m)$.
If $\gb=\nu+1$ and $\gg=<_T\mbox{-pred}(\gb)$,
then let $\mbox{w}_m=(a_m,f_m)$ witness the
illfoundedness  of $M_{\gb+1}^{\cal S}$.

By absoluteness $\seq{\mbox{w}_m\mid m<\go}$ 
can be taken to be inside $P^{\ast}_0 =V_{\gz}$, and
by Corollary \ref{truncation}, we can assume
that $\seq{\mbox{w}_m\mid m<\go}\in V_{\gk}$. 
Repeating the argument for every position $q$ as above,
a set $X\subseteq V_{\gk}$ of all such
$\mbox{w}_m(q)$'s is obtained.
$|X|\leq 2^{\aleph_0}$, as there are at most $2^{\aleph_0}$
such $q$'s. 
Working inside $V_{\gz}$, let 
$$
H=\hull^{V_{\gk}}\big(X\cup HC\cup\{\eta,\pi\}\big).
$$
By construction $|H|=2^{\aleph_0}$. Let 
$\psi^{\ast}_0:B_0\to H$ be the inverse of the
transitive collapse and let $\pi$, $\eta$ be 
the images of $\pi_0$ and $\ge_0$ via $\psi^{\ast}_0$. 

\medskip
Thus in both cases ${\cal B}^0$, $\Pi^0$ and 
$\psi^{\ast}_0$ are defined.
Finally, choose a position $p_1$ according
to $\gS$ that can be copied on $V_{\gk}$
via $\pi$, and such that $M(p_0)\notin I_{p_1}$. 

\bigskip\noindent
$\bullet$ {\bf Inductive step.}

\noindent
Let $n\geq 0$ and suppose we are given
$$
\begin{array}{l}
\mbox{\bf b}\re n+1 \in{\cal U}^-\\
\seq{(P_{\ga},\eta_{\ga},\gk_{\ga},\gf_{\ga},
\br_{\ga})\mid\ga<\gth_n}\\
\seq{\psi^n_{\ga}\mid\ga<\gth_n}\qquad\mbox{and}\\
\seq{(P_m^{\ast},\eta^{\ast}_m,\gk^{\ast}_m,\gf^{\ast}_m,
\psi^{\ast}_m,\br^{\ast}_m)\mid m\leq n}
\end{array}
$$
satisfying conditions (1)---(9) above.

\bigskip
\noindent
$\bullet$ {\bf Construction of 
${\cal P}^{n+1}$ and ${}^{\ast}\gF^{n+1}$.}

\noindent
We will now build 
$$
\begin{array}{lc}
\seq{(P_{\ga},\eta_{\ga},\gk_{\ga},\gf_{\ga},
\br_{\ga})\mid\gth_n\leq\ga<\gth_{n+1}} &
\qquad\mbox{and}\\
{}&{}\\
(P_{n+1}^{\ast},\eta^{\ast}_{n+1},\gk^{\ast}_{n+1},
\gf^{\ast}_{n+1},\br^{\ast}_{n+1}).&{}
\end{array}
$$

Let $\cal W$ be the pseudo-iteration tree obtained by
copying $p_{n+1}$ on $({\cal P}^n,{}^{\ast}\gF^n)$ and
denote its $\ga$th model by $W_{\ga}$.
Let's also agree that for $\gb\leq\gth_{n+1}$
$$
G_{\gb}:M_{\gb}^{(\Pi^n p_{n+1},{\cal B}^n)}\to W_{\gb},
\qquad\qquad
f_{\gb}:M_{\gb}^{(p_{n+1},M)}\to M_{\gb}^{(\Pi^n p_{n+1},{\cal B}^n)}
$$
are the copy map induced by
${}^{\ast}\Psi^n$ and $\Pi^n$, respectively.

\begin{definition}
\label{a}
Let $\gth_n\leq\gb\leq\gth_{n+1}$ and let 
$\ga\leq\gth_n$ be its root in $\cal W$.
Let 
$$
(\eta)^{\gb}, (\gk)^{\gb}, ({\cal U})^{\gb}, \hbox{\bf t}_{\ga} 
=\left\{ \begin{array}{ll}
i_{\ga,\gb}(\eta_{\ga}), i_{\ga,\gb}(\gk_{\ga}), 
i_{\ga,\gb}({\cal U}_{\ga}),
\br_{\ga} &\quad\mbox{if $\ga<\gth_n$,}\\
{}&{}\\
i_{\ga,\gb}(\eta^{\ast}_n), i_{\ga,\gb}(\gk^{\ast}_n), 
i_{\ga,\gb}({\cal U}^{\ast}_n), 
\br_n^{\ast} &\quad\mbox{if $\ga=\gth_n$.}
\end{array}\right.
$$

The node $\hbox{\bf s}_{\gb}\in ({\cal U})^{\gb}$
is defined as follows.

\begin{enumerate}

\item
If $n=0$ and $\gb\geq\ga=\gth_0=0$, then
$$
\mbox{\bf s}_{\gb}=\seq{(p_1,{\cal B}^0,\Pi^0,\nu_0,\gs_0)}
$$
where $\nu_0=\gb$ and $\gs_0=G_{\nu_0}\re
\chunk{M_{\nu_0}^{(\Pi^0 p_1, {\cal B}^0)};\nu_0}$
is the restriction of the copy map induced 
by $\psi^{\ast}_0:B_0\to V_{\gk}$.

\item
If $n>0$ and $\ga=\gb=\gth_n$, then 
$\hbox{\bf s}_{\gb}=\br^{\ast}_n$.

\item
If $n>0$ and $\gb>\gth_n$, then let
$\seq{(\gn_0,m_1),\ldots,(\gn_k,m_{k+1}),(\gn_{k+1},m_{k+2})}$ 
be the backward sequence of $\gb$ relative 
to $p_1,\ldots ,p_{n+1}$.  Hence
$\nu_k=\ga<\gth_{m_{k+1}}$ and 
$\seq{(\gn_0,m_1),\ldots,(\gn_k,m_{k+1})}$ 
is the backward sequence of $\ga$ relative
to $p_1,\ldots, p_{m_{k+1}}$.
Then the $\nu$'s of $\mbox{\bf s}_{\gb}$
are $\nu_0,\ldots,\nu_{k+1}$
and let  
$$
\gs_{k+1}=G_{\gb}\re
\chunk{M_{\gb}^{(\Pi^{m_{k+1}} p_{n+1},{\cal B}^{m_{k+1}})};
\gn_0,\ldots,\gn_{k+1}}
$$
Let also
$$
\hbox{\bf s}_{\gb}=
i_{\ga,\gb}^{\cal W}({\hbox{\bf t}_{\ga}})^{\frown}\langle 
(p_{n+1},{\cal B}^{m_{k+1}}, \Pi^{m_{k+1}},\Pi^{m_k,m_{k+1}}, 
{\cal C}^{m_{k+1}},\mbox{H}^{m_{k+1}},\gF^{m_{k+1}},
\Psi^{m_{k+1}},\gn_{k+1},\gs_{k+1})\rangle.
$$
\end{enumerate}
\end{definition}

\begin{lemma}
\label{A}
With $\ga$ and $\gb$ as in the definition above, 
and $m=m_{k+1}$, $l=m_k$, if $k>0$, or
$l=0$ otherwise
\begin{enumerate}

\item
$\mbox{\bf t}_{\ga}^-\in W_{\ga}\cap 
V_{\gr^{\cal W}_{\gth_l}+1}$ and
$i^{\cal W}_{\ga,\gb}\big({\mbox{\bf t}_{\ga}}^-\big)=
{\mbox{\bf t}_{\ga}}^-$; 

\item
${\mbox{\bf s}_{\gb}}^-,\mbox{\bf b}\re n+1\in W_{\gb}\cap 
V_{\gr_{\gth_m}^{\cal W}+1}$ 
and ${\mbox{\bf s}_{\gb}}^-$ is the contraction 
of $\mbox{\bf b}\re n+1$ relative to $\gb$.
\end{enumerate}
\end{lemma}

\begin{proof}
If $n=0$ then $\mbox{\bf t}_{\ga}=\emptyset$
and the Lemma is immediate. So we may assume $n>0$,
hence $k>0$.

(1) By the definition of backward sequence
$\gth_l\leq\gth_{m-1}<\ga\leq\gth_m$ and the 
critical point of $i=i^{\cal W}_{\ga,\gb}$
is $>\gr=\gr_{\gth_l}^{\cal W}$, so it is 
enough to show that ${\mbox{\bf t}_{\ga}}^-\in W_{\ga}
\cap V_{\gr +1}$.
For any $j$, $p_j, {\cal B}^j, \Pi^j,\nu_j \in H_{(2^{\go})^+}$, 
hence their rank is certainly less than $\gr$,
so the only possible source of problems are the
$\mbox{H}^j,\gF^j,\Psi^j$ and ${\cal C}^j$, for $j\leq l$.
Clause (2) in the definition of $\cal U$, 
implies that the rank of ${\cal C}^l=\seq{C_{\ga}\mid\ga<\gth_l}$
is $\leq\gr$, so 
${\mbox{\bf t}_{\ga}}^-\in W_{\gth_l}\cap V_{\gr +1}$.
But $W_{\gth_l}$ agrees with $W_{\ga}$ through $\gr +2$,
so the result follows at once.

(2) By part (1) ${\mbox{\bf s}_{\gb}}^-$ 
extends ${\mbox{\bf t}_{\ga}}^-$, so
$({\mbox{\bf s}_{\gb}}^-)\re n\in W_{\ga}\cap V_{\gr^{\cal W}_{\ga}}$.
By an argument as the one in part (1) we can conclude that 
$(p_{n+1},{\cal B}^m,\Pi^m,\Pi^{l,m},{\cal C}^m,\mbox{H}^m,
\gF^m,\Psi^m)\in W_{\gb}\cap V_{\gr_{\gth_m}^{\cal W}+1}$.
By the agreement between $W_{\ga}$ and $W_{\gb}$,
${\mbox{\bf s}_{\gb}}^-\in W_{\gb}\cap V_{\gr^{\cal W}_{\gth_m +1}}$.
Note that in the course of the proof we 
also managed to prove that ${\mbox{\bf s}_{\gb}}^-$ is the
contraction of $\mbox{\bf b}\re n+1$, and that
$\mbox{\bf b}\re n+1\in W_{\gb}\cap V_{\gr^{\cal W}_{\gth_m +1}}$.
\end{proof}

\begin{lemma}
\label{B}
$(\forall\gb\geq\gth_n)
\mbox{\bf s}_{\gb}\in ({\cal U})^{\gb}.$
\end{lemma}

\begin{proof}
We will use the same notation as in the proof
of the previous lemma.
If $n=0$ or if $\ga=\gb=\gth_n$ and $n>0$, 
then the result follows easily,
so we may assume $\gb>\gth_n$.

Let's show first that 
${\mbox{\bf s}_{\gb}}^-\in({\cal U}^-)^{\gb}$.
By part (1) of Lemma \ref{A} and ${\mbox{\bf t}_{\ga}}^-\in 
({\cal U}^-)^{\ga}$, it follows that 
${\mbox{\bf t}_{\ga}}^-\in({\cal U}^-)^{\gb}$,
so we only have to check clause (3)
in the definition of $\cal U$, relativized to $W_{\gb}$.
By Corollary \ref{truncation} and Remark (5),
$M(p_m)=$ the least $k\in I_{p_m}$ such that 
there is $q\supset p_m$ according to $\gS$,
and $k\notin I_q$ and $q$ can be copied on
$({\cal B}^m,\Pi^m)$. Moreover $p_{n+1}$
is such a $q$. 
As $p_{n+1},{\cal B}^m,\Pi^m\in H_{(2^{\aleph_0})^+}
\in W_{\gb}$, by Remark (5) relativized to 
$W_{\gb}$, we have that
$$
W_{\gb}\models M(p_m)=N(p_m)\mbox{ and }
{\mbox{\bf s}_{\gb}}^-\in({\cal U}^-)^{\gb}.
$$

As $\mbox{\bf s}_{\gb}$ extends 
$i_{\ga,\gb}(\mbox{\bf t}_{\ga})\in(\cal U)^{\gb}$ 
we only have to take care of $\gb=\gn_{k+1}$ and $\gs_{k+1}$.
Recall that $\seq{(\nu_0, m_1),\ldots,(\nu_k,m_{k+1})}$ is
the backward sequence of $\ga=\nu_k$, $m=m_{k+1}$,
$l=m_k$, and that
$\gth_l\leq\gth_{m-1}<\ga\leq\gth_m$.
Let us verify clause (6) in the definition of $\cal U$. 
Let $\gs_k=\gs_k(\mbox{\bf s}_{\gb})$ and
$\gt=\gs_k(\mbox{\bf t}_{\ga})$. By definition
of $\mbox{\bf s}_{\gb}$, $\gs_k=i_{\ga,\gb}(\gt)$.
By condition (6.c) or (8.c) 
$$
\gt=\psi\circ\pi^{l,m}_{\ga}\re
\chunk{M_{\ga}^{\Pi^l p_m};\nu_0,\ldots,\nu_k},
$$ 
where $\psi=\psi^{\ast}_n$, if $\ga=\gth_n$, or
$\psi=\psi^m_{\ga}$, otherwise.
For $\gg\leq\gth_{n+1}$, let 
$G_{\gg}:M_{\gg}^{\Pi^n p_{n+1}}\to W_{\gg}$
be the copy map  induced by the embedding
${}^{\ast}\Psi^n:{\cal B}^n\to {\cal P}^n$.
Condition (1.d) implies that 
$\psi^m_{\ga}=\psi^n_{\ga}\circ\pi^{m,n}_{\ga}$, 
when $\ga<\gth_n$, thus the diagram
$$
\setlength{\dgARROWLENGTH}{2.5em}
\begin{diagram}
\node{
\chunk{M_{\ga}^{\Pi^l p_m};
\gn_0,\ldots,\gn_k}
}
\arrow[2]{e,t}{\gt}
\arrow{ese,b}{\pi^{l,n}_{\ga}}
\node[2]{
W_{\ga}\cap V_{(\gk)^{\ga}}
} 
\arrow[2]{e,t}{i}
\node[2]{W_{\gb}\cap V_{(\gk)^{\gb}}} \\

\node[3]{
M_{\ga}^{\Pi^n p_{n+1}}
}
\arrow{n,r}{G_{\ga}=\psi^n_{\ga}}
\arrow[2]{e,b}{j} 
\node[2]{M_{\gb}^{\Pi^n p_{n+1}}}
\arrow{n,r}{G_{\gb}}
\end{diagram}
$$
commutes, where $i$ and $j$ are the embeddings 
of the pseudo-iteration trees $\cal W$ and
$\Pi^n p_{n+1}$, respectively.
As $\gs_k \in W_{\ga}$, and 
$\mbox{dom}(\gt)$ is hereditarily countable, 
then $i\circ\gt= i(\gt)$.
Thus letting 
$\gs_{k+1}=\gs_{k+1}(\mbox{\bf s}_{\gb})$ be $G_{\gb}$ 
restricted to the appropriate support, 
we have the commutative diagram
$$
\setlength{\dgARROWLENGTH}{3.5em}
\begin{diagram}
\node{
\chunk{M_{\ga}^{\Pi^l p_m};\gn_0, \ldots,\gn_k} }
\arrow[2]{e,t}{\gs_k}
\arrow{se,b}{(i^{\Pi^m p_{n+1}}_{\ga,\gb})
\circ\,\pi^{l,m}_{\ga}}
\node[2]{W_{\gb}\cap V_{(\gk)^{\gb}} } \\
\node[2]{\chunk{M_{\gb}^{\Pi^m p_{n+1}};
\gn_0,\ldots, \gn_{k+1}}}
\arrow{ne,r}{\gs_{k+1}}
\end{diagram}
$$
which is what we had to prove.
\end{proof}

Thus $\mbox{\bf s}_{\gb}$ is defined and 
belongs to $\big({\cal U}\big)^{\gb}$, for all 
$\gth_n\leq\gb\leq\gth_{n+1}$. 

\begin{lemma}
For all $n\geq 0$,
\begin{enumerate}

\item $P^{\ast}_n=W_{\gth_n}$ has at least 
$||\mbox{\bf s}_{\gth_n}||\cdot 2+1$ cut-off points
above $\gk^{\ast}_n$.

\item
For $\gth_n<\gb\leq\gth_{n+1}$,
$W_{\gb}$ has at least $||\mbox{\bf s}_{\gb}||\cdot 2+2$
cut-off points above $(\gk)^{\gb}$.

\end{enumerate}
\end{lemma}

\begin{proof}
Let $\ga$, $\gb$ and $\mbox{\bf t}_{\ga}$ be
as in Definition \ref{a}.

When $n=0$ and $0\leq\gb\leq\gth_n$, then
$\ga=0$ and $\mbox{\bf t}_{\ga}=\mbox{\bf r}_0^{\ast}=\emptyset$.
$V_{\gz}=P_0^{\ast}=W_{\gth_0}$ has at least
$||\mbox{\bf r}^{\ast}_0||\cdot 2+1$ cut-off points
above $\gk=\gk^{\ast}_0=(\gk)^0$, hence
$W_{\gb}$ has at least 
$i_{0,\gb}(||\mbox{\bf r}^{\ast}_0||\cdot 2+1)$
cut-off points above $(\gk)^{\gb}$.
As $\mbox{\bf s}_{\gb}$ properly extends
$i_{\ga,\gb}(\mbox{\bf r}^{\ast}_0)=\mbox{\bf r}^{\ast}_0$
for any $0\leq\gb\leq\gth_1$, (1) and (2) follow at once.

Suppose now $n>0$. Part (1) follows from
condition (8.d) and the fact that 
$\mbox{\bf s}_{\gth_n}=\mbox{\bf r}^{\ast}_n$, so
we may assume $\gb>\gth_n$. 
By condition (6.d), $W_{\ga}$ has at least 
$||\mbox{ \bf t}_{\ga}||\cdot 2$ cut-off points
above $(\gk)^{\ga}$, so $W_{\gb}$ has at least 
$i_{\ga,\gb}(||\mbox{\bf t}_{\ga}||\cdot 2)=
||i_{\ga,\gb}(\mbox{\bf t}_{\ga})||\cdot 2$
cut-off points above 
$(\gk)^{\gb}=i_{\ga,\gb}\big((\gk)^{\ga}\big)$.
As $\mbox{\bf s}_{\gb}$ is a proper extension
of $i_{\ga,\gb}(\mbox{\bf t}_{\ga})$,
$||\mbox{\bf s}_{\gb}||\cdot 2+2\leq
||i_{\ga,\gb}(\mbox{\bf t}_{\ga})||\cdot 2$.
\end{proof}

Let's introduce one more piece of notation.
For any $\gth_n<\gb\leq\gth_{n+1}$, let 
$q(\mbox{\bf s}_{\gb})=\gs_k(\mbox{\bf s}_{\gb})$,
where $k=k(\gb)=lh(\mbox{\bf s}_{\gb})-1$. 
As
$$
q(\mbox{\bf s}_{\gb})=
G_{\gb}\re\chunk{M_{\gb}^{(\Pi^n p_{n+1},{\cal B}^n)};
\nu_0,\ldots,\nu_k}
$$
is the restriction to a chunk 
of the $\gb$th map of the embedding 
${}^{\ast}\Psi:(\Pi^n p_{n+1},{\cal B}^n)\to
({\cal W},{\cal P}^n)$,
$q(\mbox{\bf s}_{\ga})$ and $q(\mbox{\bf s}_{\gb})$
are compatible below $\gr_{\ga}^{\cal W}+2$, 
for $\gth_n<\ga\leq\gb\leq\gth_{n+1}$.
(Recall that two functions $f$ and $g$ are 
compatible below an ordinal $\gr$ iff
$f\re V_{\gr}\cup g\re V_{\gr}$ is still a function.)

Here is the plan of what comes next.
First we construct $P_{\gth_n}$, $\gr_{\gth_n}$ 
and $\gf_{\gth_n}$.
Then fix $\gb$, $\gth_n<\gb<\gth_{n+1}$, 
and work inside $W_{\gb}$. 
Let $\xi$ be the $||\hbox{\bf s}_{\gb}||\cdot 2+1$st 
cut-off point above $(\gk)^{\gb}$ and let 
$Q_{\gb}\supset V_{\gr_{\gb}^{\cal W}}$ 
be the transitive collapse of a hull of $V_{\xi}$.
By exercising proper care $Q_{\gb}$ can 
be taken to be $2^{\aleph_0}$-closed,
of size $|V_{\gr_{\gb}^{\cal W}+1}|$ 
and such that there are 
$||\bar{\mbox{\bf s}}_{\gb}||\cdot 2$
cut-off points above $(\bar{\gk})^{\gb}$,
where $\bar{\mbox{\bf s}}_{\gb}$ and $(\bar{\gk})^{\gb}$
are the collapses of $\mbox{\bf s}_{\gb}$ and $(\gk)^{\gb}$.
We also want embeddings $q_{\gb}:M_{\gb}^{\Pi^n p_{n+1}}
\to Q_{\gb}\cap V_{(\bar{\gk})^{\gb}}$ such that
$q_{\gb}\supseteq q(\bar{\mbox{\bf s}}_{\gb})$ and
agree through $\pi_{\gb}(\gr_{\gb})+2$,
for $\gb\leq\gg<\gth_{n+1}$, $q_{\gb}$ and $q_{\gg}$.
As each $(Q_{\gb}, q_{\gb},\bar{\mbox{\bf s}}_{\gb})$ 
can be coded as an element
of $V_{\gr_{\gb}^{\cal W}+2}\cap W_{\gb}$, 
it belongs to $W_{\gth_{n+1}}$. Finally, working inside
$W_{\gth_{n+1}}$, choose a $2^{\aleph_0}$-closed 
Skolem hull $H$ of $V_{\xi}\cap W_{\gth_{n+1}}$,
where $\xi$ is the $||\mbox{\bf s}_{\gth_{n+1}}||\cdot 2+2$nd
cut-off point, so that $H$ contains all of the $Q_{\gb}$,
$q_{\gb}$, $\bar{\mbox{\bf s}}_{\gb}$, etc.
By letting $h:H\to P^{\ast}_{n+1}$ be the transitive
collapse and $h(Q_{\gb})=P_{\gb}$,
$h(\bar{\mbox{\bf s}}_{\gb})=\mbox{\bf r}_{\gb}$
and $\gf_{\gb}=h(q_{\gb})\circ f_{\gb}$,
the construction would be completed.
Unfortunately we must be more ingenious than that as
it is not clear that embeddings $q_{\gb}$ as above 
can be found inside $W_{\gb}$: the problem is
that it is difficult to maintain the agreement between 
the $q$'s past the first $\go$ of them.
And even if $q_{\gb}$'s as above were available, 
there is no guarantee that the sequence
$\seq{(Q_{\gb},q_{\gb},\bar{\mbox{\bf s}}_{\gb})
\mid\gth_n<\gb<\gth_{n+1}}$ belongs to $W_{\gth_{n+1}}$.
In order to overcome these problems, a 
sequence of approximations
$$
\seq{(Q^{\gb}_{\ga},q^{\gb}_{\ga},\mbox{\bf s}^{\gb}_{\ga})
\mid\gth_n<\ga\leq\gb}\in W_{\gb}
$$
will be built inductively, for $\gth_n\leq\gb\leq\gth_{n+1}$.

\bigskip

We now construct $P_{\gth_n}$, 
$\gf_{\gth_n}$, and $\br_{\gth_n}$.
Working inside $W_{\gth_n}=P_n^{\ast}$ let 
$\xi$ be the $||\mbox{\bf s}_{\gth_n}||\cdot 2+1$st
cut-off point above $(\gk)^{\gth_n}$, and let
$$
\begin{array}{rlr}
H^0= &\hull^{V_{\xi}}\Big(V_{\gr_{\gth_n}^{\cal W}+1}\cup
\{\hbox{\bf s}_{\gth_n},\psi^{\ast}_n,(\eta)^{\gth_n},(\gk)^{\gth_n},
\mbox{\bf b}\re n\}\Big), & \\

 & & \\

H^{\gg}=&\hull^{V_{\xi}}\Big({}^{2^{\aleph_0}}
\big(\bigcup_{\nu<\gg}H^{\nu}\big)\Big) & 
\mbox{ for }\gg\leq (2^{\aleph_0})^+ \\
\end{array}
$$
and set $H=H^{(2^{\aleph_0})^+}$.
It is easy to see that $H$ is $2^{\aleph_0}$-closed
and of size $|V_{\gr_{\gth_n}^{\cal W}+1}|$. 
(This is why plus-2 trees are used:
had we taken hulls of size $|V_{\gr_{\gth_n}}|$ we could
have run into problems with $|H^{\gg}|$'s, 
if cof$(|V_{\gr_{\gth_n}}|)\leq 2^{\go}$.)
Let $h:H\to P_{\gth_n}$ be the transitive collapse,
let $\br_{\gth_n}=h(\mbox{\bf s}_{\gth_n})$, 
$\gk_{\gth_n}=h(\gk^{\ast}_n)$, 
$\eta_{\gth_n}=h(\eta^{\ast}_n)$.
Also set $g_{\gth_n}=h(\psi^{\ast}_n)$ and
$\gf_{\gth_n}=f_{\gth_n}\circ g_{\gth_n}$.

\begin{definition}
For $\gth_n\leq\gb\leq\gth_{n+1}$, let ${\cal R}^{\gb}\in W_{\gb}$ 
be the set defined as follows.

$\seq{(Q^{\gb}_{\ga},q^{\gb}_{\ga},
\eta^{\gb}_{\ga},\gk^{\gb}_{\ga},\mbox{\bf s}^{\gb}_{\ga})
\mid\gth_n\leq\ga\leq\gb}\in {\cal R}^{\gb}$
if and only if:
\begin{description}

\item[(i)]
$(Q_{\gth_n}^{\gb}, q_{\gth_n}^{\gb}, \eta_{\gth_n}^{\gb},
\gk_{\gth_n}^{\gb},\mbox{\bf s}^{\gb}_{\gth_n})=
(P_{\gth_n}, g_{\gth_n},\eta_{\gth_n},\gk_{\gth_n},\br_{\gth_n})$.

\item[(ii)]
$Q^{\gb}_{\ga}$ is a premouse, $\gd(Q^{\gb}_{\ga})=
\eta^{\gb}_{\ga}$, $\gk^{\beta}_{\ga}$ is a cut-off point and 
$Q^{\gb}_{\ga}\models |V_{\eta^{\gb}_{\ga}}|
<\gk^{\gb}_{\ga}$. 

\item[(iii)]
$({\cal P}^n\re\gth_n,\Psi^n)^{\frown}\seq{(Q^{\gb}_{\ga},
q^{\gb}_{\ga})\mid \gth_n\leq\ga\leq\gb}$ is an enlargement 
for $(\Pi^n p_{n+1}\re\gb+1,{\cal B}^n)$ with bounds
$\seq{\gk_{\ga}\mid\ga<\gth_n}^{\frown}
\seq{\gk^{\gb}_{\ga}\mid\gth_n\leq\ga\leq\gb}$.

\item[(iv)]
$({\cal P}^n\re\gth_n,\gF^n)^{\frown}\seq{(Q^{\gb}_{\ga},
q^{\gb}_{\ga}\circ f_{\ga})\mid\gth_n\leq\ga\leq\gb}$ 
is an enlargement for $(p_{n+1}\re\gb+1,M)$ with bounds
$\seq{\eta_{\ga}\mid\ga<\gth_n}^{\frown}
\seq{\eta^{\gb}_{\ga}\mid\gth_n\leq\ga\leq\gb}$.

\item[(v)]
$\mbox{\bf s}^{\gb}_{\ga}\in({\cal U})^{\gb}_{\ga}$,
where $({\cal U})^{\gb}_{\ga}$ is the 
relativization of ${\cal U}$ to $Q^{\gb}_{\ga}$,
with $\gk$, $\eta$ and $\pi$  replaced by 
$\gk^{\gb}_{\ga}$, $\eta^{\gb}_{\ga}$ and
$q^{\gb}_{\ga}\circ f_{\ga}\circ i_{0,\ga}$; 
moreover for $\ga>\gth_n$,
$(\mbox{\bf s}^{\gb}_{\ga})^-$ is the contraction
of $\mbox{\bf b}\re n+1$ 
relative to $\ga$ and $p_1,\ldots,p_{n+1}$.

\item[(vi)]
$q^{\gb}_{\gb}\re V_{\gr_{\gb}^{\cal W}+2}\subseteq G_{\gb}$,
the copy map, and $q^{\gb}_{\ga}\supseteq
q(\mbox{\bf s}^{\gb}_{\ga})$. 

\item[(vii)]
$W_{\gb}\models\lq\lq Q^{\gb}_{\ga}\supseteq V_{\gr+1}$ 
and is of size $|V_{\gr+1}|"$, where 
$\gr=\gr_{\ga}^{\cal W}=
q^{\gb}_{\ga}(\pi_{\ga}(\gr_{\ga}))$.

\item[(viii)]
$Q^{\gb}_{\ga}$ has at least $||\mbox{\bf s}^{\gb}_{\ga}||\cdot 2$
cut-off points above $\gk^{\gb}_{\ga}$, and 
$Q^{\gth_{n+1}}_{\gth_{n+1}}$ has at least
$||\mbox{\bf s}^{\gth_{n+1}}_{\gth_{n+1}}||\cdot 2+1$ cut-off points
above $\gk^{\gth_{n+1}}_{\gth_{n+1}}$.
\end{description}
\end{definition}

\begin{lemma}
For every $\gb$ with $\gth_n\leq\gb\leq\gth_{n+1}$, 
${\cal R}^{\gb}$ is non-empty. In fact, for any
$\gth_n\leq\gg<\gb$ and any
sequence in ${\cal R}^{\gg}$
there is sequence in ${\cal R}^{\gb}$
extending it.
\end{lemma}

\begin{proof}
By induction on $\gb$.
Condition {\bf (vii)} implies that every element of
${\cal R}^{\gg}$ can be coded as a subset
of $V_{\gr_{\gg}^{\cal W}+1}\cap W_{\gg}$, hence
${\cal R}^{\gg}\subseteq W_{\gb}$, for $\gg\leq\gb$.
If $\gb=\gth_n$ then {\bf (i)} implies ${\cal R}^{\gb}$
has one element only, namely 
$\seq{(P_{\gth_n}, q_{\gth_n},\eta_{\gth_n},
\gk_{\gth_n},\br_{\gth_n})}$. 

Assume the lemma holds for some $\gb\geq\gth_n$ and let's
prove it for $\gb+1$. 
By the inductive hypothesis, it is enough to
show that any
$\seq{(Q^{\gb}_{\ga},q^{\gb}_{\ga},\eta^{\gb}_{\ga},
\gk^{\gb}_{\ga},\mbox{\bf s}^{\gb}_{\ga})
\mid\gth_n\leq\ga\leq\gb}\in {\cal R}^{\gb}$
can be extended to a sequence in ${\cal R}^{\gb+1}$.
By compatibility of $G_{\gb}$ and 
$G_{\gb+1}$  below $\gr_{\gb}^{\cal W}+2$
and by {\bf (vi)}, 
$q^{\gb}_{\gb}\re V_{\gr_{\gb}^{\cal W}+2}
\subseteq G_{\gb+1}$, and as 
$q(\mbox{\bf s}_{\gb+1})\subseteq G_{\gb+1}$,
the maps
$q^{\gb}_{\gb}\re V_{\gr_{\gb}^{\cal W}+2}$ 
and $q(\mbox{\bf s}_{\gb+1})$ are compatible.
First we must find an embedding $q\in W_{\gb}$
such that
$$
q:M_{\gb+1}^{\Pi^n p_{n+1}}\to
W_{\gb+1}\cap V_{(\gk)^{\gb+1}}
\qquad\mbox{and}\qquad
q\supseteq q^{\gb}_{\gb}\re V_{\gr_{\gb}^{\cal W}+2}
\cup q(\mbox{\bf s}_{\gb+1}).
$$
Let $\seq{S_m\mid m<\go}\in W_{\gb+1}$ be
an increasing family of finite supports
for $(\Pi^n p_{n+1},{\cal B}^n)$
such that 
$$
\bigcup_m S_m=\gth_{n+1}+1\quad\mbox{ and }\quad
(M_{\gb+1}^{(\Pi^n p_{n+1},{\cal B}^n)})_{S_0}=
\mbox{dom}(q(\mbox{\bf s}_{\gb+1}))
$$
and let $\cal V$ be the tree, set of finite sequences
closed under initial segment, whose nodes of length $m+1$ are 
embeddings $q_0\subseteq\ldots\subseteq q_m$
$$
q_0=q(\mbox{\bf s}_{\gb+1})\quad\mbox{ and }\quad
q_m:(M_{\gb+1}^{(\Pi^n p_{n+1},{\cal B}^n)})_{S_m}
\to W_{\gb+1}\cap V_{(\gk)^{\gb+1}}
$$
such that $q_m$ is compatible with $q_{\gb}^{\gb}$
below $V_{\gr_{\gb}^{\cal W}+2}$.
$\cal V\in W_{\gb+1}$ as
$q^{\gb}_{\gb}\re V_{\gr^{\cal W}_{\gb}+2}\cup
q(\mbox{\bf s}_{\gb+1})\in W_{\gb+1}$ and it is illfounded
in $V$, hence there is a $q\in W_{\gb+1}$ as desired.
We now proceed as in the construction of $P_{\gth_n}$.
Working inside $W_{\gb+1}$ let 
$\xi$ be the $||\mbox{\bf s}_{\gb+1}||\cdot 2+i$th
cut-off point above $(\gk)^{\gb+1}$, where
$i=1$ if $\gb+1<\gth_{n+1}$, or $i=2$ if $\gb+1=\gth_{n+1}$.
Let
$$
\begin{array}{rlr}
H^0= &\hull^{V_{\xi}}\Big(V_{\gr_{\gb+1}^{\cal W}+1}\cup
\{\hbox{\bf s}_{\gb+1},q,(\eta)^{\gb+1},(\gk)^{\gb+1},
\mbox{\bf b}\re n+1\}\Big), & \\

 & & \\

H^{\gm}=&\hull^{V_{\xi}}\Big({}^{2^{\aleph_0}}
\big(\bigcup_{\nu<\gm}H^{\nu}\big)\Big) & 
\mbox{ for }\gm\leq (2^{\aleph_0})^+ \\
\end{array}
$$
and set $H=H^{(2^{\aleph_0})^+}$.
Let $h:H\to Q_{\gb+1}^{\gb+1}$ be the transitive collapse,
let $\mbox{\bf s}_{\gb+1}^{\gb+1}=h(\mbox{\bf s}_{\gb+1})$, 
$\gk_{\gb+1}^{\gb+1}=h((\gk)^{\gb+1})$, 
$\eta_{\gb+1}^{\gb+1}=h((\eta)^{\gb+1})$.
Also set $q_{\gb+1}^{\gb+1}=h(q)$. 
It is easy to verify that 
$$
\seq{(Q^{\gb}_{\ga},q^{\gb}_{\ga},\eta^{\gb}_{\ga},
\gk^{\gb}_{\ga},\mbox{\bf s}^{\gb}_{\ga})
\mid\gth_n\leq\ga\leq\gb}^{\frown}
\seq{(Q^{\gb+1}_{\gb+1},q^{\gb+1}_{\gb+1},
\eta^{\gb+1}_{\gb+1},\gk^{\gb+1}_{\gb+1},
\mbox{\bf s}^{\gb+1}_{\gb+1})}\in {\cal R}^{\gb+1}
$$
and it extends 
$\seq{(Q^{\gg}_{\ga},q^{\gg}_{\ga},\eta^{\gg}_{\ga}
\gk^{\gg}_{\ga},\mbox{\bf s}^{\gg}_{\ga})
\mid\gth_n\leq\ga\leq\gg}$.
The verification of {\bf (i)}---{\bf (viii)} 
is straightforward.
As an example, let us check that each 
$Q^{\gb+1}_{\ga}$ is $2^{\aleph_0}$-closed
inside $Q^{\gb+1}_{\gb+1}$. 
If $\ga=\gb$, then, by inductive hypothesis 
$Q^{\gb+1}_{\gb}=Q^{\gb}_{\gb}$ is
$2^{\aleph_0}$-closed inside $W_{\gb}$.
But any $2^{\aleph_0}$-sequence of elements
of $Q^{\gb+1}_{\gb}$ can be coded as an element
of $V_{\gr^{\cal W}_{\gb}+2}\cap W_{\gb}=
V_{\gr^{\cal W}_{\gb}+2}\cap W_{\gb+1}$,
hence the result follows.
If $\ga<\gb$, then $Q^{\gb}_{\ga}=Q^{\gb+1}_{\ga}$
is $2^{\aleph_0}$-closed
inside $Q^{\gb}_{\gb}$, which is, as we just showed,
$2^{\aleph_0}$-closed inside $Q^{\gb+1}_{\gb+1}$.
Thus the lemma holds for $\gb+1$.

\bigskip

Suppose now that $\gb>\gth_n$ is limit and that
the result holds for $\gg<\gb$..
Fix a $\gg$ with $\gth_n\leq\gg<\gb$ and a sequence
$$
\bar{\cal Q}=
\seq{(Q^{\gg}_{\ga},q^{\gg}_{\ga},\eta^{\gg}_{\ga},
\gk^{\gg}_{\ga},\mbox{\bf s}^{\gg}_{\ga})
\mid\gth_n\leq\ga\leq\gg}\in{\cal R}^{\gg}.
$$
As ${\cal R}^{\gg}\subseteq W_{\gb}$,
this sequence belongs to $W_{\gb}$.
Working inside $W_{\gb}$,
choose a sequence $\gb_i\to\gb$ with
$\gg<\gb_i<\gb_{i+1}$ and let
$\seq{S_m\mid m<\go}\in W_{\gb+1}$ be
an increasing family of finite supports
for $(\Pi^n p_{n+1},{\cal B}^n)$
such that 
$$
\bigcup_m S_m=\gth_{n+1}+1\quad\mbox{ and }\quad
(M_{\gb}^{(\Pi^n p_{n+1},{\cal B}^n)})_{S_0}=
\mbox{dom}(q(\mbox{\bf s}_{\gb})).
$$
Let $\cal V$ be the tree on $V_{\gr_{\gb}^{\cal W}}$
searching for a sequence like 
$$
\seq{(Q^{\gb}_{\ga},q^{\gb}_{\ga},\eta^{\gb}_{\ga},
\gk^{\gb}_{\ga},\mbox{\bf s}^{\gb}_{\ga})
\mid\gth_n\leq\ga<\gb}
$$
such that 
$$
\bar{\cal Q} \subseteq
\seq{(Q^{\gb}_{\ga},q^{\gb}_{\ga},\eta^{\gb}_{\ga},
\gk^{\gb}_{\ga},\mbox{\bf s}^{\gb}_{\ga})
\mid\gth_n\leq\ga\leq\gb_i}\in {\cal R}^{\gb_i},
$$
together with embeddings 
$q_0\subseteq\ldots\subseteq q_i$ 
$$
q_i:(M_{\gb}^{\Pi^n p_{n+1}})_{S_i}\to
W_{\gb}\cap V_{(\gk)^{\gb}}
$$
such that $q=\bigcup_i q_i\supseteq q(\mbox{\bf s}_{\gb+1})$
and $q$ is compatible with 
$q^{\gb_i}_{\gb_i}$ below $\gr_{\gb_i}^{\cal W}+2$, for $i<\go$.
Using the inductive hypothesis we can choose
$$
{\cal Q}^i=
\seq{(Q^{\gb_i}_{\ga},q^{\gb_i}_{\ga},\eta^{\gb_i}_{\ga},
\gk^{\gb_i}_{\ga},\mbox{\bf s}^{\gb_i}_{\ga})
\mid\gth_n\leq\ga\leq\gb_i} \in{\cal R}^{\gb_i}
$$
such that 
$\bar{\cal Q}\subset{\cal Q}^i\subset{\cal Q}^{i+1}$.
Thus 
$$
\seq{({\cal Q}^i,G_{\gb}\re 
(M_{\gb}^{(\Pi^n p_{n+1},{\cal B}^n)})_{S_i})
\mid i<\go}
$$
is an infinite branch of $\cal V$. By absoluteness
there is ${\cal Q}^{\go}\in W_{\gb}$ such that
$$
\bar{\cal Q} \subset {\cal Q}^{\go}=
\seq{(Q^{\gb}_{\ga},q^{\gb}_{\ga},\eta^{\gb}_{\ga},
\gk^{\gb}_{\ga},\mbox{\bf s}^{\gb}_{\ga})
\mid\gth_n\leq\ga<\gb}
$$
and an embedding
$$
q:M_{\gb}^{\Pi^n p_{n+1}}\to W_{\gb}\cap V_{(\gk)^{\gb}}
$$
compatible with $q^{\gb}_{\ga}$ below 
$\gr_{\ga}^{\cal W}+2$ and such that
$q\supset q(\mbox{\bf s}_{\gb})$.
${\cal Q}^{\go}$ fails to be an element of ${\cal R}^{\gb}$
in that it has no $Q^{\gb}_{\gb}$, $q^{\gb}_{\gb}$, etc.
We now proceed as before: working in $W_{\gb}$
we take $(2^{\aleph_0})^+$ many hulls of
$V_{\xi}$, where $\xi$ is the 
$||\mbox{\bf s}_{\gb}||\cdot 2 +i$th cut-off
point, $i=1$ if $\gb<\gth_{n+1}$ or
$i=2$ if $\gb=\gth_{n+1}$. By collapsing we construct
$(Q^{\gb}_{\gb},q^{\gb}_{\gb},\eta^{\gb}_{\gb},
\gk^{\gb}_{\gb},\mbox{\bf s}^{\gb}_{\gb})$, so that
${\cal Q}^{\go\,\frown}
\seq{(Q^{\gb}_{\gb},q^{\gb}_{\gb},\eta^{\gb}_{\gb},
\gk^{\gb}_{\gb},\mbox{\bf s}^{\gb}_{\gb})}\in{\cal R}^{\gb}$ 
is the desired sequence extending $\bar{\cal Q}$.
This concludes the proof of the lemma.
\end{proof}

We are now ready to define ${\cal P}^{n+1}$ 
and ${}^{\ast}\gF^{n+1}$. Fix a sequence
$$
\seq{(Q^{\gth_{n+1}}_{\ga},q^{\gth_{n+1}}_{\ga},
\eta^{\gth_{n+1}}_{\ga},\gk^{\gth_{n+1}}_{\ga},
\mbox{\bf s}^{\gth_{n+1}}_{\ga})
\mid\gth_n\leq\ga\leq\gth_{n+1}}\in{\cal R}^{\gth_{n+1}}
$$ 
and set
$$
\begin{array}{rlrlrl}
P^{\ast}_{n+1}=& Q^{\gth_{n+1}}_{\gth_{n+1}}\qquad &
g^{\ast}=&q^{\gth_{n+1}}_{\gth_{n+1}}\qquad &
\gf^{\ast}_{n+1}=&g^{\ast}\circ f_{\gth_{n+1}} \\
 & & & & & \\
\eta^{\ast}_{n+1}=&\eta^{\gth_{n+1}}_{\gth_{n+1}}\qquad &
\gk^{\ast}_{n+1} =&\gk^{\gth_{n+1}}_{\gth_{n+1}} \qquad &
\br^{\ast}_{n+1}=&\mbox{\bf s}^{\gth_{n+1}}_{\gth_{n+1}}
\end{array}
$$
and for $\gth_n\leq\gb<\gth_{n+1}$
$$
\begin{array}{rlrlrl}
P_{\gb}=&Q^{\gth_{n+1}}_{\gb}\qquad &
g_{\gb}=&q^{\gth_{n+1}}_{\gb}\qquad &
\gf_{\gb}=&g_{\gb}\circ f_{\gb} \\
& & & & & \\
\eta_{\gb}=&\eta^{\gth_{n+1}}_{\gb} &
\gk_{\gb}=&\gk^{\gth_{n+1}}_{\gb} &
\mbox{\bf s}_{\gb}=&\mbox{\bf s}^{\gth_{n+1}}_{\gb}.
\end{array}
$$
The verification of the conditions for 
$({\cal P}^{n+1},{}^{\ast}\gF^n)$
is straightforward. As an example let us check (9).
As $({\cal W},{\cal P}^n)$ is internal,
$$
P_{n+1}^{\ast}\in W_{\gth_{n+1}}=M_{\gth_{n+1}}^{\cal W}
\subseteq W_{\gth_n}=P_n^{\ast}
$$
so $P^{\ast}_{n+1}\in P^{\ast}_n$.

\bigskip\noindent
$\bullet$ {\bf Construction of $p_{n+2},{\cal B}^{n+1},
\Pi^{n+1},\Pi^{n,n+1}$ and ${}^{\ast}\Psi^{n+1}$.}

\noindent
The construction is very similar to what we did before.
Pick $p_{n+2}\supset p_{n+1}$ according 
to $\gS$ such that $M(p_{n+1})\notin I_{p_{n+2}}$ and
$p_{n+2}$ can be copied on $({\cal P}^{n+1},{}^{\ast}\gF^{n+1})$.
Such a position must exist by Corollary \ref{corsuccessor}.
Thus clause (3) holds by fiat.

\noindent{\bf Case 1: $M(p_{n+1})=\min I_{p_{n+1}}$.}

\noindent Then set ${\cal B}^{n+1}={\cal B}^n$, 
$\Pi^{n,n+1}=$ the identity, $\Pi^{n+1}=\Pi^n$ and
let ${}^{\ast}\Psi^{n+1}=\seq{g_{\gb}\mid\gb<\gth_{n+1}}^{\frown}
\seq{g^{\ast}}$. 

\medskip

\noindent{\bf Case 2: $M(p_{n+1})>\min I_{p_{n+1}}$.}

\noindent For every integer $k\in I_{p_{n+1}}$, 
$k<M(p_{n+1})$ and for every
position $q$ according to $\gS$ such that 
$\mbox{\bf c}(q\re\gb)=k$ for some $\gb\leq lh(q)$,
there are witnesses to the fact that such $q$ cannot be copied
on $({\cal P}^{n+1},{}^{\ast}\gF^{n+1})$. 
For any such $q$, there is an ordinal $\gb$ 
such that $q\re\gb$ can be copied on $({\cal P}^{n+1},{}^{\ast}\gF^{n+1})$,
but $q\re\gb+1$ cannot. 
Let $\cal S$ be the pseudo-iteration tree of height $\gb$
obtained by copying $q\re\gb$ on $({\cal P}^{n+1},{}^{\ast}\gF^{n+1})$.

If $\gb$ is a limit ordinal, then fix an increasing sequence
$\gb_m\to \gb$ and ordinals $\xi_m\in M^{\cal S}_{\gb_m}$
such that $i_{\gb_m,\gb_{m+1}}(\xi_m)>\xi_{m+1}$,
and let $\mbox{w}_m=(\xi_m,\gb_m)$.
If $\gb=\nu+1$ and $\gg=<_T\mbox{-pred}(\gb)$,
then let $\mbox{w}_m=(a_m,f_m)$ witness the
illfoundedness  of $\ult(M_{\gg}^{\cal S},E_{\nu})$.

By absoluteness $\seq{\mbox{w}_m\mid m<\go}$
can be taken to be inside $P^{\ast}_{n+1}$, and
by Corollary \ref{truncation}, we can assume that 
$\seq{\mbox{w}_m\mid m<\go}\in 
P^{\ast}_{n+1}\cap V_{\gk^{\ast}_{n+1}}$.
Repeating the argument for every position $q$ as above, a set 
$X\subseteq P^{\ast}_{n+1}\cap V_{\gk^{\ast}_{n+1}}$ of all such
$\mbox{w}_m(q)$ is obtained.
$|X|\leq 2^{\aleph_0}$, as there are at most $2^{\aleph_0}$
such $q$'s.

Working inside $P^{\ast}_{n+1}$, let 
$$
H=\hull^{V_{\gk_{n+1}^{\ast}}}\big(
X\cup HC\cup\mbox{ran}(g^{\ast})\cup
\{{\cal C}^{n+1},\mbox{H}^{n+1},{}^{\ast}\Phi^{n+1},
\seq{g_{\gb}\mid\gb<\gth_{n+1}}\}
\big).
$$
By construction $|H|=2^{\aleph_0}$, and let 
$h:H\to B^{n+1}_{\gth_{n+1}}$ be the transitive collapse. 
For $\gth_n<\ga\leq\gth_{n+1}$, let
$$ 
\begin{array}{rlrlrlrl}
B_{\ga}=&h(C_{\ga})\qquad &
\ge_{\ga}=&h(\eta_{\ga})\qquad &
\ge_{\gth_{n+1}}=&h(\eta_{n+1}^{\ast})\qquad &
\psi^{n+1}_{\ga}=&h^{-1}\re B_{\ga} \\
 & & & & & & & \\
\pi^{n,n+1}_{\ga}=&h(g_{\ga}) &
\pi^{n,n+1}_{\gth_{n+1}}=&h\circ g^{\ast} &
\pi_{\ga}=&h(\gf_{\ga}) &
\pi_{\gth_{n+1}}=&h(\gf^{\ast}_{n+1}).
\end{array}
$$
This completes the construction of
${\cal B}^{n+1}$, $\Pi^{n+1}$, $\Pi^{n,n+1}$ and 
${}^{\ast}\Psi^{n+1}=
\seq{\psi^{n+1}_{\ga}\mid\ga<\gth_{n+1}}^{\frown}\seq{h^{-1}}$, 
hence Theorem \ref{main} is proved.